\documentclass[11pt]{article}

\usepackage{epsfig}
    \usepackage{amstext}
    \usepackage{amssymb}
    \usepackage{latexsym}
    \usepackage{dsfont}

\def\be{\begin{equation}}
\def\ee{\end{equation}}
\def\bea{\begin{eqnarray}}
\def\eea{\end{eqnarray}}
\def\bes{\begin{eqnarray*}}
\def\ees{\end{eqnarray*}}

\def\nn{\nonumber}
\def\<{\langle}
\def\>{\rangle}
\def\lb{\label}

\def\R{\mathbb{ R}}
\def\C{\mathbb{ C}}
\def\Z{\mathbb{ Z}}
\def\N{\mathbb{ N}}

\def\Q{\mathbb{ Q}}
\def\P{{ \mathbb{P}}}

\def\sg{{\sigma}}

\def\Sg{{\Sigma}}

\def\im{{\rm im}}

\def\hb{\vrule height0.18cm width0.14cm $\,$}

\def\mapright#1{\smash{\mathop{\longrightarrow}\limits^{#1}}}

\title{Reduced genus-two Gromov-Witten Invariants \\for complex manifolds}

\author{Wei Wang\thanks{Partially supported by National Natural
Science Foundation of China No.10801002, China Postdoctoral Science Foundation
No.20070420264 and LMAM in Peking University.
E-mail: alexanderweiwang@yahoo.com.cn, wangwei@math.pku.edu.cn  }\\
School of Mathematical Science \\ Peking University, Beijing 100871 \\
PEOPLES REPUBLIC OF CHINA \\ }

\date{March 5th, 2009}

\begin{document}

\maketitle

\begin{abstract}
{\it In this article, we construct the reduced genus-two Gromov-Witten
invariants for certain almost K\"{a}hler manifold $(X, \omega, J)$
such that $J$ is integrable and satisfies some regularity conditions.
In particular, the standard projective space $(\P^n, \omega_0, J_0)$ of dimension $n \le 7$
satisfies these conditions. This invariant counts the number of
simple genus-two $J$-holomorphic curves that satisfy appropriate number of constraints.  }
\end{abstract}

{\bf Key words}:  Reduced Gromov-Witten invariant,
pseudocycle, orbifold, obstruction, gluing.

{\bf AMS Subject Classification}: 53D45, 53D30, 37K20.

{\bf Running head}: Reduced genus-two Gromov-Witten Invariants

\renewcommand{\theequation}{\thesection.\arabic{equation}}
\renewcommand{\thefigure}{\thesection.\arabic{figure}}

\setcounter{equation}{0}
\section{Introduction and main results}

This article is devoted to a study on the reduced genus-two Gromov-Witten
invariants.

In symplectic topology and algebraic geometry,
Gromov-Witten invariants are rational numbers that, in certain situations,
count pseudo holomorphic curves satisfying prescribed conditions in a given
symplectic manifold. These invariants have been used to
distinguish symplectic manifolds.
They also play a crucial role in enumerative algebraic geometry
and were inspired by the closed type IIA string theory.
An early form of the invariants was used by Gromov in \cite{G} (also see \cite{M1})
to obtain important results on symplectic manifolds.
The genus zero Gromov-Witten invariants for semi-positive
symplectic manifolds were first studied by Ruan in distinguishing
symplectic manifolds in \cite{R1} and \cite{R2}. The first general
construction of the Gromov-Witten invariants were constructed
by Ruan and Tian in \cite{RT1} and \cite{RT2} for semi-positive
symplectic manifolds. They constructed these invariants by using
solutions of the inhomogeneous Cauchy-Riemann equation.
The invariants with fixed marked points also arose in the context of sigma models
and were considered by Witten in \cite{W} in early 90's. Later, in \cite{KM},
Kontsevich and Manin formulated the Gromov-Witten invariants in a more algebraic setting.
In 1995/1996, using the  new technique of virtual cycle constructions,
the Gromov-Witten invariants were constructed for general
algebraic manifolds first by Li-Tian \cite{LT1} and then
for general symplectic manifolds by Fukaya-Ono \cite{FO},
Li-Tian \cite{LT2}, Siebert \cite{Si} and Ruan \cite{R3}.

Let $(X,\, \omega)$ be a compact symplectic manifold of
dimension $2n$. The Gromov-Witten invariants are given as
homomorphisms
$$GW_{g, k, A}^X: H^\ast(X, \Q)^{\otimes k}\rightarrow
H^\ast(\overline{\mathcal{M}}_{g, k}; \Q).$$
To construct them, let $J$ be an almost complex structure on $X$ which
is tamed by $\omega$. In this case $(X, \omega, J)$ is called
a almost K\"{a}hler manifold.
For $A\in H_2(X, \Z)$ and a pair $(g, k)$ of
nonnegative integers, denote by $\overline{\mathfrak{M}}_{g, k}(X,
A; J)$ the moduli space of equivalence classes of stable
$J$-holomorphic maps from nodal genus-$g$ Riemann surfaces with $k$ marked
points in the homology class $A$. This determines a rational virtual fundamental class of dimension
\be\dim \overline{\mathfrak{M}}^{vir}_{g,k}(X, A)\equiv\dim^{vir}\overline{\mathfrak{M}}_{g, k}(X, A)
=2\langle c_1(TX), A\rangle+2(n-3)(1-g)+2k.\lb{1.1}\ee
Then one can pull back cohomology classes on $X$ and integrate them against the
virtual fundamental class to get the invariants.

Denote by $\mathfrak{M}^0_{g, k}(X, A; J)$ the subspace of
$\overline{\mathfrak{M}}_{g, k}(X, A; J)$ consisting of
stable maps $[\mathcal{C},\, u]$ whose domain $\mathcal{C}$
is a smooth Riemann surface of genus $g$. If $(\P^n,\omega_0, J_0) $ is the
complex projective space of complex dimension $n$ with the standard
K\"{a}hler structure and $\ell$ is the homology class of a
complex line in $\P^n$, the space
$\mathfrak{M}^0_{g, k}(\P^n,\, d)\equiv \mathfrak{M}^0_{g, k}(\P^n, d\ell; J_0)$
is a smooth orbifold of dimension $\dim \overline{\mathfrak{M}}^{vir}_{g,k}(\P^n, d\ell)$
at least for $d\ge 2g-1$. Moreover,
$\overline{\mathfrak{M}}_{0, k}(\P^n,\, d)\equiv
\overline{\mathfrak{M}}_{0, k}(\P^n,\, d\ell,\, J_0)$ is
a compact topological orbifold stratified by smooth orbifolds of
even dimensions and  $\mathfrak{M}^0_{0, k}(\P^n,\, d)$ is its main
stratum. In particular, $\mathfrak{M}^0_{0, k}(\P^n,\, d)$ is
a dense open subset of $\overline{\mathfrak{M}}_{0, k}(\P^n,\, d)$.

When $g\ge 1$, the moduli space
$\overline{\mathfrak{M}}_{g, k}(\P^n,\, d)$ has many irreducible
components of various dimensions. In particular,
$\mathfrak{M}^0_{g, k}(\P^n,\, d)$ is not dense in
$\overline{\mathfrak{M}}_{g, k}(\P^n,\, d)$. In fact,
some components of $\overline{\mathfrak{M}}_{g, k}(\P^n,\, d)$
have dimensions strictly larger than $\dim \overline{\mathfrak{M}}^{vir}_{g,k}(\P^n, d\ell)$.
Thus in general we can not use the space
$\overline{\mathfrak{M}}_{g, k}(X, A; J)$ directly to define
Gromov-Witten  invariants.

In \cite{Z3}, A. Zinger constructed the  reduced genus-one Gromov-Witten invariants for
$(X, \omega, J)$ under some regularity conditions which are satisfied for
the standard $(\P^n,\omega_0, J_0)$. In fact, he proved that the closure
$\overline{\mathfrak{M}}^0_{1, k}(X, A; J)$
of the subspace $\mathfrak{M}^0_{1, k}(X, A; J)$
in $\overline{\mathfrak{M}}_{1, k}(X, A; J)$
has  the structure of a compact topological orbifold stratified by smooth orbifolds of
even dimensions and  $\mathfrak{M}^0_{1, k}(X, A; J)$ is the  main
stratum of $\overline{\mathfrak{M}}^0_{1, k}(X, A; J)$. Thus one can use the space $\overline{\mathfrak{M}}^0_{1, k}(X, A; J)$
to define the  reduced genus-one Gromov-Witten invariants.

For the higher genus case, the space $\overline{\mathfrak{M}}_{g, k}(X, A; J)$
has many more irreducible components with dimensions strictly larger than
$\dim \overline{\mathfrak{M}}^{vir}_{g,k}(X, A)$. Thus in order to define the
reduced genus-g Gromov-Witten type invariants, one need to
construct a suitable subspace of $\overline{\mathfrak{M}}_{g, k}(X, A; J)$
containing $\mathfrak{M}^0_{g, k}(X, A; J)$ which has
better properties, eg. it represents a pseudocycle of the desired dimension.
Moreover, the newly defined invariants should have precisely geometric
meaning, eg. it counts the number of simple genus-g pseudo holomorphic curves
that pass appropriate number of constraints.

It is well-known from algebraic geometry that
$\mathfrak{M}^0_{2, k}(\P^n,\, d)$ is smooth of
expected dimension for $d\ge 3$ and its complement in
$\overline{\mathfrak{M}}^0_{2, k}(\P^n,\, d)$ is certainly of
complex codimension at least one. While for general symplectic
manifolds, this may not be true. In order to derive a
structure analogous to the moduli space
$\overline{\mathfrak{M}}^0_{2, k}(\P^n,\, d)$,
we introduce the following regularity conditions:

If $u:\mathcal{C}\rightarrow X$ is a smooth map from
a Riemann surface and $A\in H_2(X, \Z)$, we write
$$u\le_\omega A\qquad{\rm if}\qquad u_\ast[\mathcal{C}]=A
\quad{\rm or}\quad \langle\omega, u_\ast[\mathcal{C}]\rangle<\langle\omega, A\rangle.$$

{\bf Definition 1.1.} {\it Suppose $(X, \omega, J)$ is
a compact almost K\"{a}hler manifold such that
$J$ is integrable and $A\in H_2(X, \Z)$.
Then complex structure $J$ is
{\bf A-regular} if the following hold:

(i) For every $J$-holomorphic map $u:\P^1\rightarrow X$ such that
$u\le_\omega A$, we have
$H^1(\P^1, u^\ast TX)=0$ and  $H^1(\P^1, u^\ast TX\otimes \mathcal{O}_{\P^1}(-1))=0$.

(ii) For every nonconstant $J$-holomorphic map $u:\P^1\rightarrow X$ such that
$u\le_\omega A$, we have

\quad (ii-a) $H^1(\P^1, u^\ast TX\otimes \mathcal{O}_{\P^1}(-2))=0$;

\quad (ii-b) if there exists $z\in\P^1$ such that $du(z)=0$,
then $H^1(\P^1, u^\ast TX\otimes \mathcal{O}_{\P^1}(-3))=0$;

\quad (ii-c) if there exist $z_1, z_2\in\P^1$ such that
$z_1\neq z_2$ and $u(z_1)=u(z_2)$, then
$H^1(\P^1, u^\ast TX\otimes \mathcal{O}_{\P^1}(-3))=0$;

\quad (ii-d) if there exists $z\in\P^1$ such that $du(z)=0$,
then either $H^1(\P^1, u^\ast TX\otimes \mathcal{O}_{\P^1}(-4))=0$
or $u$ factor through a branched covering
$\widetilde{u}: S^2\rightarrow X$, i.e., there exists a holomorphic
branched covering $\phi: S^2\rightarrow S^2$ such that
$u=\widetilde{u}\circ\phi$ and $\deg(\phi)\ge 2$;

\quad (ii-e) if there exist $z_1, z_2\in\P^1$ such that
$z_1\neq z_2$ and $u(z_1)=u(z_2)$, then either $H^1(\P^1, u^\ast TX\otimes \mathcal{O}_{\P^1}(-4))=0$
or $u$ factor through a branched covering
$\widetilde{u}: S^2\rightarrow X$ as (ii-d).

(iii) For every nonconstant $J$-holomorphic map $u:T^2\rightarrow X$ such that
$u\le_\omega A$, we have

\quad (iii-a) $H^1(T^2, u^\ast TX)=0$ and  $H^1(T^2, u^\ast TX\otimes \mathcal{O}_{T^2}(-1))=0$;

\quad (iii-b) either $H^1(T^2, u^\ast TX\otimes \mathcal{O}_{T^2}(-2))=0$
or $u$ factor through a branched covering
$\widetilde{u}: \Sg\rightarrow X$, i.e., there exists a holomorphic
branched covering $\phi: T^2\rightarrow \Sg$ such that
$u=\widetilde{u}\circ\phi$ and $\deg(\phi)\ge 2$.

(iv) For every nonconstant $J$-holomorphic map $u:\Sigma\rightarrow X$ from a smooth Riemann surface
of genus two such that $u\le_\omega A$, either
$H^1(\Sg, u^\ast TX)=0$ or $u$ factor through a branched covering
$\widetilde{u}: \Sg^\prime\rightarrow X$, i.e., there exists a holomorphic
branched covering $\phi: \Sg\rightarrow \Sg^\prime$ such that
$u=\widetilde{u}\circ\phi$ and $\deg(\phi)\ge 2$.

(v) For every pair of nonconstant $J$-holomorphic maps $u_1, u_2:\P^1\rightarrow X$ such that
$u_1, u_2\le_\omega A$, we have:

\quad (v-a) for all $z_1, z_2\in\P^1$,
satisfying $u_1(z_1)=u_2(z_2)$ and $du_2(z_2)=\lambda du_1(z_1)$ for
some $\lambda\in\C\setminus\{0\}$, one of the following holds:
$H^1(\P^1, u_1^\ast TX\otimes \mathcal{O}_{\P^1}(-3))=0$,
$H^1(\P^1, u_2^\ast TX\otimes \mathcal{O}_{\P^1}(-3))=0$,
$u_1(\P^1)=u_2(\P^1)$.

\quad (v-b) for all $z_1, z_1^\prime, z_2, z_2^\prime\in\P^1$
and $z_1\neq z_1^\prime, z_2\neq z_2^\prime$
satisfying $u_1(z_1)=u_2(z_2)$, $u_1(z_1^\prime)=u_2(z_2^\prime)$,
one of the following holds:
$H^1(\P^1, u_1^\ast TX\otimes \mathcal{O}_{\P^1}(-3))=0$,
$H^1(\P^1, u_2^\ast TX\otimes \mathcal{O}_{\P^1}(-3))=0$,
$u_1(\P^1)=u_2(\P^1)$.
}

It is well-known that the standard complex structure $J_0$ on $\P^n$
satisfies the above regularity conditions.

For any almost K\"{a}hler manifold $(X, \omega, J)$, let
$N=\inf\{\langle c_1(TX),\,u_\ast[\Sg]\rangle
:\overline{\partial}_Ju=0,\,u_\ast[\Sg]\neq 0\}$,
where $\Sg$ is a smooth Riemann surface of genus
at most one.

The following are main results in this paper:

{\bf Theorem 1.2.} {\it Suppose $(X, \omega, J)$ is
a compact almost K\"{a}hler manifold such that $J$ is integrable,
$A\in H_2(X, \Z)\setminus\{0\}$ and the complex
structure $J$ is  A-regular in the sense of Definition 1.1.
Denote the closure of the space
$\mathfrak{M}^0_{2, k}(X, A, J)$ in $\overline{\mathfrak{M}}_{2,
k}(X, A, J)$ under the stable map topology by
$\overline{\mathfrak{M}}^0_{2, k}(X, A, J)$.
Then one of the following must hold:

(i) An element $b\equiv[\mathcal{C},\, u]\in
\overline{\mathfrak{M}}^0_{2, k}(X, A, J)\setminus\mathfrak{M}^0_{2, k}(X, A, J)$
if and only if $\{D^{(1)}_{h_i}b,\ldots,D_{h_i} ^{(k)}b\}_{h_i\in\Omega}$
satisfy a set of linear equations, where $D_{h_i} ^{(l)}b$
denotes the covariant derivative of $u_{h_i}$ of order $l$ at the corresponding
nodal point of $\mathcal{C}$ with respect to some metric on $\Sg_P$
and $X$, $\Omega$ is a set consisting of certain irreducible components
of $\mathcal{C}$, and $k$ is at most $3$.
In particular,
$\overline{\mathfrak{M}}^0_{2, k}(X, A, J)\setminus\mathfrak{M}^0_{2, k}(X, A, J)$
is a smooth orbifold of dimension at most $\dim
\overline{\mathfrak{M}}^{vir}_{2,k}(X, A)-2$ in a
neighborhood of  $[\mathcal{C}, u]$.

(ii) An element $b\equiv[\mathcal{C},\, u]\in
\overline{\mathfrak{M}}^0_{2, k}(X, A, J)\setminus\mathfrak{M}^0_{2, k}(X, A, J)$
must satisfy: there exists $b_1\equiv[\mathcal{C}_1,\, u_1]\in
\overline{\mathfrak{M}}_{g, k}(X, A_1, J)$
and $b_2\equiv[\mathcal{C}_2,\, u_2]\in
{\mathfrak{M}}^0_{g, 0}(X, A_2, J)$ such that
$u(\mathcal{C})=u_1(\mathcal{C}_1)$, $g\le 1$,
and $A-A_1=mA_2\neq0$, where $m$ is a positive integer.

}

{\bf Theorem 1.3.} {\it Under the assumption of Theorem 1.2,
suppose $\dim X\equiv 2n<\min\{N+6, 2N+6\}$.  Then the evaluation map
$$ev: \overline{\mathfrak{M}}^0_{2, k}(X, A, J)\rightarrow X^k$$
represents a pseudocycle of dimension $\dim
\overline{\mathfrak{M}}^{vir}_{2,k}(X, A)$, which can be used
to define the reduced genus-two Gromov-Witten invariants
$\mathrm{GW}_{2, k}^{0, X}(A; \cdot)$ for $(X,\omega, J)$.}

{\bf Theorem 1.4.} {\it Under the assumption of Theorem 1.3,
suppose $(\mu_1, \ldots, \mu_k)$ is a $k$-tuple of proper
submanifolds of $X$ of total codimension
$\dim \overline{\mathfrak{M}}^{vir}_{2,k}(X, A)$ in general position.
Then the invariant $\mathrm{GW}_{2, k}^{0, X}(A; (\mu_1, \ldots, \mu_k))$
counts the signed number of simple (cf. \S2.5 of \cite{MS}) genus-two
$J$-holomorphic curves with smooth domains that pass
$(\mu_1, \ldots, \mu_k)$. }

{\bf Remark 1.5.} The proof of Theorem 1.2 is based on understanding the conditions under which
a stable map $[\mathcal{C},\, u]$ lies in $\overline{\mathfrak{M}}^0_{2, k}(X, A, J)$.
The precise description of these conditions are contained in \S5 and \S6.
In the genus one case, the conditions were found in \cite{Z3}.
We believe that the reduced genus-two Gromov-Witten invariants are closely
related to the standard genus-two Gromov-Witten nvariants.
We are going to study their relation in a separate paper.

In this paper, let $\N$, $\N_0$, $\Z$, $\Q$, $\R$, and $\C$ denote
the sets of natural integers, non-negative integers, integers, rational
numbers, real numbers, and complex numbers respectively.

\noindent {\bf Acknowledgements.} I would like to sincerely thank Professor Gang Tian
for introducing me to Symplectic geometry and for his valuable
helps to me in all ways. I would like to say that how enjoyable
it is to work with him.

\setcounter{equation}{0}
\section{ Structure of the moduli space $\overline{\mathfrak{M}}_{2, k}(X, A, J)$   }

In this section, we study the structure of the moduli space
$\overline{\mathfrak{M}}_{2, k}(X, A, J)$ and the obstruction bundle
on it. In the following of this paper, we assume $(X, A, J)$
satisfies the regularity conditions in Definition 1.1.

An element $[\mathcal{C},\, u]$ in $\overline{\mathfrak{M}}_{2,
k}(X, A, J)$ is the equivalence class of a pair consisting of a
connected $k$-pointed nodal genus-two Riemann surface $\mathcal{C}$
and a $J$-holomorphic map $u: \mathcal{C}\rightarrow X$ such that
every contracted genus-$0$ component of $\mathcal{C}$ contains at
least three special points (i.e. node-branches and marked points)
and every contracted genus-$1$ component contains at least one
special point, cf. Chapter 24 of \cite{MirSym}. In general, one can
use the associated graph $T_\mathcal{C}$ of $\mathcal{C}$ to
describe $[\mathcal{C},\, u]$ as in Chapter 2 of \cite{FO}.

In order to study the structure of the moduli space, we make the
following definition:

{\bf Definition 2.1.} {\it The principle component $\Sg_P$ of a
stable map $[\mathcal{C},\, u]$ in $\overline{\mathfrak{M}}_{2,
k}(X, A, J)$ is the union of irreducible components
$\{\Sigma_i\}_{1\le i \le l}$ of  $\mathcal{C}$ such that
$\bigcup_{1\le i\le l}\Sg_i$ is a connected nodal surface of genus
two  and $l$ is the least number satisfying this property. In other
words, $\Sg_P$ is the smallest connected nodal surface in
$\mathcal{C}$ of genus-two. }

{\bf Remark 2.2.} Note that by Definition 2.1, $\mathcal{C}$ is
obtained from $\Sg_P$ by attaching bubble trees, (cf. \S3 below).
It is easy to see that $\Sg_P$ belongs to one of the following cases:

(i) A smooth Riemann surface of genus two.

(ii) Two smooth tori and a set of spheres.

(iii) A  torus with only one node.

(iv) A smooth torus and a set of spheres, they together contains exactly one circle.

(v) A set of spheres contain exactly two circles.

We illustrate each case in the following propositions separately.

In the following, we denote by $n_{nod}$ the number of nodes in
$\mathcal{C}$ and $\mathfrak{M}_T(X, A, J)$ the stratum of
$\overline{\mathfrak{M}}_{2, k}(X, A, J)$  of type $T$, where $T$ is
the combinatorial type of $[\mathcal{C}, u]$ as in Chapter 2 of
\cite{FO}.

By (iv) in Definition 1.1, we have the following:

{\bf Proposition 2.3.} {\it The muduli space $\mathfrak{M}^{0, simp}_{2, k}(X, A, J)$
consisting of simple genus-two $J$-holomorphic curves with smooth domains
is a smooth orbifold with the desired dimension
$\dim \overline{\mathfrak{M}}^{vir}_{2,k}(X, A)$.\hfill\hb }

In the following, we separate our study into several cases according
the behavior of a stable map $[\mathcal{C},\, u]$ in
$\overline{\mathfrak{M}}_{2, k}(X, A, J)$ restricted to its
principle component $\Sg_P$.

{\bf Proposition 2.4.} {\it Suppose the principle component $\Sg_P$
of a stable map $[\mathcal{C},\, u]$ in $\overline{\mathfrak{M}}_{2,
k}(X, A, J)$ is described in (i) of Remark 2.2, i.e., $\Sg_P$ is a
smooth Riemann surface of genus two. Then we have the following:

(i)  If $u|_{\Sg_P}\neq const$ and $H^1(\Sg, u^\ast TX)=0$,
then $\mathfrak{M}_T(X, A, J)$ is a
$\dim \overline{\mathfrak{M}}^{vir}_{2,k}(X, A)-2n_{nod}$
dimensional smooth orbifold.

(ii) If $u|_{\Sg_P}\neq const$ and $H^1(\Sg, u^\ast TX)\neq0$,
then $u|_{\Sg_P}$ factor through a branched covering
$\widetilde{u}: \Sg^\prime\rightarrow X$, i.e., there exists a holomorphic
branched covering $\phi: \Sg_P\rightarrow \Sg^\prime$ such that
$u|_{\Sg_P}=\widetilde{u}\circ\phi$ and $\deg(\phi)\ge 2$.

(iii)  If $u|_{\Sg_P}=const$, then $\mathfrak{M}_T(X, A, J)$ is a
$\dim \overline{\mathfrak{M}}^{vir}_{2,k}(X, A)+2(2n-n_{nod})$ dimensional smooth orbifold. }

\begin{figure}
\resizebox{9cm}{6.9cm}{\includegraphics*[0cm,0cm][18.06cm,13.55cm]{figure1.eps}}
\resizebox{9cm}{6.9cm}{\includegraphics*[0cm,0cm][18.06cm,13.55cm]{figure2.eps}}
\caption{Domains in Propositions 2.4 and 2.5}
\end{figure}

{\bf Proof.} The proposition follows from (i) and (iv) in Definition 1.1
and the implicit function theorem.
For the reader's convenience, here we give the proof of (iii). We prove the simplest case,
the general case follows similarly. Suppose there are $m$ bubbles $\{\mathcal{C}_i\}_{1\le i\le m}$
attached directly to $\Sg_P$ and $u_i\equiv u|_{\mathcal{C}_i}$
is non-constant for $1\le i\le m$. Then there is a natural isomorphism
\bea &&\mathfrak{M}_T(X, A, J)\lb{2.1}\\
\cong&&\left(\mathcal{M}_{2, k_0+m}\times
\left\{\prod_{i=1}^m\mathfrak{M}^0_{0, k_i+1}(X, A_i, J):
ev_{k_i+1}(u_i)=ev_{k_j+1}(u_j),\;1\le i, j\le m\right\}\right)\left/\frac{}{}\right.S_m,
\nn\eea
where $\mathcal{M}_{g, l}$ denotes the moduli space of
smooth Riemann surfaces of genus $g$ with $l$ marked points.
$ k_0$ denotes the number of marked points  on
$\Sg_P$ and $ k_i$ denotes the number of marked points
on $\mathcal{C}_i$ for $1\le i\le m$. In particular, we have
$\sum_{i=0}^m k_i=k$. $A_i\in H_2(X,\Z)$ and $\sum_{i=1}^m A_i=A$.
$ev_{k_i+1}(u_i)$ is the evaluation map of $u_i$ at the $(k_i+1)$-th
marked point. $S_m$ is the permutation group of order $m$.
By (i) of Definition 1.1, the evaluation map
$$ev_{k_1+1}\times\cdots\times ev_{k_m+1}: \prod_{i=1}^m\mathfrak{M}^0_{0, k_i+1}(X, A_i, J)\rightarrow X^m$$
is transversal to the diagonal $\Delta\equiv\{(x,\ldots,x)\in X^m\}$.
Hence the right hand side of (\ref{2.1}) is a smooth orbifold
by the implicit function theorem. By the index theorem, we have
\bea \dim\mathfrak{M}_T(X, A, J)&=&\dim\mathcal{M}_{2, k_0+m}
+\sum_{i=1}^m\dim\mathfrak{M}^0_{0, k_i+1}(X, A_i, J)-{\rm codim}\Delta\nn\\
&=& 2(3+k_0+m)+\sum_{i=1}^m 2(\langle c_1(TX), A_i\rangle+n-3+k_i+1)-2n(m-1)\nn\\
&=&2(\langle c_1(TX), A\rangle+k+3+n-m)=\dim \overline{\mathfrak{M}}^{vir}_{2,k}(X, A)+2(2n-m).\nn \eea
Hence (iii) holds in this case. The general case follows by a similar
argument: once there is one more node, the dimension decreases by $2$.
The proof of the proposition is complete.\hfill\hb

{\bf Proposition 2.5.} {\it Suppose the principle component $\Sg_P$
of a stable map $[\mathcal{C},\, u]$ in $\overline{\mathfrak{M}}_{2,
k}(X, A, J)$ is described in (ii) of Remark 2.2, i.e., it consists
of two smooth tori $T_1$ and $T_2$ and a set of spheres
$\{S_i\}_{1\le i\le l}$. Then we have the following:

(i)  If $\deg(u|_ {T_1})\neq 0$ and $\deg(u|_ {T_2})\neq 0$, then
$\mathfrak{M}_T(X, A, J)$ is a  $\dim
\overline{\mathfrak{M}}^{vir}_{2,k}(X, A)-2n_{nod}$
dimensional smooth orbifold.

(ii)   If $\deg(u|_ {T_1})= 0$ and $\deg(u|_ {T_2})\neq 0$ or
 $\deg(u|_ {T_1})\neq 0$ and $\deg(u|_ {T_2})= 0$, then $\mathfrak{M}_T(X, A, J)$
is a  $\dim
\overline{\mathfrak{M}}^{vir}_{2,k}(X, A)+2(n-n_{nod})$
dimensional smooth orbifold.

(iii)  If $\deg(u|_ {T_1})= 0$ and $\deg(u|_ {T_2})= 0$, then
$\mathfrak{M}_T(X, A, J)$ is a $\dim
\overline{\mathfrak{M}}^{vir}_{2,k}(X, A)+2(2n-n_{nod})$
dimensional smooth orbifold. }

{\bf Proof.} The proposition follows from (i) and (iii-a) in Definition 1.1
and the implicit function theorem.\hfill\hb

{\bf Proposition 2.6.} {\it Suppose the principle component $\Sg_P$ of a
stable map $[\mathcal{C},\, u]$ in $\overline{\mathfrak{M}}_{2, k}(X, A, J)$
is described in (iii) of Remark 2.2, i.e., a torus with only one node.
Let $(T^2, z_1, z_2)$ be the normalization of $\Sg_P$.
Then we have the following:

(i)  If $u|_{\Sg_P}\neq const$ and
$H^1(T^2, u^\ast TX\otimes \mathcal{O}_{T^2}(-z_1-z_2))=0$,
then $\mathfrak{M}_T(X, A, J)$ is a
$\dim \overline{\mathfrak{M}}^{vir}_{2,k}(X, A)-2n_{nod}$
dimensional smooth orbifold .

(ii) If $u|_{\Sg_P}\neq const$ and
$H^1(T^2, u^\ast TX\otimes \mathcal{O}_{T^2}(-z_1-z_2))\neq0$,
then $u|_{\Sg_P}$ factor through a branched covering
$\widetilde{u}: \Sg^\prime\rightarrow X$, i.e., there exists a holomorphic
branched covering $\phi: T^2\rightarrow \Sg^\prime$ such that
$u|_{\Sg_P}=\widetilde{u}\circ\phi$ and $\deg(\phi)\ge 2$.

(iii)  If $u|_{\Sg_P}=const$, then $\mathfrak{M}_T(X, A, J)$ is a
$\dim \overline{\mathfrak{M}}^{vir}_{2,k}(X, A)+2(2n-n_{nod})$
dimensional smooth orbifold. }

\begin{figure}
\resizebox{9cm}{6.9cm}{\includegraphics*[0cm,0cm][18.06cm,13.55cm]{figure3.eps}}
\resizebox{9cm}{6.9cm}{\includegraphics*[0cm,0cm][18.06cm,13.55cm]{figure4.eps}}
\caption{Domains in Propositions 2.6 and 2.7}
\end{figure}

{\bf Proof.} The proposition follows from (i) and (iii-b) in Definition 1.1
and the implicit function theorem.\hfill\hb

{\bf Proposition 2.7.} {\it Suppose the principle component $\Sg_P$
of a stable map $[\mathcal{C},\, u]$ in $\overline{\mathfrak{M}}_{2,
k}(X, A, J)$ is described in (iv) of Remark 2.2, i.e., a smooth
torus $\Sg\equiv S_0$ together with a set of spheres $\{S_i\}_{1\le
i\le l}$ and $\Omega=\{S_{i_0},\ldots, S_{i_s}\}\subset
\{S_0,\ldots, S_l\}$ form a circle. Let $A_1=u_\ast[S_0]$ and
$A_2=\sum_{i\in\Omega\setminus S_0}u_\ast[S_i]$. Then we have the
following:

(i)  If $A_1\neq 0$ and $A_2\neq 0$, then $\mathfrak{M}_T(X, A, J)$
is a  $\dim \overline{\mathfrak{M}}^{vir}_{2,k}(X, A)-2n_{nod}$ dimensional smooth orbifold.

(ii)  If $A_1= 0$ and $A_2\neq 0$,, then $\mathfrak{M}_T(X, A, J)$
is a  $\dim \overline{\mathfrak{M}}^{vir}_{2,k}(X, A)+2(n-n_{nod})$ dimensional smooth orbifold.

(iii)  If $A_1\neq 0$,  $A_2= 0$ and $S_0\notin\Omega$, then
$\mathfrak{M}_T(X, A, J)$ is a  $\dim
\overline{\mathfrak{M}}^{vir}_{2,k}(X, A)+2(n-n_{nod})$
dimensional smooth orbifold.

(iv)  If $A_1\neq 0$,  $A_2= 0$ and $S_0\in\Omega$, let
$\{z_1, z_2\}= S_0\cap\overline{\Sg_P\setminus S_0}$,
then we have:

\qquad(iv-a)  If $H^1(S_0, u^\ast TX\otimes \mathcal{O}_{S_0}(-z_1-z_2))=0$,
then $\mathfrak{M}_T(X, A, J)$ is a
$\dim \overline{\mathfrak{M}}^{vir}_{2,k}(X, A)-2n_{nod}$
dimensional smooth orbifold .

\qquad(iv-b)  If $H^1(S_0, u^\ast TX\otimes \mathcal{O}_{S_0}(-z_1-z_2))\neq0$,
then $u|_{S_0}$ factor through a branched covering
$\widetilde{u}: \Sg^\prime\rightarrow X$, i.e., there exists a holomorphic
branched covering $\phi: S_0\rightarrow \Sg^\prime$ such that
$u|_{S_0}=\widetilde{u}\circ\phi$ and $\deg(\phi)\ge 2$.

(v)  If $A_1=0=A_2$, then $\mathfrak{M}_T(X, A, J)$ is a $\dim
\overline{\mathfrak{M}}^{vir}_{2,k}(X, A)+2(2n-n_{nod})$
dimensional smooth orbifold. }

{\bf Proof.} The proposition follows from (i), (ii-a),
(ii-c), (iii-a) and (iii-b) in Definition 1.1
and the implicit function theorem. The simplest case that the domain $\mathcal{C}$
is a smooth torus $\Sg\equiv S_0$ and only one sphere $S_1$
with $S_0\in\Omega$ or $S_0\notin\Omega$ are illustrated
in Figure 2.3.\hfill\hb

\begin{figure}
\resizebox{9cm}{6.9cm}{\includegraphics*[0cm,0cm][18.06cm,13.55cm]{figure5.eps}}
\resizebox{9cm}{6.9cm}{\includegraphics*[0cm,0cm][18.06cm,13.55cm]{figure6.eps}}
\caption{Domains for $S_0\in\Omega$ and $S_0\notin\Omega$ in Proposition 2.7}
\end{figure}

{\bf Proposition 2.8.} {\it Suppose the principle component $\Sg_P$ of a
stable map $[\mathcal{C},\, u]$ in $\overline{\mathfrak{M}}_{2, k}(X, A, J)$
is described in (v) of Remark 2.2, i.e., a set of
spheres $\{S_i\}_{1\le i\le l}$ and
$\Omega_1=\{S_{i_1},\ldots, S_{i_s}\}$ and
$\Omega_2=\{S_{j_1},\ldots, S_{j_t}\}$ form two circles.
Then we have the following:

{\bf Case 1.}  If $\Omega_1\subset\Omega_2$, we let $A_1=\sum_{i\in\Omega_1}u_\ast[S_i]$ and
$A_2=\sum_{i\in\Omega_2\setminus \Omega_1}u_\ast[S_i]$.
Then we have

(i)  If $A_1\neq 0$ and $A_2\neq 0$, then $\mathfrak{M}_T(X, A, J)$
is a  $\dim \overline{\mathfrak{M}}^{vir}_{2,k}(X, A)-2n_{nod}$ dimensional smooth orbifold.

(ii)  If $A_1=0$,  $A_2\neq 0$, then $\mathfrak{M}_T(X, A, J)$
is a  $\dim \overline{\mathfrak{M}}^{vir}_{2,k}(X, A)+2(n-n_{nod})$
dimensional smooth orbifold.

(iii)  If $A_1\neq 0$,  $A_2= 0$,
then one of the following holds:

\qquad(iii-a)  $\mathfrak{M}_T(X, A, J)$ is a
$\dim \overline{\mathfrak{M}}^{vir}_{2,k}(X, A)-2n_{nod}$
dimensional smooth orbifold .

\qquad(iii-b) There exists $i_0\in \Omega_1$ such that $A_1=u_\ast[s_{i_0}]$
and $u|_{S_{i_0}}$ factor through a branched covering
$\widetilde{u}: S^2\rightarrow X$, i.e., there exists a holomorphic
branched covering $\phi: S_{i_0}\rightarrow S^2$ such that
$u|_{S_{i_0}}=\widetilde{u}\circ\phi$ and $\deg(\phi)\ge 2$.

\qquad(iii-c) There exist $i_1, i_2\in \Omega_1$ such that
$A_1=u_\ast[S_{i_1}]+u_\ast[S_{i_2}]$ and $u_{i_1}(S_{i_1})=u_{i_2}(S_{i_2})$.

(iv)  If $A_1=0=A_2$, then $\mathfrak{M}_T(X, A, J)$
is a $\dim \overline{\mathfrak{M}}^{vir}_{2,k}(X, A)+2(2n-n_{nod})$ dimensional smooth orbifold.

{\bf Case 2.}  If $\Omega_1\setminus\Omega_2\neq\emptyset$
and $\Omega_2\setminus\Omega_1\neq\emptyset$ hold,
Let $A_1=\sum_{i\in\Omega_1}u_\ast[S_i]$ and $A_2=\sum_{i\in\Omega_2\setminus \Omega_1}u_\ast[S_i]$.
Then we have the following:

(v)  If $A_1\neq 0$ and $A_2\neq 0$, then $\mathfrak{M}_T(X, A, J)$
is a  $\dim \overline{\mathfrak{M}}^{vir}_{2,k}(X, A)-2n_{nod}$ dimensional smooth orbifold.

(vi)  If $A_1\neq 0$,  $A_2= 0$ and  $\sum_{i\in\Omega_2}u_\ast[S_i]=0$,
then $\mathfrak{M}_T(X, A, J)$
is a  $\dim \overline{\mathfrak{M}}^{vir}_{2,k}(X, A)+2(n-n_{nod})$
dimensional smooth orbifold.

(vii)  If $A_1\neq 0$,  $A_2= 0$ and $\sum_{i\in\Omega_2}u_\ast[S_i]\neq0$,

then one of the following holds:

\qquad(vii-a)  $\mathfrak{M}_T(X, A, J)$ is a
$\dim \overline{\mathfrak{M}}^{vir}_{2,k}(X, A)-2n_{nod}$
dimensional smooth orbifold .

\qquad(vii-b) There exists $i_0\in \Omega_1$ such that
$A_1=u_\ast[s_{i_0}]$ and $u|_{S_{i_0}}$ factor through a branched covering
$\widetilde{u}: S^2\rightarrow X$, i.e., there exists a holomorphic
branched covering $\phi: S_{i_0}\rightarrow S^2$ such that
$u|_{S_{i_0}}=\widetilde{u}\circ\phi$ and $\deg(\phi)\ge 2$.

\qquad(vii-c) There exist $i_1, i_2\in \Omega_1$ such that
$A_1=u_\ast[S_{i_1}]+u_\ast[S_{i_2}]$ and $u_{i_1}(S_{i_1})=u_{i_2}(S_{i_2})$.

(viii)  If $A_1=0=A_2$, then $\mathfrak{M}_T(X, A, J)$
is a $\dim \overline{\mathfrak{M}}^{vir}_{2,k}(X, A)+2(2n-n_{nod})$ dimensional smooth orbifold.
}

\begin{figure}
\resizebox{9cm}{6.9cm}{\includegraphics*[0cm,0cm][18.06cm,13.55cm]{figure7.eps}}
\resizebox{9cm}{6.9cm}{\includegraphics*[0cm,0cm][18.06cm,13.55cm]{figure8.eps}}
\caption{Domains in Cases 1 and 2 of Proposition 2.8}
\end{figure}

{\bf Proof.} The proposition follows from (i), (ii-a), (ii-d) and (v-b) in Definition 1.1
and the implicit function theorem. The simplest case that the domain $\mathcal{C}$
are exactly two sphere $S_1$ and $S_2$ are illustrated
in Figure 2.5.\hfill\hb

\begin{figure}
\resizebox{9cm}{6.9cm}{\includegraphics*[0cm,0cm][18.06cm,13.55cm]{figure9.eps}}
\resizebox{9cm}{6.9cm}{\includegraphics*[0cm,0cm][18.06cm,13.55cm]{figure10.eps}}
\caption{Domains for (i) and (v) with exactly two spheres in Proposition 2.8}
\end{figure}

Summing up the results above, we have the following:

{\bf Theorem 2.9.} {\it Suppose the regularity conditions in Definition 1.1 hold.
Then we have the following:

(i) Suppose $T$ belongs to (i) of Propositions 2.4-2.8 or (iv-a) of
Proposition 2.7 or (iii-a), (v), (vii-a)  of Proposition 2.8. then $\mathfrak{M}_T(X, A, J)$ is a
$\dim \overline{\mathfrak{M}}^{vir}_{2,k}(X, A)-2n_{nod}$
dimensional smooth orbifold .

(ii) Suppose $T$ belongs to (ii) of Propositions 2.4, 2.6 or (iv-b) of
Proposition 2.7 or (iii-b), (iii-c), (vii-b), (vii-c) of Proposition 2.8.
Denote by $A_P=u_\ast[\Sg_P]$. Then we have:
there exists $b_1\equiv[\mathcal{C}_1,\, u_1]\in
\overline{\mathfrak{M}}_{g, k}(X, A_1, J)$
and $b_2\equiv[\mathcal{C}_2,\, u_2]\in
{\mathfrak{M}}_{g, 0}^0(X, A_2, J)$ such that
$u(\Sg_P)=u_1(\mathcal{C}_1)$, $g\le 1$,
and $A_P-A_1=mA_2\neq0$.

(iii) If $T$ belongs to (ii) of Propositions 2.5 or (ii), (iii) of
Proposition 2.7 or (ii), (vi) of Proposition 2.8, then $\mathfrak{M}_T(X, A, J)$ is a
$\dim \overline{\mathfrak{M}}^{vir}_{2,k}(X, A)-2(n-n_{nod})$
dimensional smooth orbifold.

(iv) If $T$ belongs to (iii) of Propositions 2.4-2.6
or (v) of Proposition 2.7 or (iv), (viii) of
Proposition 2.8, then $\mathfrak{M}_T(X, A, J)$ is a
$\dim \overline{\mathfrak{M}}^{vir}_{2,k}(X, A)-2(2n-n_{nod})$
dimensional smooth orbifold.
}

Note that for the special case $(\P^n, \omega_0, J_0)$, we have

{\bf Theorem 2.10.} {\it Suppose $d\ge 3$, then each stratum
$\mathfrak{M}_T(\P^n, d)$ of $\overline{\mathfrak{M}}_{2, k}(\P^n,\,
d)$ is a smooth orbifold. More precisely, we have the following:

(i) If $T$ belongs to (i) of Theorem 2.9,
then $\dim\mathfrak{M}_T(\P^n,\, d)=\dim \overline{\mathfrak{M}}^{vir}_{2,k}(\P^n, d\ell)-2n_{nod}$.

(ii) If $T$ belongs to (ii) of Theorem 2.9,
then $\dim\mathfrak{M}_T(\P^n,\, d)=\dim \overline{\mathfrak{M}}^{vir}_{2,k}(\P^n, d\ell)+2(n-1-n_{nod})$.

(iii) If $T$ belongs to (iii) of Theorem 2.9,
then $\dim\mathfrak{M}_T(\P^n,\, d)=\dim \overline{\mathfrak{M}}^{vir}_{2,k}(\P^n, d\ell)+2(n-n_{nod})$.

(iv) If $T$ belongs to (iv) of Theorem 2.9,
then $\dim\mathfrak{M}_T(\P^n,\, d)=\dim \overline{\mathfrak{M}}^{vir}_{2,k}(\P^n, d\ell)+2(2n-n_{nod})$.
}

{\bf Proof.} It remains to prove (ii).
We only prove the case that $\Sg_P$ is smooth, the proof can be
generalized to the case that $\Sg_P$ is not smooth.
It is sufficient to show that the muduli space $\mathfrak{M}^0_{2, k}(\P^n, 2)$ is a smooth orbifold of
dimension $4n+8+2k=\dim \overline{\mathfrak{M}}^{vir}_{2,k}(\P^n, 2\ell)+2(n-1)$.

 By Castelnuovo's bound (cf. P116 of \cite{ACGH}),
the image of each
$[\mathcal{C},\, u]\in \mathfrak{M}^0_{2, k}(\P^n, 2)$ in $\P^n$
has genus zero. Hence it must factor through a degree-one map
$\widetilde{u}: S^2\rightarrow \P^n$, i.e., there exists a holomorphic
branched covering $\phi: \Sg\rightarrow S^2$ such that
\be u=\widetilde{u}\circ\phi,\qquad \deg(\phi)=2.\lb{2.2}\ee
Thus there is a natural identification of
\be\mathfrak{M}^0_{2, k}(\P^n, 2)\cong
\mathfrak{M}^0_{0, 0}(\P^n, 1)\times \mathfrak{M}^0_{2, k}(\P^1, 2).\lb{2.3}\ee
In fact, the first factor describes the position of $\im(u)$ in $\P^n$
and the second factor describes the branched covering from $\Sg$ to $S^2$.

By the index theorem, we have
\bea \dim\mathfrak{M}^0_{2, k}(\P^n, 2)&=&
\dim\mathfrak{M}^0_{0, 0}(\P^n, 1)+\dim\mathfrak{M}^0_{2, k}(\P^1, 2)
\nn\\&=& 4n-4+12+2k=4n+8+2k.\nn\eea
Note that in the last equality, we have used the fact that
$$H^1(\Sg, \mathcal{O}(\phi^\ast T\P^1))\cong
H^0(\Sg, \mathcal{O}(\phi^\ast T\P^1)^\ast\otimes K_\Sg)^\ast\cong
H^0(\Sg, \mathcal{O}(-2))^\ast=0.$$
for any $\phi\in\mathfrak{M}^0_{2, k}(\P^1, 2)$.
Thus the linearization $D_\phi$ of the
$\overline{\partial}$-operator for the bundle $\phi^\ast T\P^1$
$$D_\phi: \Gamma(\Sg, \phi^\ast T\P^1)\rightarrow \Gamma(\Sg, \Lambda^{0, 1}T^\ast\Sg \otimes \phi^\ast  T\P^1)$$
is surjective. Hence $\mathfrak{M}^0_{2, k}(\P^1, 2)$ is a smooth orbifold of
dimension $12$ by the index theorem and the implicit function theorem.
 Clearly the smooth case has a rather simple proof:
every degree $2$ map from a smooth genus $2$ curve to $\P^n$
is a double cover of a line; such a double cover is determined by $6$
points on the line, and the dimension of $G(2, n+1)$ is $4(n-1)$.
\hfill\hb

Note that the linearized operator $D_u$ of the $\overline{\partial}$-operator
at $[\mathcal{C}, u]\in\overline{\mathfrak{M}}_{2, k}(X, A, J)$
is not surjective in general. Hence we need to study the obstructions
$H^1_{\overline{\partial}}(\mathcal{C}, u^\ast TX)$. We have the following:

{\bf Theorem 2.11.} {\it Suppose the regularity conditions in Definition 1.1 hold.
then the obstruction at $[\mathcal{C},\, u]\in\mathfrak{M}_T(X, A, J)$
is one of the following cases:

(i) If $T$ belongs to (i) of Theorem 2.9, then $H^1_{\overline{\partial}}(\mathcal{C}, u^\ast TX)=0$.

(ii) If $T$ belongs to (ii) of Theorem 2.9, then
$H^1_{\overline{\partial}}(\mathcal{C}, u^\ast TX)=\mathrm{coker} D_u$.

(iii) If $T$ belongs to (iii) of Theorem 2.9, then
$H^1_{\overline{\partial}}(\mathcal{C}, u^\ast TX)\cong \C^{n}$.

(iv) If $T$ belongs to (iv) of Theorem 2.9, then
$H^1_{\overline{\partial}}(\mathcal{C}, u^\ast TX)\cong \C^{2n}$.
}

{\bf Proof.} It follows directly from the regularity conditions
in Definition 1.1. The proof for the special case $(\P^n, \omega_0, J_0)$
is given below. The proof of the general case is the same.\hfill\hb

Note that for the special case $(\P^n, \omega_0, J_0)$, we have

{ \bf Lemma 2.12.} (cf. Corollary 6.5 of \cite{Z1})
{ \it Let $\Sg$ be a smooth Riemann surface. If $u: \Sg\rightarrow \P^n$
is a holomorphic map, then the linearization $D_u$ of the
$\overline{\partial}$-operator for the bundle $u^\ast T\P^n$
$$D_u: \Gamma(\Sg, u^\ast T\P^n)\rightarrow \Gamma(\Sg, \Lambda^{0, 1}T^\ast\Sg \otimes u^\ast  T\P^n)$$
is surjective provided $d+\chi(\Sg)>0$, where $d$ is the degree
of $u$.}

{ \bf Lemma 2.13.} {\it Let $\Sg$ be a smooth Riemann surface. If $u:
\Sg\rightarrow \P^n$ is a holomorphic map of degree $d$, then for
any tuple of pairwise distinct points $\{p_0,\ldots,
p_l\}\in\Sg^{l+1}$,  the map
$$\varphi^{(l)}: \ker D_u\rightarrow\bigoplus_{0\le m\le l} T_{u(p_m)}\P^n,
\qquad\varphi^{(l)}(\xi)=(\xi(p_0), \xi(p_1),\ldots,\xi(p_l))$$
is surjective provided $d+\chi(\Sg)\ge l+2$. }

{\bf Theorem 2.14.} {\it Suppose $d\ge 3$, then the obstruction at
$[\mathcal{C},\, u]\in\mathfrak{M}_T(\P^n,\, d)$
is one of the following cases:

(i) If $T$ belongs to (i) of Theorem 2.10, then $H^1_{\overline{\partial}}(\mathcal{C}, u^\ast T\P^n)=0$.

(ii) If $T$ belongs to (ii) of Theorem 2.10, then $H^1_{\overline{\partial}}(\mathcal{C}, u^\ast T\P^n)\cong \C^{n-1}$.

(iii) If $T$ belongs to (iii) of Theorem 2.10, then $H^1_{\overline{\partial}}(\mathcal{C}, u^\ast T\P^n)\cong \C^{n}$.

(iv) If $T$ belongs to (iv) of Theorem 2.10, then $H^1_{\overline{\partial}}(\mathcal{C}, u^\ast T\P^n)\cong \C^{2n}$.
}

{\bf Proof.} Note that any $[\mathcal{C},\, u]\in\overline{\mathfrak{M}}_{2, k}(\P^n,\, d)$
is obtained from an element in $\mathfrak{M}_{T_P}(\P^n, \deg(u|_{\Sg_P}))$
by attaching bubble trees. While there are no obstructions for
attaching bubble trees. Thus the obstructions comes from
$u_{\Sg_P}\equiv u|_{\Sg_P}$.

By the proof of Propositions 2.4-2.8, the operator
$D_u$ is surjective when $T$ belongs to (i) of Theorem 2.9.
Hence (i) holds.

We prove (ii) as follows. By Propositions 2.4-2.8,
when $T$ belongs to (ii) of Theorem 2.9, we must have
$\deg(u|_{\Sg_P} )=2$ and we can write
\be\mathfrak{M}_{T_P}(\P^n, 2)
\cong \mathfrak{M}^0_{0, 0}(\P^n, 1)\times \mathfrak{M}_{T_P}(\P^1, 2).\lb{2.4}
\ee
First we consider the case that $\Sg_P$ is smooth.
In this case, any $[\mathcal{C}, u]$ in $\mathfrak{M}_{T_P}(\P^n, 2)$
 must factor through a degree-one map
$\widetilde{u}: S^2\rightarrow \P^n$, i.e., there exists a holomorphic
branched covering $\phi: \Sg_P\rightarrow S^2$ such that
\be u=\widetilde{u}\circ\phi,\qquad \deg(\phi)=2.\lb{2.5}\ee
Since $\im(\widetilde{u})$ is a line in $\P^n$, it is a complex submanifold
of $\P^n$. We denote it by $\P^1_{\widetilde{u}}$.
Hence we have
\be u_{\Sg_P}^\ast T\P^n = u_{\Sg_P}^\ast (T\P^1_{\widetilde{u}}\oplus N_{\P^1_{\widetilde{u}}, \P^n})
= u_{\Sg_P}^\ast (T\P^1_{\widetilde{u}}\oplus \oplus_{i=1}^{n-1}H_{\P^1_{\widetilde{u}}, i}),
\lb{2.6}\ee
where we denote by $N_{\P^1_{\widetilde{u}}, \P^n}$ the normal bundle of
$\P^1_{\widetilde{u}}$ in $\P^n$ and $H_{\P^1_{\widetilde{u}}, i}$ its
decomposition into $n-1$ line bundles.
Thus by Dolbeault Theorem, we have
\bea H^1_{\overline{\partial}}(\Sg_P, u_{\Sg_P}^\ast T\P^n)&\cong&
H^1(\Sg_P, \mathcal{O}(u_{\Sg_P}^\ast T\P^n))\nn\\&\cong&
H^1(\Sg_P, \mathcal{O}(u_{\Sg_P}^\ast (T\P^1_{\widetilde{u}}\oplus \oplus_{i=1}^{n-1}H_{\P^1_{\widetilde{u}}, i}))\nn\\
&\cong& H^1(\Sg_P, \mathcal{O}(u_{\Sg_P}^\ast (T\P^1_{\widetilde{u}})))\bigoplus \oplus_{i=1}^{n-1}H^1(\Sg_P, \mathcal{O}(u_{\Sg_P}^\ast(H_{\P^1_{\widetilde{u}}, i})))\nn\\
&\cong& H^0(\Sg_P, \mathcal{O}(-2))^\ast\bigoplus \oplus_{i=1}^{n-1}H^0(\Sg_P, \mathcal{O})^\ast\cong \C^{n-1}
\lb{2.7}\eea
by Kodaira-Serre duality. The proof of the other cases are similar.
For the reader's convenience, here we give the proof
of case (ii) in Proposition 2.6 and omit the proofs of the others.
Suppose $\Sg_P$ is a torus $\Sg$ with only one node
and $(\Sg, x_1, x_2)$ is the  normalization of $\Sg_P$.
Let $\xi\in L^p(\Sg_P, \Lambda^{0, 1}T^\ast\Sg_P \otimes u_{\Sg_P}^\ast T\P^n)$
and we want to find $\sigma\in L_1^p(\Sg_P, u_{\Sg_P}^\ast T\P^n)$
such that $D_{u_{\Sg_P}}\sigma=\xi$. Since
$\xi\in L^p(\Sg, \Lambda^{0, 1}T^\ast\Sg \otimes u_{\Sg}^\ast T\P^n)$,
we can find $\widehat{\sigma}\in L_1^p(\Sg, u_{\Sg}^\ast T\P^n)$
such that $D_{u_\Sg}\hat\sigma=\xi$ by Lemma 2.12.
Since $\deg(u_\Sg)=2$, we may assume $\widehat{\sigma}(x_1)=0$ by Lemma 2.13. Hence if $\sigma$ exists, we must have
$\widehat{\sigma}(x_2)=0$ also.
We consider the short exact sequence of sheaves on $\Sg$
\be 0\rightarrow\mathcal{O}(u_{\Sg}^\ast T\P^n\otimes(-x_1-x_2))
\rightarrow\mathcal{O}(u_{\Sg}^\ast T\P^n\otimes(-x_1))
\mapright{\gamma}\mathcal{O}((u_{\Sg}^\ast T\P^n\otimes(-x_1))_{x_2})\rightarrow0
\lb{2.8}\ee
where we view $\mathcal{O}((u_{\Sg}^\ast T\P^n\otimes(-x_1))_{x_2})$
as a sheaf on $\Sg$ via extension by $0$ (cf. p.38 of \cite{GH}).
Taking the corresponding long exact sequence in cohomology,
we obtain
\bea
&&\mapright{} H^0(\Sg, \mathcal{O}(u_{\Sg}^\ast T\P^n\otimes(-x_1)))
\mapright{\gamma}
H^0(\Sg, \mathcal{O}((u_{\Sg}^\ast T\P^n\otimes(-x_1))_{x_2})) \nn\\
&&\mapright{\partial}
H^1(\Sg, \mathcal{O}(u_{\Sg}^\ast T\P^n\otimes(-x_1-x_2)))
\mapright{}
H^1(\Sg, \mathcal{O}(u_{\Sg}^\ast T\P^n\otimes(-x_1))).
\lb{2.9}\eea
Note that
$ H^1(\Sg, \mathcal{O}(u_{\Sg}^\ast T\P^n\otimes(-x_1)))=0$
by Lemma 2.13. Thus $\partial$ is surjective. Thus we have
\bea&&\mathrm{coker} D_{u_{\Sg_P}}\cong\mathrm{coker}\gamma
\cong H^1(\Sg, \mathcal{O}(u_{\Sg}^\ast T\P^n\otimes(-x_1-x_2)))\nn\\
\cong&& H^1(\Sg, \mathcal{O}(u_{\Sg}^\ast (T\P^1_{\widetilde{u}})\otimes(-x_1-x_2)))\bigoplus \oplus_{i=1}^{n-1}H^1(\Sg, \mathcal{O}(u_{\Sg}^\ast(H_{\P^1_{\widetilde{u}}, i})\otimes(-x_1-x_2)))\nn\\
\cong&& H^0(\Sg, \mathcal{O}(-2))^\ast\bigoplus \oplus_{i=1}^{n-1}H^0(\Sg, \mathcal{O})^\ast\cong \C^{n-1}
\lb{2.10}
\eea
The proof of (iii) is obvious and we omit it here.

We prove  (iv). Denote by $\Sg$ the union of components
of $\Sg_P$ which are mapped to constants such that each connected component
of $\Sg$ has genus greater than zero. Then $\Sg$ contains one or two
connected components and each one is mapped to a constant. In the second case,
we can write $\Sg=\Sg_1\cup\Sg_2$.
Then we have
\be H^1_{\overline{\partial}}(\mathcal{C}, u^\ast T\P^n)\cong \mathcal{H}^{0, 1}_\Sg\otimes T_{ev(\Sg)}\P^n\cong\C^{2n}
\lb{2.11}\ee
provided $\Sg$ is connected. Here we denote by $\mathcal{H}^{0, 1}_\Sg$ the space of harmonic
$(0, 1)$-forms on $\Sg$, cf. \S 22.3 of \cite{MirSym}.
\be H^1_{\overline{\partial}}(\mathcal{C}, u^\ast T\P^n)\cong (\mathcal{H}^{0, 1}_{\Sg_1}\otimes T_{ev(\Sg_1)}X)\oplus(\mathcal{H}^{0, 1}_{\Sg_2}\otimes T_{ev(\Sg_2)}X)\cong\C^{2n}
\lb{2.12}\ee
provided $\Sg$ is disconnected. Now the theorem follows.\hfill\hb

\setcounter{equation}{0}
\section{Gluing construction }

Given a stable map $[\mathcal{C},\, u]\in\mathfrak{M}_T(X, A, J)
\subset\overline{\mathfrak{M}}_{2, k}(X, A, J)$. Our goal in this
section is to construct approximately $J$-holomorphic maps $[\Sg_v,
u_v]\in \mathfrak{X}^0_{2, k}(X, A, J)$ by using the gluing
technique, where  $\mathfrak{X}^0_{2, k}(X, A, J)$ denotes the space
of equivalence classes of smooth maps from $\Sg_v$ ro $X$ with
$k$ marked points in the homology class $A$ and $\Sg_v$ is a
smooth Riemann surface of genus two depending on the gluing
parameter $v$. Roughly speaking, $\Sg_v$ is obtained from
$\mathcal{C}$ by replacing each attaching node of the corresponding
two components by thin necks connecting them. Thus geometrically
$\Sg_v$ is a smooth Riemann surface of genus two, but should be
viewed as a Riemann surface  close to $\mathcal{C}$. While $u_v$
equals to $u$ away from the thin necks and $u_v$ is close to $u$ in
an appropriate sense. Thus $u_v$ is $J$-holomorphic away from the thin
necks.

\setcounter{equation}{0}
\subsection{Gluing in bubble trees  }

In this section, we describe the gluing construction in bubble trees.
We proceed as in \cite{Z2} and \cite{Z3} in this section.
Let $q_N, q_S:\C\rightarrow S^2\subset \R^3$ be the
stereographic projections mapping the origin of
$\C$ to the north and south poles respectively. Explicitly, we have
\be q_N(z)=\left(\frac{2z}{1+|z|^2}, \,\frac{1-|z|^2}{1+|z|^2}\right),
\quad q_S(z)=\left(\frac{2z}{1+|z|^2}, \,\frac{-1+|z|^2}{1+|z|^2}\right).
\lb{3.1}\ee
We denote the south pole of $S^2$, i.e., the point $(0, 0, -1)\in\R^3$
by $\infty$ and $e_\infty=dq_S|_0(\frac{\partial}{\partial s})$, where
we write $z=s+it\in\C$. We  identify  $\C$ with $S^2\setminus\{\infty\}$ via the map
$q_N$.

{\bf Definition 3.1.} {\it  A rooted tree $I$ is a finite  partially ordered
set satisfying: if $h, h_1, h_2\in I$ such that
$h_1, h_2<h$, either $h_1\le h_2$ or $h_2\le h_1$ holds;
moreover, $I$ has a unique minimal element $\hat{0}$, i.e., $\hat{0}<h$
for all $h\in \hat{I}\equiv I\setminus\{\hat0\}$.     }

For any $h\in\hat I$, denote by $\iota_h\in I$ the largest
element of $I$ which is smaller that $h$.
We call $\iota: \hat I\rightarrow I$ the attaching map of $I$.

{\bf Definition 3.2.} {\it Suppose $M$ is a finite set.
A $X$-valued bubble tree with  $M$-marked points is a
tuple
\bea b=(M, I, x, (j, y), u),\quad {\rm and}\quad &&x:\hat{I}\rightarrow S^2\setminus\{\infty\},\;
j:M\rightarrow I,\nn\\ &&y:M\rightarrow S^2\setminus\{\infty\},\; u:I\rightarrow C^\infty(S^2, X)
\lb{3.2}\eea
such that $u_h(\infty)=u_{\iota_h}(x_h)$ for all $h\in\hat I$.
The special points on each bubble $\Sg_h\equiv\{h\}\times S^2$, i.e.,
$(j_l, y_l)\in\Sg_h$ and $(\iota_l, x_l)$
with $\iota_l=h$ together with the point $(h, \infty)$, are pairwise distinct. In addition,
if ${u_h}_\ast [S^2]=0\in H_2(X, \Z)$, then
$\Sg_h$ should contain at least three special points.
$u$ is $J$-holomorphic if its restriction to each component is.  }

We associate such a tuple with a nodal Riemann surface
\be \Sg_b=\left(\bigsqcup(\{h\}\times S^2)_{h\in I}\right)\left/\frac{}{}\right.\sim
\lb{3.3}\ee
where $(h, \infty)\sim(\iota_h, x_h)$ for $h\in \hat I$. We call $x_h$ the attaching node
of the bubble $h$. Clearly we obtain a continuous map $u_b:\Sg_b\rightarrow X$.

The general structure of bubble trees is described by tuples
$T_B=(M, I, j, \underline{A})$, where $I$ and $j$ are maps as described
in Definitions 3.1 and 3.2, while $A_h={u_h}_\ast[\Sg_h]\in H_2(X, \Z)$ for $h\in I$.
We call such tuples {\it bubble types.}
Denote by $\mathcal{H}_{T_B}$ the space of holomorphic maps
of type $T_B$ and  $\mathfrak{M}_{T_B}$ its equivalence classes.

For each $h\in I$, let
\be \chi_{T_B}h=\left\{\matrix{0,\qquad &&{\rm if}\quad A_i=0\quad\forall i\le h; \cr
1,\qquad &&{\rm if}\quad A_h\neq0,\quad{\rm and} \quad A_i=0\quad\forall i< h; \cr
2,\qquad &&{\rm otherwise.}\cr }\right.
\lb{3.4}\ee
Let $\beta: \R\rightarrow[0, 1]$ be a smooth function such that
\be \beta(t)=\left\{\matrix{0\quad {\rm if}\quad t\le 1, \cr
             1\quad {\rm if}\quad t\ge 2 \cr }\right.
\qquad{\rm and} \qquad \beta^\prime(t)>0,\quad {\rm for}\quad t\in(1, 2)\lb{3.5}\ee
and $\beta_r(t)=\beta(r^{-\frac{1}{2}}t)$ for any $r>0$.

Given a bubble type $T_B=(M, I, j, \underline{A})$, let $d(T_B):
I\rightarrow\R$ be given by \be d_i(T_B)=\langle\omega, A_i\rangle+|\{l\in
M:j_l=i\}|+\sum_{\iota_h=i}d_h(T_B), \qquad \forall i\in
I.\lb{3.6}\ee
Note that $d_i(T_B)$ is uniquely determined by
Definition 3.1 and (\ref{3.6}). A bubble tree in (\ref{3.2}) is
called {\it balanced} if for all $i\in \hat{I}$ the following conditions hold:
\bea &&(B1) \int_{\C}|du_i\circ q_N|^2z+\sum_{\iota_h=i}d_h(T_B)x_h+\sum_{j_l=i}y_l=0;\nn\\
&&(B2) \int_{\C}|du_i\circ
q_N|^2\beta(|z|)+\sum_{\iota_h=i}d_h(T_B)\beta(|x_h|)+\sum_{j_l=i}\beta(|y_l|)=\frac{1}{2}.
\nn\eea It is called {\it completely balanced} if (B1) and (B2) hold
for all $i\in I$.

Denote by $PSL(2, \C)$ the group of M\"{o}bius transformations. Let
\be PSL(2, \C)^{(0)}=\{g\in PSL(2, \C): g(\infty)=\infty\},\qquad
\mathcal{G}_{T_B}=\prod_{h\in I}PSL(2, \C)^{(0)}. \lb{3.7}\ee For
$b=(M, I, x, (j, y), u)\in\mathcal{H}_{T_B}$ and $g\in PSL(2,
\C)^{(0)}$, define $gb=(M, I, gx, (j, gy), gu)$ by \bea
(gx)_h=g_{\iota_h}x_h,\quad (gy)_l=g_{j_l}y_l,\quad (gu)_i=g_iu_i,
\lb{3.8}\eea where for a map $u:S^2\rightarrow X$ and $g\in PSL(2,
\C)$, we define $(gu)(z)=u(g^{-1}z)$.

Let $\mathcal{M}^{(0)}_{T_B}\subset\mathcal{H}_{T_B}$ denote the
subset of completely balanced bubble trees. Then the group
$\prod_{h\in I}S^1\times  \mathrm{Aut}(T_B)$ acts on
$\mathcal{M}^{(0)}_{T_B}$ and all the stabilizers are finite. Then
we have
\be\mathfrak{M}_{T_B}\cong\mathcal{M}^{(0)}_{T_B}\left/\frac{}{}\right.\left(\prod_{h\in
I}S^1\times  \mathrm{Aut}(T_B)\right). \lb{3.9}\ee
By Proposition 3.3 in \cite{Z2},  $\mathcal{M}^{(0)}_{T_B}$ is a smooth oriented
manifold and then $\mathfrak{M}_{T_B}$ is a smooth orbifold. One may
think of elements in $\mathcal{M}^{(0)}_{T_B}$ as good
representatives of $\mathfrak{M}_{T_B}$. In fact,
$\mathcal{M}^{(0)}_{T_B}=\Psi_{T_B}^{-1}((0, \frac{1}{2})^I)$, where
$\Psi_{T_B}\equiv(\Psi_{T_B, i})_{i\in
I}:\mathcal{H}_{T_B}\rightarrow (\C\times\R)^I$ is defined by \bea
\Psi_{T_B, i}(M, I, x, (j, y), u)=&&\left( \int_{\C}|du_i\circ
q_N|^2z+\sum_{\iota_h=i}d_h(T_B)x_h+\sum_{j_l=i}y_l,\right.
\nn\\&&\left.\int_{\C}|du_i\circ
q_N|^2\beta(|z|)+\sum_{\iota_h=i}d_h(T_B)\beta(|x_h|)+\sum_{j_l=i}\beta(|y_l|)\right)
\lb{3.10}\eea and $\Psi_{T_B}$ is smooth and transversal to every
point $(0, r_i)_{i\in I}$ such that $|r_i-\frac{1}{2}|\le
\frac{1}{4}$ for all $i\in I$. Let \bea
\widetilde{\mathcal{M}}^{(0)}_{T_B}&=&\Psi_{T_B}^{-1}\left(\{(0,
r_i)_{i\in I}: r_i\in\left(\frac{1}{4}, \frac{3}{4}\right)\;\;{\rm
if }\;\;\chi_{T_B}i=1, \quad r_i=\frac{1}{2} \;\;{\rm
otherwise}\}\right),
\lb{3.11}\\
\widetilde{\mathcal{FT}}_B&=&\widetilde{\mathcal{M}}^{(0)}_{T_B}\times
\C^{|\hat{I}|}\rightarrow\widetilde{\mathcal{M}}^{(0)}_{T_B}, \qquad
\widetilde{\mathcal{FT}}_B^\emptyset=\widetilde{\mathcal{M}}^{(0)}_{T_B}\times
\{\C\setminus\{0\}\}^{|\hat{I}|}.\lb{3.12} \eea If
$\pi_{\mathfrak{F}}: \mathfrak{F}\rightarrow\mathfrak{X}$ is a
normed  vector bundle and $\delta:\mathfrak{X}\rightarrow\R$ is any
function, let $\mathfrak{F}_{\delta}=\{v\in\mathfrak{F}:
|v|<\delta(\pi_{\mathfrak{F}}(v))\}$ be the $\delta$-disk bundle.

Now we describe the basic gluing construction in bubble trees. For
each sufficiently small element $v=(b,
v)\in\widetilde{\mathcal{FT}}_B^\emptyset$, where $(\Sg_b, u_b)$ is
an element of $\mathcal{M}^{(0)}_{T_B}$, let $q_v:
\Sg_v\rightarrow\Sg_b$ be the basic gluing map constructed in
\cite{Z2}.  Let \be b(v)=(\Sg_v, u_v), \qquad u_v=u_b\circ q_v
\lb{3.13}\ee be the approximately $J$-holomorphic map corresponding to
$v$. The primary marked point $y_0(v)$ of $\Sg_v$ is the point
$\infty$ of $\Sg_v\cong S^2$. By the construction of $q_v$, it
factors through each of the maps $q_{v, i}: \Sg_v\rightarrow\Sg_b$
for $i\in I$. Let $g_b$ be the Riemannian metric on $\Sg_b$ such
that its restriction to each component is the standard metric on
$S^2$. By \S3.3 of \cite{Z2}, We can construct a Riemannian metric
$g_v$ on  $\Sg_v$  such that:

(G1) $q_v: (\Sg_v, g_v)\rightarrow(\Sg_b, g_b)$ is an isometry (and
thus holomorphic) outside of the annuli
\bea A^+_{v. h}&=&q_{v, \iota_h}^{-1}\left(\{z\in\Sg_{b, \iota_h}: 1\le|v_h|^{-\frac{1}{2}}|\phi_{x_h}z|\le 2 \}\right),\lb{3.14}\\
 A^-_{v. h}&=&q_{v, \iota_h}^{-1}\left(\left\{z\in\Sg_{b, \iota_h}: \frac{1}{2}\le|v_h|^{-\frac{1}{2}}|\phi_{x_h}z|\le 1 \right\}\right),\lb{3.15}\eea
where $\phi_xz=z-x\equiv q_N^{-1}(z)-q_N^{-1}(x)\in\C$ for $x, z\in
S^2\setminus\{\infty\}$.

(G2) $q_{v, \iota_h}: (A^{\pm}_{v, h}, \,g_v)\rightarrow(q_{v,
\iota_h}(A^{\pm}_{v, h}),\, g_b)$ is an isometry.

Moreover, the map $q_v$ collapses $\hat{I}$ disjoint circles on
$\Sg_v$ and is a diffeomorphism away from them. These circles are
mapped to the  $\hat{I}$ nodal branches. Alternatively, $(\Sg_v,
g_v)$ can be viewed as the surface obtained by smoothing the nodes
of $\Sg_b$. An explicit construction of $q_{v, i}$ may be described
as follows. For a  rooted tree  $I$ and a tuple $v\equiv(v_h)_{h\in
\hat{I}}\in\C^{\hat{I}}$ such that $\sum_{h\in\hat{I}}|v_h|$ is
sufficiently small, choose any ordering $\prec$ of $I$ consistent
with its partial ordering. If  $v_h\in\C$ with $0<|v_h|< \delta$,
let $p_{h, (x_h, v_h)}: B_{x_h}(\delta^{\frac{1}{2}})
\equiv\{\phi_{x_h}z<\delta^{\frac{1}{2}}\}\rightarrow \C\cup
\{\infty\}$ be given by \be p_{h, (x_h,
v_h)}(z)=(1-\beta_{|v_h|}(2|\phi_{x_h}z|))\overline{\left(\frac{v_h}{\phi_{x_h}z}\right)}
\lb{3.16}\ee and define $q_{v,\, (x_h,\, v_h)}:
\Sg_{T_B^h}\rightarrow\Sg_{T_B^h}\cup\Sg_h$ by \be q_{v,\, (x_h,\,
v_h)}(z)= \left\{\matrix{(h,\, q_S(p_{h, (x_h, v_h)}(z))),&&\quad
{\rm if}\quad  |v_h|^{-\frac{1}{2}}|\phi_{x_h}z|\le 1; \cr (\iota_h,
\,\phi_{x_h}^{-1}(\beta_{|v_h|}(|\phi_{x_h}z|)\phi_{x_h}z)),&&\quad
{\rm if}\quad 1\le|v_h|^{-\frac{1}{2}}|\phi_{x_h}z|\le 2;\cr
(\iota_h,\, \;z),&&\quad {\rm otherwise},\cr}\right.\lb{3.17}\ee
where $\Sg_{T_B^h}$ is obtained from $T_B$ by dropping the
bubble $h$ together with all bubbles descendent from it. Thus $q_{v,\,
(x_h,\, v_h)}$ is a diffeomorphism except on the circle
$|v_h|^{-\frac{1}{2}}|\phi_{x_h}z|= 1$ and  the circle is mapped to
the point $(h, \infty)=(\iota_h, x_h) $. Moreover, $q_{v,\, (x_h,\,
v_h)}$ is holomorphic outside the annulus
$\frac{1}{2}\le|v_h|^{-\frac{1}{2}}|\phi_{x_h}z|\le 2$. Taking
$q_{v,\,\hat 0}=Id$ and $q_{v,\,h}=q_{v,\, (x_h,\, v_h)}\circ
q_{v,\,\iota_h}$ inductively according to the ordering $\prec$,
we obtain $q_{v,\,h}$ for all $h\in I$.

By (G1), $u_v$ is $J$-holomorphic outside $A^{\pm}_{v. h}$. For $p>2$,
we define norms $\|\cdot\|_{v, p, 1}$ and $\|\cdot\|_{v, p}$ on
$\Gamma(v)\equiv L_1^p(\Sg_v,\, u_v^{\ast}TX)$ and  $\Gamma^{0,
1}(v)\equiv L^p(\Sg_v, \,\Lambda^{0, 1}T^\ast\Sg\otimes
u_v^{\ast}TX)$ respectively as in \S3.3 of \cite{Z2}. These norms
are equivalent to the ones used in \cite{LT2}. Let $D_v:
\Gamma(v)\rightarrow\Gamma^{0, 1}(v)$ be the linearization of the
$\overline{\partial}_J$-operator at $b(v)$. Since the linearization
$D_b$ of the $\overline{\partial}_J$-operator at $b$ is surjective by
(i) of Definition 1.1, if $v\in\widetilde{\mathcal{FT}}_B^\emptyset$ is
sufficiently small, $D_v$ is  also surjective. In particular, we
have a decomposition \bea  \Gamma(v)=\Gamma_-(v)&\oplus& \Gamma_+(v)
\equiv\{\xi\circ q_v:\xi\in\Gamma_-(b)\equiv\ker D_b\}\nn\\
&\oplus&\{\zeta\in \Gamma(v): \zeta(\hat {0}, \infty)=0;
\;\langle\zeta, \xi\rangle_{v, 2}=0,\; \forall\xi\in\Gamma_-(v)\;
{\rm s.t.}\;\xi(\hat {0}, \infty)=0\}, \lb{3.18}\eea where $(\hat 0,
\infty)\equiv y_0(v)$ is the primary marked point of $\Sg_v$. Note
that the choice of $\Gamma_+(v)$ is permissible by (i) of Definition
1.1. Moreover, the operator $D_v:\Gamma_+(v)\rightarrow \Gamma^{0,
1}(v)$ is an isomorphism and the norms of $D_v$ and of the inverse
of its restriction to $\Gamma_+(v)$  depend only on $b$ and not on
$v$. Let \be \pi_{v, -}:
\Gamma(v)\rightarrow\Gamma_-(v),\qquad\pi_{v, +}:
\Gamma(v)\rightarrow\Gamma_+(v) \lb{3.19}\ee be the projection maps
corresponding to the decomposition (\ref{3.18}).

Denote by $\mathfrak{X}_{0, M}(X, A, J)$ the space of equivalence
classes of smooth maps into $X$ from genus zero Riemann surfaces
with marked points indexed by the set $\{0\}\cup M$ in the homology
class $A$ and by $\mathfrak{X}^0_{0, M}(X, A, J)$ its subset
consisting of those maps with smooth domains. If
$K\subset\mathfrak{M}_{T_B}$ denote by $K^{(0)}$ and
$\widetilde{K}^{(0)}$ the preimages of $K$ under the projections
${\mathcal{M}}^{(0)}_{T_B}\rightarrow\mathfrak{M}_{T_B}$ and
$\widetilde{\mathcal{M}}^{(0)}_{T_B}\rightarrow\mathfrak{M}_{T_B}$
respectively. We have the following Lemma for gluing in bubble
trees.

{ \bf Lemma 3.3.} (cf. Lemma 3.3 of \cite{Z3}) {\it For every
precompact open subset $K$ of  $\mathfrak{M}_{T_B}$, there exist
$\delta_K, \epsilon_K, C_K\in\R^+$ and an open neighborhood $U_K$ of
$K$ in $\mathfrak{X}_{0, M}(X, A, J)$ such that

(i) For all $v=(b, v)\in\widetilde{\mathcal{FT}}_{B,
\delta_K}^\emptyset|_{\widetilde{K}^{(0)}}$, the equation
$$ \overline{\partial}_J\exp_{u_v}\zeta=0, \qquad\zeta\in\Gamma_+(v),\quad
\|\zeta\|_{v, p, 1}<\epsilon_K,
$$
has a unique solution $\zeta_v$.

(ii) The map $\widetilde{\phi}:\widetilde{\mathcal{FT}}_{B,
\delta_K}^\emptyset|_{\widetilde{K}^{(0)}} \rightarrow
\mathfrak{M}^0_{0, \{0\}\cup M}(X, A, J)\cap U_K,\; v\mapsto
[\exp_{b(v)}\zeta_v]$ is smooth.

(iii) For all $v=(b, v)\in\widetilde{\mathcal{FT}}_{B,
\delta_K}^\emptyset|_{\widetilde{K}^{(0)}}$, we have
$ev_0(\widetilde{\phi}(v))=ev_0(b)$.

(iv) For all $v=(b, v)\in\widetilde{\mathcal{FT}}_{B,
\delta_K}^\emptyset|_{\widetilde{K}^{(0)}}$, we have
$||\zeta_v\|_{v, p, 1},\; \|\nabla^T\zeta_v\|_{v, p, 1}\le
C_K|v|^{1/p}$, where $\nabla^T\zeta_v$ denotes the covariant
derivative with respect to the connection defined in \S3 in
\cite{Z3}.. }

\subsection{Gluing in the principle component }

The general structure of genus-two bubble maps is described as a
tuple \be T=((I_1\cup M_P, \Sg_P, A_P), (T_B^{(l)})_{l\in I_1} ),
\lb{3.20}\ee where $\Sg_P$ is a nodal Riemann surface of genus-two
as in Remark 2.2, $M_P$ denotes the  marked points on $\Sg_P$,
$A_P\in H_2(X, \Z)$ and $T_B^{(l)}$s are bubble
trees defined in \S 3.1 for $l\in I_1$. We denote by $\{x_l\}_{l\in
I_1}$ the $I_1$ points on $\Sg_P$ where the corresponding bubble
trees are attached. Let \be\mathcal{C}=\left(\Sg_P\bigsqcup
(T_B^{(l)})_{l\in I_1}\right)\left/\frac{}{}\right.\sim, \qquad
(\Sg_P,\, x_l)\sim(T_B^{(l)},\,(\hat 0, \infty)),\quad \forall l\in
I_1, \lb{3.21}\ee where $(\hat 0, \infty)$ is the primary marked
point on $\Sg_{T_B^{(l)}}$ and $\Sg_{T_B^{(l)}}$ is the nodal
Riemann surface of genus zero corresponding to the bubble tree
$T_B^{(l)}$. We call $\Sg_P$ the  principle component
of $\mathcal{C}$.

Note that we have a natural isomorphism
\bea &&\mathfrak{M}_T(X, A, J)\lb{3.22}\\
\cong &&\{(b_P, (b^{(l)})_{l\in I_1})\in\mathfrak{M}_{T_P}(X, A_P, J)
\times \prod_{l\in I_1} \mathfrak{M}_{T_B^{(l)}}:
ev_0(b^{(l)})=ev_{\iota_l}(b_P),\;\forall l\in
I_1\}\left/\frac{}{}\right. \mathrm{Aut}^{\ast}(T), \nn\eea where
$T_P=(I_1\cup M_P, \Sg_P, A_P)$,  $ev_0(b^{(l)})$ is the evaluation
map at the primary marked point $(\hat 0, \infty)$  and
$ev_{\iota_l}(b_P)$ is the evaluation map at the attaching node
$x_l$ of the bubble tree $T_B^{(l)}$ and
$\mathrm{Aut}^{\ast}(T)=\mathrm{Aut}(T)/\{g\in \mathrm{Aut}(T):
g\cdot h=h,\;\forall h\in I_1$\}. Let
\bea{\mathcal{M}}^{(0)}_{T}&=&\{(b_P, (b^{(l)})_{l\in
I_1})\in\mathfrak{M}_{T_P}(X, A_P, J)\times \prod_{l\in I_1}
{\mathcal{M}}^{(0)}_{T_B^{(l)}}:
ev_0(b^{(l)})=ev_{\iota_l}(b_P),\;\forall l\in I_1\},
\lb{3.23}\\
\widetilde{\mathcal{M}}^{(0)}_{T}&=&\{(b_P, (b^{(l)})_{l\in
I_1})\in\mathfrak{M}_{T_P}(X, A_P, J)\times \prod_{l\in I_1}
\widetilde{\mathcal{M}}^{(0)}_{T_B^{(l)}}:
ev_0(b^{(l)})=ev_{\iota_l}(b_P),\;\forall l\in I_1\}. \lb{3.24}\eea
Then clearly we have $\mathfrak{M}_T(X, A, J)
={\mathcal{M}}^{(0)}_{T}/(\mathrm{Aut^\ast}(T)\times(S^1)^{|I|})$,
where $I=\cup_{l\in I_1} I_B^{(l)}$ and $I_B^{(l)}$
is the partially ordered set associated to the bubble tree
$T_B^{(l)}$ as in Definition 3.1. Let \bea
&\widetilde{\mathcal{FT}}=\widetilde{\mathcal{F}_P\mathcal{T}}\oplus
\widetilde{\mathcal{F}_0\mathcal{T}}\oplus\widetilde{\mathcal{F}_1\mathcal{T}}
\rightarrow\widetilde{\mathcal{M}}^{(0)}_{T},\nn\\
&\widetilde{\mathcal{F}_P\mathcal{T}}=\widetilde{\pi}_P^\ast\mathcal{FT}_P,\quad
\widetilde{\mathcal{F}_0\mathcal{T}}=\bigoplus_{h\in
I_1}\widetilde{\mathcal{F}_h\mathcal{T}},\;
\widetilde{\mathcal{F}_h\mathcal{T}}=\widetilde{\pi}_P^\ast
L_h\mathcal{T}_0,\quad
\widetilde{\mathcal{F}_1\mathcal{T}}=\widetilde{\mathcal{M}}^{(0)}_{T}\times
\C^{|I\setminus I_1|}, \lb{3.25}\eea where $\widetilde{\pi}_P:
\widetilde{\mathcal{M}}^{(0)}_T\rightarrow \mathfrak{M}_{T_P}(X, A_P, J)$
is the projection map, $L_h\mathcal{T}_0$ is the universal
tangent line bundle at the marked point $x_h$ for $h\in I_1$ and
$\mathcal{FT}_P$ is the bundle of gluing parameters in the principle
component. In particular, $\mathrm{rank}
\mathcal{FT}_P=n_{nod}(\Sg_P)$, i.e., the number of nodes in
$\Sg_P$. As before, let $\widetilde{\mathcal{FT}}^{\emptyset}$ be the
subset of $\widetilde{\mathcal{FT}}$ consisting of those elements with
all components nonzero.

Now we describe the gluing construction in the principle component
$\Sg_P$. For each sufficiently small element $v=(b,
v)\in\widetilde{\mathcal{FT}}^\emptyset$, we have \be v=(b,
v)\equiv(b, v_P, v_0, \{v^{(l)}\}_{l\in
I_1})\in\widetilde{\mathcal{F}_P\mathcal{T}}\oplus
\widetilde{\mathcal{F}_0\mathcal{T}}\oplus\widetilde{\mathcal{F}_1\mathcal{T}}
\lb{3.26}\ee We smooth out all the nodes in $\Sg_P$ by the parameter
$v_P$. The bundle of gluing parameters in the principle component
$\mathcal{FT}_P$ over $\mathfrak{M}_{T_P}(X, A_P, J)$ has the form
\be \mathcal{FT}_P= \bigoplus_{x\in \mathrm{nod}(\Sg_P)} T_{x,
0}\Sg_{x, 0}\otimes T_{x, 1}\Sg_{x, 1} ,\lb{3.27}\ee where we denote
by $\Sg_{x ,0}$ and $\Sg_{x ,1}$ the two components corresponding to
the node $x$. Here in order to simplify notations, we omit the
use of a finite cover of $\mathfrak{M}_{T_P}(X, A_P, J)$ as in
\cite{RT2}. For any nonzero $v_x\in T_{x, 0}\Sg_{x, 0}\otimes
T_{x, 1}\Sg_{x, 1}$, define the map \bea\Phi_{x, v_x}: T_{x,
0}\Sg_{x, 0}\setminus\{0\}\rightarrow T_{x, 1}\Sg_{x,
1}\setminus\{0\}, \qquad X\otimes\Phi_{x, v_x}X=v_x.\lb{3.28}\eea
Now let  $\phi_{x, 0}: \Sg_{x, 0}\rightarrow T_{x, 0}\Sg_{x, 0}$ and
$\phi_{x, 1}: \Sg_{x, 1}\rightarrow T_{x, 1}\Sg_{x, 1}$ be
holomorphic coordinates near $x$ on the two components respectively.
Let $\widetilde{p}_{h, (x, v_x)}: \{\phi_{x,
0}z<\delta^{\frac{1}{2}}\}\rightarrow T_{x, 1}\Sg_{x, 1}$ be given
by \be \widetilde{p}_{h, (x, v_x)}(z)=(1-\beta_{|v_x|}(2|\phi_{x,
0}z|))\Phi_{x, v_x}(\phi_{x, 0}z) \lb{3.29}\ee and define
$$q_{v,\, (x,\, v_x)}:
(\Sg_{x, 0}\setminus\{|v_x|^{-\frac{1}{2}}|\phi_{x,
0}z|\le\frac{1}{2}\}) \cup(\Sg_{x,
1}\setminus\{|v_x|^{-\frac{1}{2}}|\phi_{x, 1}z|\le2\})
\rightarrow\Sg_{x, 0}\cup\Sg_{x, 1}$$ by \be q_{v,\, (x,\, v_x)}(z)=
\left\{\matrix{(0,\, \;z),&&\quad {\rm if}\quad z\in\Sg_{x,
0}\setminus\{|v_x|^{-\frac{1}{2}}|\phi_{x, 0}w|\le2\};\cr (0,
\,\phi_{x, 0}^{-1}(\beta_{|v_x|}(|\phi_{x, 0}z|)\phi_{x,
0}z)),&&\quad {\rm if}\quad 1\le|v_x|^{-\frac{1}{2}}|\phi_{x,
0}z|\le 2;\cr (1,\, \phi_{x, 1}^{-1}(\widetilde{p}_{h, (x,
v_x))}(z))),&&\quad {\rm if}\quad
\frac{1}{2}\le|v_x|^{-\frac{1}{2}}|\phi_{x, 0}z|\le 1; \cr (1,\,
\;z),&&\quad {\rm if}\quad z\in\Sg_{x,
1}\setminus\{|v_x|^{-\frac{1}{2}}|\phi_{x, 1}w|\le2\}.\cr
}\right.\lb{3.30}\ee Note that by (\ref{3.28}) and (\ref{3.29}), the
map $q_{v,\, (x,\, v_x)}$ is well-defined. We smooth out all the
nodes in $\Sg_P$ as above and obtain $q_{b_P, v_P}$. Then we define
\be q_{v_P}: \Sg_{(b,\, v_P)}\rightarrow\Sg_{b}\equiv\mathcal{C}
\lb{3.31}\ee to be the extension of $q_{b_P, v_P}$  by identity to
the bubble components. Geometrically, $\Sg_{(b,\, v_P)}$ is obtained
from $\Sg_b$ by replacing all the nodes in the principle component
by thin necks connecting the corresponding two components.

Let $v=(b, v)\equiv(b, v_P, v_0, \{v^{(l)}\}_{l\in I_1})$ be given
by (\ref{3.26}). If $i\in {I}_B^{(h)}$, we put \be
\rho_i(v)=\left(b,\;v_{0, h}\prod_{\{i^\prime \in\hat{I}_B^{(h)}:
\;i^\prime \le i\}}v^{(h)}_{i^\prime}\right)\in T_{x_h}\Sg_P,
\lb{3.32}\ee where the product term is defined to be $1$ if
$\{i^\prime \in\hat{I}_B^{(h)}: \;i^\prime \le i\}=\emptyset$. We
denote by \be \chi(T)=\{i: \chi_{{T_B}^{(h)}}i=1,\;h\in
I_1\}.\lb{3.33}\ee

\setcounter{equation}{0}
\section{  Study for $\mathfrak{M}_T(X, A, J)$ in
(i) of Theorem 2.9}

In the remaining of this paper, we prove the main theorem by looking for
conditions under which the approximately $J$-holomorphic map $[\Sg_v,
u_v]\in \mathfrak{X}^0_{2, k}(X, A, J)$ can be deformed
into a $J$-holomorphic map, where $\Sg_v$ is the smooth Riemann
surface of genus two and $u_v$ is the approximately $J$-holomorphic map constructed
below. We will separate the proof into several sections according to
the classification of stable maps in Theorem 2.9.

In this section we study stable maps in (i) of Theorem 2.9.
In these cases we have $H^1_{\overline{\partial}}(\mathcal{C}, u^\ast TX)=0$ by Theorem 2.10.

Let  $\Sg_{(b,\, v_P)}$ be the  Riemann surface
constructed via (\ref{3.30}) in \S 3.2.
By construction, its principle component
is a smooth Riemann surface $\Sg_{b_P, v_P}$ of genus-two and
$\Sg_{(b,\, v_P)}$ is obtained from  $\Sg_{b_P, v_P}$ by attaching
$|I_1|$ bubble trees at the points $\{x_h\}_{h\in I_1}$.
For each $h\in I_1$, we identify a small neighborhood $U(x_h)$ of
$x_h$ in $\Sg_{b_P, v_P}$ with a neighborhood $\widetilde{U}(x_h)$ of
$0$ in $T_{x_h}\Sg_{b_P, v_P}$
biholomorphicly and isometrically. In fact, we can choose a
K\"{a}hler metric  $g_{b_P, v_P}$  on $\Sg_{b_P, v_P}$
to be flat on each $U(x_h)$.
We assume that all of these neighborhoods are disjoint from each
other and from the $n_{nod}(\Sg_P)$ thin necks of $\Sg_{b_P, v_P}$.
if $z\in U(x_h)$, denote by $|z-x_h|$ its norm
with respect to the metric  $g_{b_P, v_P}$.
Then we define the map
\be q_{v_0}: \Sg_{(b, v_P, v_0)}\rightarrow\Sg_{(b,\, v_P)}
\lb{4.1}\ee
via the formula (\ref{3.17}) by replacing the term $\phi_{x_h}z$
there by $z-x_h\in T_{x_h}\Sg_{b_P, v_P}$.
Then we smooth out all the nodes in the bubble trees as in \S3.1
to obtain
\be q_{v_1}: \Sg_v\equiv\Sg_{(b, v)}\rightarrow\Sg_{(b,\, v_P, v_0)}
\lb{4.2}\ee
At last, we define
\be q_v=q_{v_P}\circ q_{v_0}\circ q_{v_1}: \Sg_{v}\rightarrow\Sg_b\equiv\mathcal{C}.
\lb{4.3}\ee
By construction, $q_v$ is a homeomorphism
outside of $n_{nod}(\Sg_b)$ circles of $\Sg_v$ and is biholomorphic
outside of $n_{nod}(\Sg_b)$ thin necks.
We take
\be b(v)=(\Sg_v, j_v, u_v), \qquad{\rm where}\qquad u_v={u}_b\circ q_v,
\lb{4.4}\ee
to be the approximately $J$-holomorphic map corresponding to the
basic gluing map $q_v$, where $j_v$ is the complex structure on $\Sg_v$.   We denote by
\be\Gamma(v)\equiv L_1^p(\Sg_v, u_v^{\ast}TX),\qquad
\Gamma^{0, 1}(v)\equiv L^p(\Sg_v, \Lambda^{0, 1}_{j_v}T^\ast\Sg_v\otimes u_v^{\ast}TX),
\lb{4.5}\ee
the Banach completions of the corresponding spaces of
smooth sections with respect to the norms
$\|\cdot\|_{v, p, 1}$ and $\|\cdot\|_{v, p}$
induced from the basic gluing map $q_v$ as in \cite{Z2}.
Let
\be \Gamma_-(v)=\{(\xi\circ q_{v}):
\xi\in\Gamma_-(b)\equiv \ker D_b\}\subset\Gamma(v)
\lb{4.6}\ee
and $\Gamma_+(v)$ the $(L^2,\,v)$-orthogonal complement of
$\Gamma_-(v)$ in  $\Gamma(v)$. Let $\pi_{v,\,\pm}$ be the
$(L^2,\,v)$-orthogonal projections onto $\Gamma_\pm(v)$ respectively.

The following is the main theorem in this section.

{ \bf Theorem 4.1.} {\it Suppose $T=((I_1\cup M_P, \Sg_P, A_P), (T_B^{(l)})_{l\in I_1} )$
is a bubble type belongs to (i) of Theorem 2.9, then for every precompact open subset $K$
of $\mathfrak{M}_T(X, A, J)$, there exist
$\delta_K, \epsilon_K, C_K\in\R^+$ and an open neighborhood
$U_K$ of $K$ in $\mathfrak{X}_{2, k}(X, A, J)$ with the following properties:

(i) For all $v=(b, v)\in\widetilde{\mathcal{FT}}_{\delta_K}^\emptyset|_{\widetilde{K}^{(0)}}$,
the equation
\be \overline{\partial}_J\exp_{u_v}\zeta=0, \qquad\zeta\in\Gamma_+(v),\quad
\|\zeta\|_{v, p, 1}<\epsilon_K,
\lb{4.7}\ee
has a unique solution $\zeta_v$.

(ii) The map $\widetilde{\phi}:\widetilde{\mathcal{FT}}_{\delta_K}^\emptyset|_{\widetilde{K}^{(0)}}
\rightarrow \mathfrak{M}^0_{2, k}(X, A, J)\cap U_K,\;
v\mapsto [\exp_{b(v)}\zeta_v]$ is smooth.

In particular, we have
$\mathfrak{M}_T(X, A, J)\subset\overline{\mathfrak{M}}^0_{2, k}(X, A, J)$
is a smooth orbifold of dimension at most
$\dim \overline{\mathfrak{M}}^{vir}_{2,k}(X, A)-2$.  }

{\bf Proof.} Note that  in these cases, the operator $D_b$ is surjective.
Hence $D_v$ is also surjective provided $v=(b, v)\in\widetilde{\mathcal{FT}}^\emptyset|_{\widetilde{K}^{(0)}}$
is sufficiently small by continuity, where $D_v: \Gamma(v)\rightarrow\Gamma^{0, 1}(v)$ is the linearization
of the $\overline{\partial}_J$-operator at $b(v)$. Moreover, by the
choice of the norms $\|\cdot\|_{v, p, 1}$ and $\|\cdot\|_{v, p}$,
we have the following estimates similar to Theorem 4.1 in \cite{Z3}
\bea &\|\pi_{v, -}\xi\|_{v, p, 1}\le C_K\|\xi\|_{v, p, 1}, \;\forall \xi\in\Gamma(v);
\qquad \|D_v\xi\|_{v, p}\le C_K|v|^{\frac{1}{p}}\|\xi\|_{v, p, 1}, \;\forall \xi\in\Gamma_-(v);
\lb{4.8}\\
&C_K^{-1}\|\xi\|_{v, p, 1}\le\|D_v\xi\|_{v, p}\le C_K\|\xi\|_{v, p, 1}, \quad\forall \xi\in\Gamma_+(v).
\lb{4.9}\eea
Thus by a standard argument as in the genus-zero case,
cf. \cite{MS} or \cite{Z3}, the theorem follows. \hfill\hb

\setcounter{equation}{0}
\section{ Study for $\mathfrak{M}_T(X, A, J)$ in (iv) of Theorem 2.9}

In this section we study stable maps in (iv) of Theorem 2.9.
In these cases we have $H^1_{\overline{\partial}}(\mathcal{C}, u^\ast TX)\cong \C^{2n}$
by Theorem 2.10. In \cite{Z1}, A. Zinger studied the enumerative
problem of genus-two curves with a fixed complex structure in $\P^n$
for $n=2, 3$. In  \cite{Z3}, A. Zinger defined the reduced genus-one Gromov-Witten
invariants. Our present paper uses  similar arguments as these two papers,
but our cases are more complicated, we need to study all the
boundary components $\mathfrak{M}_T(X, A, J)$, while in \cite{Z1},
the author only need to consider two cases for $\P^2$ and five cases for $\P^3$.

First we study the case that $\Sg_P$ is a smooth Riemann surface
of genus two in \S5.1-5.3. In this case the obstruction bundle has the form
$\mathcal{H}^{0, 1}_{\Sg_P}\otimes T_{ev(\Sg_P)}X\cong\C^{2n}$,
where $\mathcal{H}^{0, 1}_{\Sg_P}$ is the space of harmonic
$(0, 1)$-forms on $\Sg_P$ and $ev(\Sg_P)$ is the evaluation map
at the principle component.  We study the general cases in \S5.4-5.5,
which are modifications of the methods in \S5.1-5.3.

Now let us assume $\Sg_P$ is smooth. Given $\psi\in\mathcal{H}^{0, 1}_{\Sg_P}$,
$b\equiv(b_P, \,(b^{(l)})_{l\in I_1})\in \widetilde{\mathcal{M}}^{(0)}_{T}$, $x\in\Sg_P$, $m\ge 1$
and  a K\"{a}hler metric $g_{b,\, \Sg_P}$ on $\Sg_P$ which is flat
near $x$. Define $D_{b,\, x}^{(m)}\psi\in T_x^{\ast 0, 1}\Sg_P^{\otimes m}$ as follows:
If $(s, \,t)$ are conformal coordinates centered at
$x$ such that $s^2+t^2$ is the square of the $g_{b,\, \Sg_P}$-distance to $x$.
Let
\be \{D_{b,\, x}^{(m)}\psi\}\left(\frac{\partial}{\partial s}\right)
\equiv\{D_{b, \,x}^{(m)}\psi\}\underbrace{\left(\frac{\partial}{\partial s},\cdots,\frac{\partial}{\partial s}\right)}_m
=\frac{\pi}{m!}\left\{\frac{D^{m-1}}{ds^{m-1}}\psi\left|_{(s,\, t)=0}\frac{}{}\right.\right\}\left(\frac{\partial}{\partial s}\right),
\lb{5.1}\ee
where the covariant derivatives are taken with respect to the
metric $g_{b, \,\Sg_P}$. Since $\psi\in\mathcal{H}^{0, 1}_{\Sg_P}$,
we have $\psi=f(ds-idt)$ for some anti-holomorphic function $f$.
Because $g_{b, \,\Sg_P}$ is flat near $x$,
it follows that $D_{b,\, x}^{(m)}\psi\in T_x^{\ast 0, 1}\Sg_P^{\otimes m}$.
For an orthogonal basis $\{\psi_1, \,\psi_2\}$ of $\mathcal{H}^{0,\, 1}_{\Sg_P}$,
let $s^{(m)}_{b, \,x}\in T_x^{\ast}\Sg_P^{\otimes m}\otimes\mathcal{H}^{0, 1}_{\Sg_P}$
be given by
\be s^{(m)}_{b, \,x}(v)\equiv s^{(m)}_{b, \,x}\underbrace{(v,\cdots,v)}_m
=\sum_{1\le j\le 2}\left(\overline{\{D_{b,\, x}^{(m)}\psi_j\}(v)}\right)\psi_j.
\lb{5.2}\ee
The section $s^{(m)}_{b,\, x}$ is independent of the choice of
a basis for $\mathcal{H}^{0, 1}_{\Sg_P}$ but is dependent on
the choice of the metric $g_{b, \,\Sg_P}$ when $m>1$. However,
$s^{(1)}_{b,\, x}$ depends only on $(\Sg_P,\, j_{\Sg_P})$,
we denote it by $s_{\Sg_P, \,x}$. By p.246 in \cite{GH},
$s_{\Sg_P, \,x}$ does not vanish and thus spans a subbundle of
$\Sg_P\times\mathcal{H}^{0,\, 1}_{\Sg_P}\rightarrow\Sg_P$.
Denote this subbundle by $\mathcal{H}^+_{\Sg_P}$ and its orthogonal complement by
$\mathcal{H}^-_{\Sg_P}$. In particular, $\mathcal{H}^+_{\Sg_P}$
is independent of the choice of the metric $g_{b, \,\Sg_P}$ on $\Sg_P$.

For $h\in \cup_{l\in I_1}I_B^{(l)}\equiv I$ and
$m\in\mathds{N}$, define
\be D^{(m)}_hb=\frac{2}{(m-1)!}\frac{D^{m-1}}{ds^{m-1}}\frac{d}{ds}(u_h\circ q_S)\left|\frac{}{}\right._{(s, \,t)=0},\lb{5.3}\ee
where the covariant derivatives are taken with respect to
the standard metric $s+it\in\C$ and a metric $g_{X,\,b}$ on $X$
which is flat near $u(\Sg_P)$, e.g. we may assume $X$ is isomorphic to $\C^n$
near $u(\Sg_P)$.

Let $\delta_T\in C^\infty(\widetilde{\mathcal{M}}^{(0)}_{T}, \R^+)$
satisfying $4\delta_T(b)\|du_i\|_{b, C^0}<r_{X}$ for any $b\in\widetilde{\mathcal{M}}^{(0)}_{T}$,
where $r_{X}$ is the injectivity radius of $X$
with respect to the metric $g_{X,\,b}$.
We use K\"{a}hler metrics
$g_{b, \Sg_P}$ on $\Sg_P$ which are flat near  $x_h$ for
$h\in I_1$

For $h\in I$ and $\epsilon>0$, denote by
\bea \widetilde{A}^-_{b, \,h}(\epsilon)&=&
\{(h, z)\in\Sg_{b,\, h}\equiv\{h\}\times S^2: |z|>\epsilon^{-\frac{1}{2}}/2\},\lb{5.4}\\
\widetilde{A}^+_{b, h}(\epsilon)&=&
\{(\iota_h, z)\in \Sg_{b,\,\iota_h}: |z-x_h|<2\epsilon^{\frac{1}{2}}\}.\lb{5.5}\\
A^\pm_{v, h}(\epsilon)&=& q_v^{-1}(\widetilde{A}^\pm_{b, h}(\epsilon))\subset\Sg_v,\lb{5.6}
\eea
where $\iota_h=\Sg_P$ for $h\in I_1$.

Now for $v=(b, v)\in\widetilde{\mathcal{FT}}_{\delta_T}^\emptyset$
sufficiently small and $V\psi\in T_{u(\Sg_P)}X\otimes\mathcal{H}^{0,
1}_{\Sg_P}$, define $R_vV\psi\in\Gamma^{0, 1}(u_v)$ as follows: If
$z\in\Sg_v$ is such that $q_v(z)\in\Sg_{b,\, h}$ for some
$h\in\chi(T)$ as defined in (\ref{3.33}) and $|q_S^{-1}(q_v(z))|\le
2\delta_T(b)$, we define $\overline{u}_v(z)\in T_{u(\Sg_P)}X$ by
$\exp_{u(\Sg_P)}\overline{u}_v(z)=u_v(z)$. Given $z\in\Sg_v$, let
$h_z$ be such that $q_v(z)\in \Sg_{b,\, h_z}$. If $w\in T_z\Sg_v$,
put \be R_vV\psi|_zw= \left\{\matrix{0 ,&&\quad {\rm if}\quad
\chi_Th_z=2;\cr
\beta(\delta_T(b)|q_vz|)(\psi|_zw)\Pi_{\overline{u}_v(z)}V,&&\quad
{\rm if}\quad \chi_Th_z=1;\cr (\psi|_zw)V,&&\quad {\rm if}\quad
\chi_Th_z=0,\cr}\right.\lb{5.7}\ee where $\chi_T$ is the natural
extension of $\chi_{T_B^{(l)}}$ to $T$ and $\Pi_{\overline{u}_v(z)}$
is the parallel transport along the geodesic
$t\mapsto\exp_{u(\Sg_P)}t\overline{u}_v(z)$ with respect to the
Levi-Civita connection of the metric $g_{X,\,b}$.

By the same proof of Lemma 4.3 of \cite{Z1}, we have the following expansion:

{\bf Lemma 5.1.} {\it Suppose $T$ is a bubble type given by (iv) of
Theorem 2.9 with $\Sg_P$ being a smooth Riemann surface of genus
two. Then there exists $\delta\in C^\infty(\mathfrak{M}_T(X, A, J),\,\R^+)$
such that for all
$v=(b,\,v)\in\widetilde{\mathcal{FT}}_{\delta}^\emptyset$, $V\in
T_{u(\Sg_P)}X$ and $\eta\in\mathcal{H}^{0, 1}_{\Sg_P}$ we have
\bea\langle\langle \overline{\partial}_Ju_v,\,
R_vV\eta\rangle\rangle_{v,\,2}= -\sum_{m\ge1,\, h\in\chi(T)}\langle
D^{(m)}_hb, \,V\rangle
\overline{\left(\{D_{b,\,\widetilde{x}_h(v)}^{(m)}\eta\}
((d\phi_{b,\,
{\mathcal{T}(h)}}|_{\widetilde{x}_h(v)})^{-1}\rho_h(v))\right)}
\nn\eea where $\mathcal{T}(h)$ is determined by $h\in
{I}_B^{(\mathcal{T}(h))}$, $\widetilde{x}_h(v)=q_{v,\,
\iota_h}^{-1}(\iota_h,\, x_h)\in\Sg_v$ and $\phi_{b, h}$ is a
holomorphic identification of neighborhoods of $x_h$ in $\Sg_{b,
\,\iota_h}$ and $T_{x_h}\Sg_{b,\, \iota_h}$, and $\rho_h(v)$ is
given by (\ref{3.32}).  } \hfill\hb

Next we estimate the formal adjoint $D_v^\ast$ of the linearization
$D_v$ of the $\overline{\partial}_J$-operator at $u_v$ with respect to
the above $(L^2,\,v)$-inner product.

{ \bf Lemma 5.2.} (cf. Lemma 2.2 of \cite{Z1}) {\it Suppose $T$ is a
bubble type given by (iv) of Theorem 2.9 with $\Sg_P$ being a
smooth Riemann surface of genus two. Then there exists $\delta\in
C^\infty(\mathfrak{M}_T(X, A, J),\,\R^+)$ such that for all
$v=(b,\,v)\in\widetilde{\mathcal{FT}}_{\delta}^\emptyset$ and
$V\eta\in T_{u(\Sg_P)}X\otimes\mathcal{H}^{0, 1}_{\Sg_P}$, we have
$D_v^\ast R_vV\eta$ vanishes outside of the annuli \bea
\widehat{A}_{v,\,h}\equiv q_v^{-1}(\{(h,\,z)\in\Sg_{b,\,h}:
\delta_T(b)\le|q_S^{-1}(q_v(z))|\le 2\delta_T(b)\}) \nn\eea with
$h\in\chi(T)$. Moreover, there exists $C\in
C^\infty(\mathfrak{M}_T(X, A, J),\,\R^+)$ such that \bea ||D_v^\ast
R_vV\eta||_{v,\,C^0}\le C(b)
\left(\sum_{h\in\chi(T)}|\rho_h(v)|\right)|V|_v||\eta||_2.
\lb{5.8}\eea  }

{\bf Proof.} Since we will use the expression of
$D_v^\ast R_vV\eta$ below, we give the proof of the lemma here.
Let $(s,\, t)$ be the conformal
coordinates on $\widehat{A}_{v,\,h}$ given by
$q_v(s,\, t)=s+it\in\C$. Write $g_v=\theta^{-2}(s,\,t)(ds^2+dt^2)$,
then we have $\theta=\frac{1}{2}(1+s^2+t^2)$ by Riemannian geometry.
Let $\xi(s,\,t)$ be given by
\bea \xi(s,\,t)=\{R_vV\eta\}_{(s,\,t)}\partial_s=
\beta\left(\delta_T(b)\sqrt{s^2+t^2}\right)(\eta_{(s,\,t)}\partial_s)
\Pi_{\overline{u}_v(s,\,t)}V.
\nn\eea
Then by Remark C.1.4 of \cite{MS} we have
\bea D_v^\ast R_vV\eta|_z=\theta^2\left(-\frac{D}{ds}\xi+J\frac{D}{dt}\xi\right),
\nn\eea
where $\frac{D}{ds}$ and $\frac{D}{dt}$ denote covariant
derivations with respect to the metric $g_{X,\,b}$ on $X$.
Since this metric is flat on the support of $\xi$ by assumption
and $\eta\in \mathcal{H}^{0, 1}_{\Sg_P}$, we have
\bea D_v^\ast R_vV\eta|_z=\frac{(1+s^2+t^2)^2}{4}
\left(\beta^\prime|_{\delta_T(b)\sqrt{s^2+t^2}}\delta_T(b)\frac{-s+it}{\sqrt{s^2+t^2}} \right)
\cdot(\eta|_{(s,\,t)}\partial_s)\Pi_{\overline{u}_v(s,\,t)}V.
\lb{5.9}\eea
Note that the right hand side of (\ref{5.9}) vanishes
unless $\delta_T(b)\le \sqrt{s^2+t^2}\le 2\delta_T(b)$
by (\ref{3.5}). Hence we have
\bea |D_v^\ast R_vV\eta|_{v,\,z}\le C(b_v)|\eta|_{(s,\,t)}\partial_s|
|V|\le C(b)|\rho_h(v)||\eta||_2|V|.
\lb{5.10}\eea
Hence (\ref{5.8}) holds.\hfill\hb

Next we describe our choice for a tangent-space model as
\S2.3 of \cite{Z1}.
Let $\Gamma_\pm(v)$ and $\pi_{v,\,\pm}$ be given by
the formula (\ref{4.6}). Now we fix an $h^\ast\in\chi(T)$ and let
\be \overline{\Gamma}_-(v)=D_v^\ast R_v(\mathcal{H}_{\Sg_P}^+(\widetilde{x}_{h^\ast}(v))\otimes T_{u(\Sg_P)}X),
\lb{5.11}\ee
where $\mathcal{H}_{\Sg_P}^+$ is defined below (\ref{5.2}) and
$\widetilde{x}_{h^\ast}(v)=q_{v,\, \iota_{h^\ast}}^{-1}(\iota_{h^\ast},\, x_{h^\ast})\in\Sg_v$.
Denote by $\overline{\Gamma}_+(v)$ the $(L^2,\,v)$-orthogonal complement of
$\overline{\Gamma}_-(v)$ in  $\Gamma(v)$ and $\overline{\pi}_{v,\,\pm}$ the
$(L^2,\,v)$-orthogonal projections onto $\overline{\Gamma}_\pm(v)$.
Let $\widetilde{\Gamma}_+(v)$ be the image of $\Gamma_+(v)$ under
$\overline{\pi}_{v,\,+}$ and $\widetilde{\Gamma}_-(v)$ be its
$(L^2,\,v)$-orthogonal complement.
Denote by $\widetilde{\pi}_{v,\,\pm}$ the
$(L^2,\,v)$-orthogonal projections onto $\widetilde{\Gamma}_\pm(v)$.
Then we have the following:

{ \bf Lemma 5.3.} {\it Suppose $T$ is a bubble type given by (iv) of
Theorem 2.9 with $\Sg_P$ being a smooth Riemann surface of genus
two. Then for every precompact open subset $K$ of $\mathfrak{M}_T(X,
A, J)$, there exist $\delta_K, \epsilon_K, C_K\in\R^+$ and an open
neighborhood $U_K$ of $K$ in $\mathfrak{X}_{2, k}(X, A, J)$ with the
following properties:

(i) For all $v=(b,
v)\in\widetilde{\mathcal{FT}}_{\delta_K}^\emptyset|_{\widetilde{K}^{(0)}}$,
we have \bea &\|\widetilde{\pi}_{v, \pm}\xi\|_{v, p, 1}\le
C_K\|\xi\|_{v, p, 1}, \qquad\forall \xi\in\Gamma(v)
\lb{5.12}\\
&C_K^{-1}\|\xi\|_{v, p, 1}\le\|D_v\xi\|_{v, p}\le C_K\|\xi\|_{v, p,
1}, \qquad\forall \xi\in\widetilde{\Gamma}_+(v). \lb{5.13}\eea

(ii) For all $[\widetilde{b}]\in \mathfrak{X}_{2, k}(X, A, J)\cap
U_K$, there exist $v=(b,
v)\in\widetilde{\mathcal{FT}}_{\delta_K}|_{\widetilde{K}^{(0)}}$ and
$\zeta\in\widetilde{\Gamma}_+(v)$ such that $\|\zeta\|_{v, p,
1}<\epsilon_K$ and $[\exp_{b(v)}\zeta_v]=[\widetilde{b}]$.

}

{\bf Proof.} The proof of (i) follows from \S2.3 in \cite{Z1}.
In fact, $\widetilde{\Gamma}_-(v)$ is a
tangent-space model in the sense of Definition 3.11 in
\cite{Z2}. Then (\ref{5.12}) and (\ref{5.13}) follow
from Lemmas 3.5, 3.12 and 3.16 in \cite{Z2}.
The argument of \S4 in \cite{Z2} can be
modified to show the existence of $(v,\,\zeta)$ satisfying
(ii) and this pair is unique up to the action of the
automorphism group $\mathrm{Aut}^\ast(T)\times (S^1)^{|I|}$, cf. the proof of
Lemma 4.4 in \cite{Z3}.\hfill\hb

For any $v=(b, v)\in\widetilde{\mathcal{FT}}^{\emptyset}$
and $h\in\chi(T)$, let
\be \alpha^{(k)}_{T,\, h}(v)=(D^{(k)}_hb)s^{(k)}_{b, \,x_{\mathcal{T}(h)}}\rho_h(v),
\qquad\alpha^{(k)}_T(v)=\sum_{h\in\chi(T)}\alpha^{(k)}_h(v).
\lb{5.14}\ee
We denote $\alpha^{(1)}_{T,\, h}(v)$ and $\alpha_T^{(1)}(v)$ by
$\alpha_{T,\, h}(v)$ and $\alpha_T(v)$ respectively.

We want to analyze the conditions under which a stable
map can be deformed to a $J$-holomorphic map whose domain is smooth.
We have the following first-order estimate:

{ \bf Lemma 5.4.} {\it Suppose $T$ is a bubble type given by (iv) of
Theorem 2.9 with $\Sg_P$ being a smooth Riemann surface of genus
two. Then for every precompact open subset $K$ of $\mathfrak{M}_T(X,
A, J)$, there exist $\delta_K,  C_K\in\R^+$ and an open neighborhood
$U_K$ of $K$ in $\mathfrak{X}_{2, k}(X, A, J)$ satisfying: For every
$v=(b, v)\in\widetilde{\mathcal{FT}}_{\delta_K}^\emptyset|_{\widetilde{K}^{(0)}}$
and $V\eta\in T_{u(\Sg_P)}X\otimes\mathcal{H}^{0, 1}_{\Sg_P}$, we
have \bea\left |\langle\langle \overline{\partial}_Ju_v,\,
R_vV\eta\rangle\rangle_{v, 2}+ \langle\langle\alpha_T(v),
\,V\eta\rangle\rangle_{v, 2}\right| \le
C_K|v|\cdot|\rho(v)|\cdot\|V\eta\|. \nn\eea}

{\bf Proof.} As in Lemma 4.5 of \cite{Z1},  we have
\bea ||s_{b,\,\widetilde{x}_h(v)}((d\phi_{b,\, {\mathcal{T}(h)}}|_{\widetilde{x}_h(v)})^{-1}\rho_h(v))
-s_{b,\,x_{\mathcal{T}(h)}}(\rho_h(v))||_2&\le&
C_K|\phi_{b,\, \mathcal{T}(h)}(\widetilde{x}_h(v))|_{b}|\rho_h(v)|
\nn\\&\le& C_K|v|\cdot|\rho_h(v)|,\nn\\
\sum_{m\ge 2}\left|D_h^{(m)}b\right||\rho_h(v)|^m&\le& C_K|\rho_h(v)|^2,
\nn\eea
for all $h\in\chi(T)$ and
$v=(b, v)\in\widetilde{\mathcal{FT}}_{\delta_K}^\emptyset|_{\widetilde{K}^{(0)}}$
with $\delta_K$ being sufficiently small.
Thus the lemma follows from Lemma 5.1.\hfill\hb

Let $\delta_K$ be given by Lemma 5.4.
For each $v=(b, v)\in\widetilde{\mathcal{FT}}_{\delta_K}^\emptyset|_{\widetilde{K}^{(0)}}$,,
we define the homomorphism
\be \pi_{v, -}^{0, 1}:
\Gamma^{0, 1}(v)\rightarrow\Gamma^{0, 1}_-(b_P),
\qquad \pi_{v, -}^{0, 1}\xi=-\sum_{1\le i\le n,\,1\le j\le 2}\langle \xi, R_ve_i\psi_j\rangle e_i\psi_j\in\Gamma^{0, 1}_-(b_P),
\lb{5.15}\ee
where  $\{\psi_1, \,\psi_2\}$ is an orthonormal basis
for $\mathcal{H}^{0, 1}_{\Sg_P}$ as in (\ref{5.2}) and
$\{e_j\}_{1\le j\le n}$ is an orthonormal basis for
$ T_{u(\Sg_P)}X$. Denote the kernel of  $\pi_{v, -}^{0, 1}$ by $\Gamma_+^{0, 1}(v)$.
Then we have the following:

{ \bf Lemma 5.5.} {\it Suppose $T$ is a bubble type given by (iv) of
Theorem 2.9 with $\Sg_P$ being a smooth Riemann surface of genus
two. Then an element $b\equiv[\mathcal{C},\, u]\in \mathfrak{M}_T(X,
A, J)\cap\overline{\mathfrak{M}}^0_{2, k}(X, A, J)$ must satisfy
$\alpha_T(v)=0$ for some $v=(b,
v)\in\widetilde{\mathcal{FT}}_{\delta_K}^\emptyset|_{\widetilde{K}^{(0)}}$,
where $\delta_K$ is given by Lemma 5.4. }

{\bf Proof.} By (ii) of Lemma 5.3, \bea U_T=\{[\exp_{u_v}\zeta]:
v=(b,
v)\in\widetilde{\mathcal{FT}}_{\delta_K}|_{\widetilde{K}^{(0)}},\;
\zeta\in\widetilde{\Gamma}_+(v),\;\|\zeta\|_{v, p, 1}<\delta_K\}
\nn\eea is an open neighborhood of $K$ in $\mathfrak{X}_{2, k}(X, A,
J)$. Now suppose $\overline{\partial}_J\exp_{u_v}\zeta=0$ for a pair
$(v,\,\zeta)$. Write \be
\overline{\partial}_J\exp_{u_v}\zeta=\overline{\partial}_Ju_v+D_v\zeta+N_v\zeta,
\lb{5.16}\ee where $N_v$ is a quadratic form satisfying (cf. Theorem
3a of \cite{F}) \be N_v0=0,\qquad \|N_v\xi-N_v\xi^\prime\|_{v,\,p}
\le C_K(\|\xi\|_{v, p, 1}+\|\xi^\prime\|_{v, p,
1})\|\xi-\xi^\prime\|_{v, p, 1}, \lb{5.17}\ee

Now we estimate $\|\overline{\partial}_Ju_v\|_{v,\,p}$. By the
construction of $q_v$, we have $\overline{\partial}_Ju_v=0$ outside
the annuli $A^-_{v, h}\left(|v_h|\right)$ for $h\in\chi(T)$ and
$A^\pm_{v, h}\left(|v_h|\right)$ for $\{h\in I: \chi_Th=2\}$. By the
construction of $q_v$, we have $\|dq_v\|_{C^0}<C(b)$ for some $C\in
C^\infty(\mathfrak{M}_T(X, A, J),\,\R^+)$. Thus we have \bea
\|\overline{\partial}_Ju_v\|_{v,\,p} \le C_K\sum_{h\in
I}\left\|du_h\left|\frac{}{}\right._{\widetilde{A}^\pm_{b,
h}\left(|v_h|\right)}\right\|_{v,\,p} \lb{5.18}\eea Moreover, we
have \bea &&\left\|du_h\left|\frac{}{}\right._{\widetilde{A}^\pm_{b,
h}\left(|v_h|\right)}\right\|_{v,\,p} \nn\\&\le&
C_K\left(\int_{\widetilde{A}^\pm_{b, h}\left(|v_h|
\right)}|du_h|^pd\mu\right)^\frac{1}{p}
+C_K\left(\int_{\widetilde{A}^\pm_{b, h}\left(|v_h|
\right)}|z|^{\frac{2(2-p)}{p}}|du_h|^2d\mu\right)^\frac{1}{2}
\nn\\&\le& C_K\left(\int_{\widetilde{A}^\pm_{b, h}\left(|v_h|
\right)}1d\mu\right)^\frac{1}{p}
+C_K\left(\int_{\widetilde{A}^\pm_{b, h}\left(|v_h|
\right)}|z|^{\frac{2(2-p)}{p}}d\mu\right)^\frac{1}{2} \le
C_K|v_h|^\frac{1}{p}.\lb{5.19}\eea Now (\ref{5.16})-(\ref{5.19}) and
(\ref{5.13}) yield \be \|\zeta\|_{v, p, 1}\le C_K|v|^\frac{1}{p}.
\lb{5.20}\ee Note that \be
\overline{\partial}_J\exp_{u_v}\zeta=0\Leftrightarrow \left\{\matrix
{\pi_{v, -}^{0, 1}(\overline{\partial}_Ju_v+D_v\zeta+N_v\zeta)=0\in
\Gamma^{0, 1}_-(b_P),\cr
\overline{\partial}_Ju_v+D_v\zeta+N_v\zeta=0\in \Gamma^{0,
1}_+(v).\cr}\right. \lb{5.21}\ee Since the structure on $X$ near
$u(\Sg_P)$ is isomorphic to $\C^n$ by assumption, we have
$N_v\zeta=0$ on the support of $R_vV\eta$. Hence we have \be
\langle\langle N_v\zeta,\, R_vV\eta\rangle\rangle_{v,\,2}=0.
\lb{5.22}\ee Thus by Lemmas 5.2 and 5.4 we have \bea \pi_{v, -}^{0,
1}(v, \zeta)\equiv\pi_{v, -}^{0,
1}(\overline{\partial}_Ju_v+D_v\zeta+N_v\zeta)
=\alpha_T(v)+\epsilon(v, \zeta), \lb{5.23}\eea where \bea
\|\epsilon(v, \zeta)\|\le C_K(|v|+\|\zeta\|_{v, p, 1})|\rho(v)| \le
C_K|v|^\frac{1}{p}|\rho(v)|. \lb{5.24}\eea Note that we may choose a
basis $\{\psi_1, \psi_2\}$ of $\mathcal{H}^{0, 1}_{\Sg_P}$ such that
$\psi_i(x_{\mathcal{T}(h)})\neq 0$ for $h\in\chi(T)$ and $i=1, 2$.
Thus in order to satisfies (\ref{5.21}), we must have
$\alpha_T(v)=0$ provided $|v|$ is sufficiently small. This proves
the lemma.\hfill\hb

Now we separate our study into several subsections according to the
number of bubble trees.

\subsection{ There is only one bubble tree}

In this subsection, we consider the case $I_1=1$.
Note that in this case $\mathrm{rank}\alpha_T(v)=2n$
which is less than the dimension of $\mathrm{coker}D_b=4n$.
Thus we need further expansions according the position of the
attaching node of the bubble tree. More precisely, we have
the following cases:

{ \bf Theorem 5.1.1.} {\it Suppose $T$ is a bubble type given by
(iv) of Theorem 2.9 with $\Sg_P$ being smooth and $|\chi(T)|=1$.
Denote by $h$ the single bubble in $\chi(T)$. Assume the attaching
node $x_{\mathcal{T}(h)}$ of the bubble tree is not one of the six
branch points of the canonical map $\Sg_P\rightarrow\P^1:\, x\mapsto
s_{\Sg_P, \,x}$. Then an element $b\equiv[\mathcal{C},\, u]\in
\mathfrak{M}_T(X, A, J)\cap\overline{\mathfrak{M}}^0_{2, k}(X, A, J)$
if and only if $D^{(1)}_hb=D^{(2)}_hb=0$. In particular,
$\mathfrak{M}_T(X, A, J)\cap\overline{\mathfrak{M}}^0_{2, k}(X, A, J)$
is a smooth orbifold of dimension at most $\dim
\overline{\mathfrak{M}}^{vir}_{2,k}(X, A)-2$.  }

{\bf Proof.} By Lemma 5.5, an element $b\equiv[\mathcal{C},\, u]\in
\mathfrak{M}_T(X, A, J)\cap\overline{\mathfrak{M}}^0_{2, k}(X, A, J)$
must satisfy $\alpha_T(v)=0$ for some $v=(b,
v)\in\widetilde{\mathcal{FT}}_{\delta_K}^\emptyset|_{\widetilde{K}^{(0)}}$,
thus we have $D^{(1)}_hb=0$. Denote by $\mathcal{S}=\{b\in
\mathfrak{M}_T(X, A, J): D^{(1)}_hb=0\}$.  Then by (ii-a) of Definition 1.1,
$\mathcal{S}$ is a smooth suborbifold of $\mathfrak{M}_T(X, A, J)$
of codimension $2n$. Let $\mathcal{NS}$ denote the normal bundle of
$\mathcal{S}$ in $\mathfrak{M}_T(X, A, J)$ and identify a small
neighborhood of its zero section with a tubular neighborhood of
$\mathcal{S}$ in $\mathfrak{M}_T(X, A, J)$.

Suppose $(b, N)\in \mathcal{NS}$ and $v=((b,
N);\;v)\in\widetilde{\mathcal{FT}}_{\delta_K}^\emptyset|_{\widetilde{K}^{(0)}}$.
We consider the second-order expansion of $\langle\langle
\overline{\partial}¡ª_Ju_{((b, N);\;v)},\, R_vV\eta\rangle\rangle_{v,
2}$. Note that we have \bea
||s^{(2)}_{b,\,\widetilde{x}_h(v)}((d\phi_{b,\,
{\mathcal{T}(h)}}|_{\widetilde{x}_h(v)})^{-1}\rho_h(v))
-s^{(2)}_{b,\,x_{\mathcal{T}(h)}}(\rho_h(v))||_2&\le& C_K|\phi_{b,\,
\mathcal{T}(h)}(\widetilde{x}_h(v))|_{b}|\rho_h(v)|^2
\nn\\&\le& C_K|v|\cdot|\rho_h(v)|^2,\nn\\
\sum_{m\ge 3}\left|D_h^{(m)}(b, N)\right||\rho_h(v)|^m&\le&
C_K|\rho_h(v)|^3, \nn\\\left|D_h^{(2)}(b, N)-D_h^{(2)}(b,
0)\right|&\le& C_K|N|, \nn\eea for $v,\, N$ sufficiently small by
continuity.

Now by Lemma 5.1 we have \bea&&\left|\frac{}{}\right.\langle\langle
\overline{\partial}¡ª_Ju_{(b, N);\;v},\, R_vV\eta\rangle\rangle_{v, 2}+
\langle\langle D^{(1)}_h(b,
N)s^{(1)}_{b,\,\widetilde{x}_h(v)}((d\phi_{b,\,
{\mathcal{T}(h)}}|_{\widetilde{x}_h(v)})^{-1}\rho_h(v))
\nn\\&&+\alpha^{(2)}_T((b, 0);\;v), \,V\eta\rangle\rangle_{v,
2}\left|\frac{}{}\right. \le C_K|\rho(v)|^2(|v|+|N|)\|V\eta\|.
\lb{5.1.1}\eea

Let $\{\psi_j\}$ be an orthogonal basis for
$\mathcal{H}_{\Sg_P}^{0,\,1}$ such that
$\psi_1\in\mathcal{H}_{\Sg_P}^+(\widetilde{x}_{h}(v))$,
$\psi_2\in\mathcal{H}_{\Sg_P}^-(\widetilde{x}_{h}(v))$ and $\{V_i\}$
an orthogonal basis for $T_{ev_P(b, N)}X$. Note that since
$\zeta\in \widetilde{\Gamma}^+(v)$, we have \bea \langle\langle
\zeta,\,D_v^\ast R_vV_i\psi_1\rangle\rangle_{v,\,2}=0 \lb{5.1.2}\eea
by the construction of $\widetilde{\Gamma}^+(v)$ before Lemma 5.3.
Here we use notations in Lemma 5.5. Since
$\psi_2\in\mathcal{H}_{\Sg_P}^-(\widetilde{x}_{h}(v))$, we have
$\psi_2(\widetilde{x}_{h}(v))=0$. Hence by (\ref{5.10}) we have \bea
|D_v^\ast R_vV_i\psi_2|_{v,\,z}\le C_K|\psi_2|_{(s,\,t)}|\partial_s|
|V_i|\le C_K|\rho(v)|^2. \lb{5.1.3}\eea

By  (\ref{5.1.1})-(\ref{5.1.3}), (\ref{5.20}) and (\ref{5.22}), we
have \bea &&\pi_{v, -}^{0, 1}(v, \zeta)\equiv\pi_{v, -}^{0,
1}(\overline{\partial}_Ju_v+D_v\zeta+N_v\zeta) \nn\\=&&D^{(1)}_h(b,
N)s^{(1)}_{b,\,\widetilde{x}_h(v)}((d\phi_{b,\,
{\mathcal{T}(h)}}|_{\widetilde{x}_h(v)})^{-1}\rho_h(v))
+\alpha^{(2)}_T((b, 0);\;v)+\epsilon(v, \zeta), \nn\\=&&D^{(1)}_h(b,
N)s^{(1)}_{b,\,\widetilde{x}_h(v)}((d\phi_{b,\,
{\mathcal{T}(h)}}|_{\widetilde{x}_h(v)})^{-1}\rho_h(v)) \nn\\&&+
D_h^{(2)}b\left(\overline{\{D_{b,\,x_{\mathcal{T}(h)}}^{(2)}\psi_1\}(\rho_h(v))}\psi_1
+\overline{\{D_{b,\,x_{\mathcal{T}(h)}}^{(2)}\psi_2\}(\rho_h(v))}\psi_2\right)+\epsilon(v,
\zeta), \lb{5.1.4}\eea where \bea \|\epsilon(v, \zeta)\|\le
C_K|\rho(v)|^2(|v|+|N|+\|\zeta\|_{v,\,p,\,1}) \le
C_K|\rho(v)|^2(|N|+|v|^\frac{1}{p}). \lb{5.1.5}\eea Since the
attaching node $x_{\mathcal{T}(h)}$ of the bubble tree is not one of
the six branch points of the canonical map $\Sg_P\rightarrow\P^1:\,
x\mapsto s_{\Sg_P, \,x}$, we have
$D_{b,\,x_{\mathcal{T}(h)}}^{(2)}\psi_2\neq0$. Hence in order to
satisfy (\ref{5.21}), we must have $D_h^{(2)}b=0$ provided $|v|$ and
$|N|$ are sufficiently small.

Conversely, suppose $D^{(1)}_hb=D^{(2)}_hb=0$. We want to construct
$\phi(v)\in\mathfrak{M}^0_{2, k}(X, A, J)$ converging to $b$. We
have the following: In a small neighborhood of $b$, for any \bea
v=((\widetilde{b},
N);\;v)\in\widetilde{\mathcal{FT}}_{\delta_K}^{\emptyset}|_{\widetilde{K}^{(0)}},
\quad \|\zeta\|_{v, p, 1}\le 2C_K|v|^\frac{1}{p} \nn\eea the
equation \bea &&\pi_{v, -}^{0, 1}(((\widetilde{b}, N);\;v),\; \zeta)
\nn\\&\equiv&D^{(1)}_h(\widetilde{b},
N)s^{(1)}_{\widetilde{b},\,\widetilde{x}_h(v)}((d\phi_{\widetilde{b},\,
{\mathcal{T}(h)}}|_{\widetilde{x}_h(v)})^{-1}\rho_h(v))
+\alpha^{(2)}_T((\widetilde{b}, 0);\;v)+\epsilon(v, \zeta)=0
\lb{5.1.6}\eea has a unique small solution $(b^\ast, N^\ast)\in
\widetilde{\mathcal{M}}^{(0)}_{T}$.

In fact, let $\pi^\pm_{\widetilde{x}_{h}(v)}$ be the projections to
$\mathcal{H}_{\Sg_P}^\pm(\widetilde{x}_{h}(v))$ respectively and
consider the equation \bea \pi^-_{\widetilde{x}_{h}(v)}\circ\pi_{v,
-}^{0, 1}(((\widetilde{b}, N);\;v),\; \zeta)
\equiv(D_h^{(2)}\widetilde{b})
\overline{\{D_{\widetilde{b},\,x_{\mathcal{T}(h)}}^{(2)}\psi_2\}(\rho_h(v))}\psi_2+\pi^-_{\widetilde{x}_{h}(v)}\circ\epsilon(v,
\zeta)=0. \lb{5.1.7}\eea Note that $D^{(1)}_h\widetilde{b}=d\widetilde{u}_h|_\infty=0$,
thus by (ii-b) of Definition 1.1,
and $D_{b,\,x_{\mathcal{T}(h)}}^{(2)}\psi_2\neq0$, the map
\bea\Phi^-: \widetilde{\mathcal{FT}}_{\delta_K}^{\emptyset}
\rightarrow
T_{u(\Sg_P)}X\otimes\mathcal{H}^{-}_{\Sg_P}(\widetilde{x}_{h}(v)),
\quad v\mapsto
(D_h^{(2)}\widetilde{b})\overline{\{D_{\widetilde{b},\,x_{\mathcal{T}(h)}}^{(2)}\psi_2\}(\rho_h(v))}\psi_2
\nn\eea is transversal to the zero section. This together with
(\ref{5.1.5}) yields a unique solution $(b^\ast, N)$ of
(\ref{5.1.7}) for each $N$ provided $|v|$ and $|N|$ are sufficiently
small.

Now consider the equation \bea
&&\pi^+_{\widetilde{x}_{h}(v)}\circ\pi_{v, -}^{0, 1}(((b^\ast,
N);\;v),\; \zeta) \equiv D^{(1)}_h(b^\ast,
N)s^{(1)}_{b^\ast,\,\widetilde{x}_h(v)}((d\phi_{b^\ast,\,
{\mathcal{T}(h)}}|_{\widetilde{x}_h(v)})^{-1}\rho_h(v))
\nn\\&&\qquad+
D_h^{(2)}b^\ast\overline{\{D_{b^\ast,\,x_{\mathcal{T}(h)}}^{(2)}\psi_1\}(\rho_h(v))}\psi_1+\pi^+_{\widetilde{x}_{h}(v)}\circ\epsilon(v,
\zeta)=0 \lb{5.1.8}\eea By (ii-a) of Definition 1.1, the map \bea\Phi^+:
\widetilde{\mathcal{FT}}_{\delta_K}^{\emptyset} \rightarrow
T_{u(\Sg_P)}X\otimes\mathcal{H}^{+}_{\Sg_P}(\widetilde{x}_{h}(v)),
\quad N\mapsto D^{(1)}_h(b^\ast,
N)s^{(1)}_{b^\ast,\,\widetilde{x}_h(v)}((d\phi_{b^\ast,\,
{\mathcal{T}(h)}}|_{\widetilde{x}_h(v)})^{-1}\rho_h(v))\nn\eea is
transversal to the zero section. Note that \bea
\left\|D_h^{(2)}b^\ast\overline{\{D_{b^\ast,\,x_{\mathcal{T}(h)}}^{(2)}\psi_1\}(\rho_h(v))}\psi_1\right\|\le
C_K|\rho(v)|^2. \nn\eea This together with (\ref{5.1.5}) yields a
unique solution $(b^\ast, N^\ast)$ of (\ref{5.1.8}) provided $|v|$
is sufficiently small. Thus $(b^\ast, N^\ast)$ is the unique
solution of (\ref{5.1.6}) as desired. Denote by $\mu(v,
\zeta)=((b^\ast, N^\ast),\, v)$.

Now we define the map \bea &\Psi_v: \{\zeta\in
\widetilde{\Gamma}_+(\mu(v, \zeta)): \zeta\le
2C_K|v|^\frac{1}{p}\}\rightarrow\Gamma^{0, 1}_+(\mu(v, \zeta)) \nn\\
&\Psi_v(\zeta)=\overline{\partial}¡ª_Ju_{\mu(v, \zeta)}+D_{\mu(v,
\zeta)}\zeta+N_{\mu(v, \zeta)}\zeta \nn\eea Since the derivative
\bea D\Psi_v(0): \widetilde{\Gamma}_+(\mu(v,
\zeta))\rightarrow\Gamma^{0, 1}_+(\mu(v, \zeta)) \nn\eea is an
isomorphism and $\|\Psi_v(0)\|_{\mu(v, 0), p, 1}\le
2C_K|v|^\frac{1}{p}$. The equation $\Psi_v(\zeta)=0$ has a unique
small solution $\zeta_v$ by the contraction principle. We define the
map \bea \phi: \{(b,
v)\in\widetilde{\mathcal{FT}}_{\delta_K}^\emptyset|_{\widetilde{K}^{(0)}}:
D^{(1)}_hb=0=D^{(2)}_hb \}\rightarrow \mathfrak{M}^0_{2,\,
k}(X, A, J),\qquad \phi(v)=[\exp_{\mu(v,\, \zeta_v) }\zeta_v].
\nn\eea Then $\phi(v)$ converges to $b$ in the stable map topology,
Hence $b\in\mathfrak{M}_T(X, A, J)\cap\overline{\mathfrak{M}}^0_{2,\,
k}(X, A, J)$. The proof of the theorem is complete.\hfill\hb

{ \bf Theorem 5.1.2.} {\it Suppose $T$ is a bubble type given by
(iv) of Theorem 2.9 with $\Sg_P$ being smooth and $|\chi(T)|=1$.
Denote by $h$ the single bubble in $\chi(T)$. Assume the attaching
node $x_{\mathcal{T}(h)}$ of the bubble tree is one of the six
branch points $\{z_m\}_{1\le m\le 6}$ of the canonical map.
Then we have the following:

(i) If $H^1(\P^1, u_h^\ast TX\otimes \mathcal{O}_{\P^1}(-4\infty))=0$,
then an element $b\equiv[\mathcal{C},\, u]\in
\mathfrak{M}_T(X, A, J)\cap\overline{\mathfrak{M}}^0_{2, k}(X, A, J)$
if and only if $\alpha_T(v)=0$ for some $v=(b,
v)\in\widetilde{\mathcal{FT}}_{\delta_K}^\emptyset|_{\widetilde{K}^{(0)}}$
and $\widetilde{\alpha}_T(w, v)=0$, where $\widetilde{\alpha}_T(w,
v)$ is a linear combination of $D^{(2)}_{h}b$ and $D_{h} ^{(3)}b$
with coefficients depending on $w\in T_{z_m}\Sg_P$, $v$ and the
position of the nodes in bubble trees. In particular,
$\mathfrak{M}_T(X, A, J)\cap\overline{\mathfrak{M}}^0_{2, k}(X, A, J)$
is a smooth orbifold of dimension at most $\dim
\overline{\mathfrak{M}}^{vir}_{2,k}(X, A)-2$
.

(ii) If $H^1(\P^1, u_h^\ast TX\otimes \mathcal{O}_{\P^1}(-4\infty))\neq0$,
then an element $b\equiv[\mathcal{C},\, u]\in
\mathfrak{M}_T(X, A, J)\cap\overline{\mathfrak{M}}^0_{2, k}(X, A, J)$
must satisfy: $u_h$ factor through a branched covering $\widetilde{u}: S^2\rightarrow
X$, i.e., there exists a holomorphic branched covering $\phi:
S^2\rightarrow S^2$ such that $u_h=\widetilde{u}\circ\phi$ and
$\deg(\phi)\ge 2$.
}

{\bf Proof.} By Lemma 5.5, an element $b\equiv[\mathcal{C},\, u]\in
\mathfrak{M}_T(X, A, J)\cap\overline{\mathfrak{M}}^0_{2, k}(X, A, J)$
must satisfy $\alpha_T(v)=0$ for some $v=(b,
v)\in\widetilde{\mathcal{FT}}_{\delta_K}^\emptyset|_{\widetilde{K}^{(0)}}$,
thus we have $D^{(1)}_hb=0$. Denote by $\mathcal{S}_m=\{b\in
\mathfrak{M}_T(X, A, J): D^{(1)}_hb=0,\;x_{\mathcal{T}(h)}=z_m\}$,
where $\{z_m\}_{1\le m\le 6}$ are the six branch points of
$s_{\Sg_P}$. Note that $z_m$ depends on the complex structure
$j_{\Sg_P}$ on $\Sg_P$, in order to simplify notations, we use $z_m$
to denote them. Clearly, each $\mathcal{S}_m$ is a smooth
suborbifold of $\mathfrak{M}_T(X, A, J)$ of codimension $2n+2$. Let
$\mathcal{NS}_m$ denote the normal bundle of $\mathcal{S}_m$ in
$\mathcal{S}$, where $\mathcal{S}$ is given by Theorem 5.1.1 and
identify a small neighborhood of its zero section with a tubular
neighborhood of $\mathcal{S}_m$ in $\mathcal{S}$.

Suppose $(b,\,w)\in \mathcal{NS}_m$, $(b, w, N)\in \mathcal{NS}$ and
$v=((b, w,
N);\;v)\in\widetilde{\mathcal{FT}}_{\delta_K}^\emptyset|_{\widetilde{K}^{(0)}}$.
We consider the third-order expansion of $\langle\langle
\overline{\partial}¡ª_Ju_{((b, w, N);\;v)},\, R_vV\eta\rangle\rangle_{v,
2}$. Note that we have \bea
||s^{(3)}_{b,\,\widetilde{x}_h(v)}((d\phi_{b,\,
{\mathcal{T}(h)}}|_{\widetilde{x}_h(v)})^{-1}\rho_h(v))
-s^{(3)}_{b,\,z_m}(\rho_h(v))||_2&\le& C_K|\phi_{b,\,
\mathcal{T}(h)}(\widetilde{x}_h(v))|_{b}|\rho_h(v)|^3
\nn\\&\le& C_K|v|\cdot|\rho_h(v)|^3,\nn\\
\sum_{m\ge 4}\left|D_h^{(m)}(b, w, N)\right||\rho_h(v)|^m&\le&
C_K|\rho_h(v)|^4, \nn\\\left|D_h^{(3)}(b, w, N)-D_h^{(3)}(b, 0,
0)\right|&\le& C_K(|N|+|w|), \nn\eea for $v,\,w,\, N$ sufficiently
small by continuity.

Now by Lemma 5.1 we have \bea&&\left|\frac{}{}\right.\langle\langle
\overline{\partial}¡ª_Ju_{((b, w, N);\;v)},\, R_vV\eta\rangle\rangle_{v,
2}+ \langle\langle D^{(1)}_h(b, w,
N)s^{(1)}_{b,\,\widetilde{x}_h(v)}((d\phi_{b,\,
{\mathcal{T}(h)}}|_{\widetilde{x}_h(v)})^{-1}\rho_h(v))
\nn\\&&\quad+D^{(2)}_h(b, w,
N)s^{(2)}_{b,\,\widetilde{x}_h(v)}((d\phi_{b,\,
{\mathcal{T}(h)}}|_{\widetilde{x}_h(v)})^{-1}\rho_h(v))
+\alpha^{(3)}_T((b, 0, 0);\;v), \,V\eta\rangle\rangle_{v,
2}\left|\frac{}{}\right. \nn\\\le&&
C_K|\rho(v)|^3(|v|+|w|+|N|)\|V\eta\|. \lb{5.1.9}\eea Note that we
have $\widetilde{x}_h(v)=\widetilde{x}_h((b, w,
N);\;v)=\widetilde{x}_h((b, w, 0);\;v)\equiv\widetilde{x}_h(w,
v)\in\Sg_P$. Identify a small neighborhood of $z_m$ in $\Sg_P$ with
a small neighborhood of $0$ in $T_{z_m}\Sg_P$, then we have
$|\widetilde{x}_h(w, v)|\le C_K(|v|+|w|)$. Let $s^{(2,
-)}\in\Gamma(\Sg_P,\,T^\ast\Sg_P^{\otimes
2}\otimes\mathcal{H}^-_{\Sg_P})$ be the projection of the section
$s^{(2)}_{b,\, x}$ onto the subbundle $\mathcal{H}^-_{\Sg_P}$. Then
we have $s^{(2, -)}$ is independent of the metric on $\Sg_P$ and has
transversal zeros at the six points $\{z_m\}_{1\le m\le 6}$ (cf.
P.402 of \cite{Z1}).

Let $\{\psi_j\}$ be an orthonormal  basis for
$\mathcal{H}_{\Sg_P}^{0,\,1}$ such that
$\psi_1\in\mathcal{H}_{\Sg_P}^+(\widetilde{x}_{h}(w, v))$,
$\psi_2\in\mathcal{H}_{\Sg_P}^-(\widetilde{x}_{h}(w, v))$ and
$\{V_i\}$ an orthonormal basis for $T_{ev_P(b, w, N)}X$. Note
that since $\zeta\in \widetilde{\Gamma}^+(v)$, we have \bea
\langle\langle \zeta,\,D_v^\ast R_vV_i\psi_1\rangle\rangle_{v,\,2}=0
\lb{5.1.10}\eea Since
$\psi_2\in\mathcal{H}_{\Sg_P}^-(\widetilde{x}_{h}(w, v))$, we have
$\psi_2(\widetilde{x}_{h}(w, v))=0$. Note that
$$s^{(2, -)}(\widetilde{x}_{h}(w, v))=\left(\overline{D^{(2)}_{b,\,\widetilde{x}_{h}(w, v)}\psi_2}\right)\psi_2,
\qquad s^{(2, -)}(z_m)=0.$$ Thus we have $\|\nabla\psi_2\|\le
C_K|\widetilde{x}_{h}(w, v)|$ for $|v|$ and $|w|$ small enough.
Hence by (\ref{5.10}) we have \bea |D_v^\ast
R_vV_i\psi_2|_{v,\,z}\le C_K|\psi_2|_{(s,\,t)}|\partial_s| |V_i|\le
C_K|\rho(v)|^2(|\rho(v)|+|\widetilde{x}_{h}(w, v)|). \lb{5.1.11}\eea
Note that since $s^{(2, -)}(z_m)=0$ and $|\pi^-_{\widetilde{x}_h(w,
v)}-\pi^-_{z_m}|\le C_K|\rho_h(v)|^2$, we have \bea \left|s^{(2,
-)}_{b,\,\widetilde{x}_h(w, v)}((d\phi_{b,\,
{\mathcal{T}(h)}}|_{\widetilde{x}_h(v)})^{-1}\rho_h(v))
-s^{(3,\,-)}_{b,\,z_m}(\widetilde{x}_h(w, v), \rho_h(v),
\rho_h(v))\right|\le C_K|\widetilde{x}_{h}(w, v)|^2|\rho_h(v)|^2.
\lb{5.1.12}\eea By  (\ref{5.1.9})-(\ref{5.1.12}), (\ref{5.20}) and
(\ref{5.22}), we have \bea &&\pi_{v, -}^{0, 1}(v,
\zeta)\equiv\pi_{v, -}^{0,
1}(\overline{\partial}_Ju_v+D_v\zeta+N_v\zeta) \nn\\=&&D^{(1)}_h(b, w,
N)s^{(1)}_{b,\,\widetilde{x}_h(v)}((d\phi_{b,\,
{\mathcal{T}(h)}}|_{\widetilde{x}_h(v)})^{-1}\rho_h(v))
+D^{(2)}_h(b, w, N)s^{(2,\,+)}_{b,\,\widetilde{x}_h(v)}((d\phi_{b,\,
{\mathcal{T}(h)}}|_{\widetilde{x}_h(v)})^{-1}\rho_h(v))
\nn\\&&\quad+D^{(2)}_h(b, 0,
0)s^{(3,\,-)}_{b,\,z_m}(\widetilde{x}_h(w, v), \rho_h(v), \rho_h(v))
+D^{(3)}_h(b, 0, 0)s^{(3,\,+)}_{b,\,z_m}(\rho_h(v))
\nn\\&&\quad+D^{(3)}_h(b, 0, 0)s^{(3,\,-)}_{b,\,z_m}(\rho_h(v))
+\epsilon(v, \zeta), \lb{5.1.13}\eea where \bea \|\epsilon(v,
\zeta)\|\le&&
C_K|\rho(v)|^3(|v|+|w|+|N|)+C_K|\rho(v)|^2(|\rho(v)|+|\widetilde{x}_{h}(w,
v)|)\cdot\|\zeta\|_{v,\,p,\,1} \nn\\&&+C_K|\widetilde{x}_{h}(w,
v)||\rho(v)|^2(|\widetilde{x}_{h}(w, v)|+|w|+|N|)\lb{5.1.14}\eea
Since $s^{(2, -)}$  has transversal zeros at $z_m$, we have
$s^{(3,\,-)}_{b,\,z_m}\neq 0$. Hence in order to satisfies
(\ref{5.21}), we must have \bea \widetilde{\alpha}_T(w, v)\equiv
D^{(2)}_h(b, 0, 0)s^{(3,\,-)}_{b,\,z_m}(\widetilde{x}_h(w, v),
\rho_h(v), \rho_h(v)) +D^{(3)}_h(b, 0,
0)s^{(3,\,-)}_{b,\,z_m}(\rho_h(v))=0 \lb{5.1.15}\eea provided $|v|$,
$|w|$ and $|N|$ are sufficiently small,

Conversely, suppose $\alpha_T(v)=0$ and $\widetilde{\alpha}_T(w, v)=0$.
We want to construct  $\phi(w, v)\in\mathfrak{M}^0_{2, k}(X, A, J)$
converging to $b$. Note that by (ii-d) of Definition 1.1, (ii) of Theorem 5.1.2 holds.
Hence we only need to consider the case
$H^1(\P^1, u_h^\ast TX\otimes \mathcal{O}_{\P^1}(-4\infty))=0$.
We have the following: In a small
neighborhood of $b$, for any \bea v=((\widetilde{b}, w,
N);\;v)\in\widetilde{\mathcal{FT}}_{\delta_K}^{\emptyset}|_{\widetilde{K}^{(0)}},
\quad |w|\le\delta_K,\quad \|\zeta\|_{v, p, 1}\le
2C_K|v|^\frac{1}{p} \nn\eea the equation (\ref{5.1.13}) with $b$
replaced by $\widetilde{b}$ has a unique small solution $(b^\ast, w,
N^\ast)\in \widetilde{\mathcal{M}}^{(0)}_{T}$.

In fact, let $\pi^\pm_{\widetilde{x}_{h}(w, v)}$ be the projections
to $\mathcal{H}_{\Sg_P}^\pm(\widetilde{x}_{h}(w, v))$ respectively
and consider the equation \bea &&\pi^-_{\widetilde{x}_{h}(w,
v)}\circ\pi_{v, -}^{0, 1}(((\widetilde{b}, w, N);\;v),\; \zeta)
\equiv D^{(2)}_h(\widetilde{b}, 0,
0)s^{(3,\,-)}_{\widehat{b},\,z_m}(\widetilde{x}_h(w, v), \rho_h(v),
\rho_h(v)) \nn\\&&\qquad+D^{(3)}_h(\widetilde{b}, 0,
0)s^{(3,\,-)}_{\widetilde{b},\,z_m}(\rho_h(v))+\pi^-_{\widetilde{x}_{h}(w,
v)}\circ\epsilon(v, \zeta)=0 \lb{5.1.16}\eea By the assumption
$H^1(\P^1, u_h^\ast TX\otimes \mathcal{O}_{\P^1}(-4\infty))=0$, the map \bea&&\Phi^-:
\widetilde{\mathcal{FT}}_{\delta_K}^{\emptyset} \rightarrow
T_{u(\Sg_P)}X\otimes\mathcal{H}^{-}_{\Sg_P}(\widetilde{x}_{h}(w,
v)), \nn\\&&v\mapsto D^{(2)}_h(\widetilde{b}, 0,
0)s^{(3,\,-)}_{\widetilde{b},\,z_m}(\widetilde{x}_h(w, v),
\rho_h(v), \rho_h(v)) +D^{(3)}_h(\widetilde{b}, 0,
0)s^{(3,\,-)}_{\widetilde{b},\,z_m}(\rho_h(v)) \nn\eea is
transversal to the zero section. This together with (\ref{5.1.14})
yields a unique solution $(b^\ast, w, N)$ of (\ref{5.1.16}) for each
$N$ small enough.

Now consider the equation \bea &&\pi^+_{\widetilde{x}_{h}(w,
v)}\circ\pi_{v, -}^{0, 1}(((b^\ast, w, N);\;v),\; \zeta) \equiv
D^{(1)}_h(b^\ast, w, N)s^{(1)}_{b^\ast,\,\widetilde{x}_h(w,
v)}((d\phi_{b^\ast,\, {\mathcal{T}(h)}}|_{\widetilde{x}_h(w,
v)})^{-1}\rho_h(v)) \nn\\&&\qquad+D^{(2)}_h(b^\ast, w,
N)s^{(2,\,+)}_{b^\ast,\,\widetilde{x}_h(w, v)}((d\phi_{b^\ast,\,
{\mathcal{T}(h)}}|_{\widetilde{x}_h(w, v)})^{-1}\rho_h(v))
\nn\\&&\qquad+D^{(3)}_h(b^\ast, 0,
0)s^{(3,\,+)}_{b^\ast,\,z_m}(\rho_h(v)) +\pi^+_{\widetilde{x}_{h}(w,
v)}\circ\epsilon(v, \zeta)=0. \lb{5.1.17}\eea
By (ii-a) of Definition 1.1, the map
\bea &&\Phi^+: \widetilde{\mathcal{FT}}_{\delta_K}^{\emptyset}
\rightarrow
T_{u(\Sg_P)}X\otimes\mathcal{H}^{+}_{\Sg_P}(\widetilde{x}_{h}(v)),
\nn\\&&\qquad N\mapsto D^{(1)}_h(b^\ast, w,
N)s^{(1)}_{b^\ast,\,\widetilde{x}_h(w, v)}((d\phi_{b^\ast,\,
{\mathcal{T}(h)}}|_{\widetilde{x}_h(w, v)})^{-1}\rho_h(v)) \nn\eea
is transversal to the zero section. Note that \bea
\left\|D^{(2)}_h(b^\ast, w,
N)s^{(2,\,+)}_{b^\ast,\,\widetilde{x}_h(w, v)}((d\phi_{b^\ast,\,
{\mathcal{T}(h)}}|_{\widetilde{x}_h(w, v)})^{-1}\rho_h(v))
+D^{(3)}_h(b^\ast, 0, 0)s^{(3,\,+)}_{b^\ast,\,z_m}(\rho_h(v))
\right\|\le C_K|\rho(v)|^2. \nn\eea This together with
(\ref{5.1.14}) yields a unique solution $(b^\ast, w, N^\ast)$ of
(\ref{5.1.17}) provided $\delta_K$ is sufficiently small. Thus
$(b^\ast, w, N^\ast)$ is the unique solution of (\ref{5.1.13}) as
desired. Let $\mu(v, \zeta)=((b^\ast, w, N^\ast),\, v)$. Then the
same argument in Theorem 5.1.1 yields the $J$-holomorphic map $\phi(w,
v)=[\exp_{\mu(v,\, \zeta_v) }\zeta_v]\in\mathfrak{M}^0_{2,
k}(\P^n,\, d)$. The proof of the theorem is complete.\hfill\hb

Now we consider the general bubble tree case. Denote by $\chi(T)=\{h_1, h_2,\ldots,h_p\}$.

{ \bf Theorem 5.1.3.} {\it Suppose $T$ is a bubble type given by
(iv) of Theorem 2.9 with $\Sg_P$ being smooth together with
$|I_1|=1$ and $|\chi(T)|\ge 2$. Assume the attaching node
$\widehat{x}$ of the bubble tree is not one of the six branch points
of the canonical map. Then an element $b\equiv[\mathcal{C},\, u]\in
\mathfrak{M}_T(X, A, J)\cap\overline{\mathfrak{M}}^0_{2, k}(X, A, J)$
if and only if $\alpha_T(v)=0$ for some $v=(b,
v)\in\widetilde{\mathcal{FT}}_{\delta_K}^\emptyset|_{\widetilde{K}^{(0)}}$
and $\widetilde{\alpha}_T(v)=0$, where $\widetilde{\alpha}_T(v)$ is
a linear combination of $\{D^{(1)}_{h_i}b,\,D_{h_i}
^{(2)}b\}_{h_i\in\chi(T)}$ with coefficients depending on $v$ and
the position of the nodes in bubble trees. Moreover, $\alpha_T(v)$
and $\widetilde{\alpha}_T(v)$ are linear independent. In particular,
$\mathfrak{M}_T(X, A, J)\cap\overline{\mathfrak{M}}^0_{2, k}(X, A, J)$
is a smooth orbifold of dimension at most $\dim
\overline{\mathfrak{M}}^{vir}_{2,k}(X, A)-2$.  }

{\bf Proof.} We only give the proof of the simplest case that all of
$\{h_1, h_2,\ldots,h_p\}$ are attached to the bubble $\widehat{h}\in
I_1$. The general case follows by a similar argument and we will
sketch its proof in the end.

By Lemma 5.5, an element $b\equiv[\mathcal{C},\, u]\in
\mathfrak{M}_T(X, A, J)\cap\overline{\mathfrak{M}}^0_{2, k}(X, A, J)$
must satisfy $\alpha_T(v)=0$ for some $v=(b,
v)\in\widetilde{\mathcal{FT}}_{\delta_K}^\emptyset|_{\widetilde{K}^{(0)}}$,
where
\bea\alpha_T(v)=\sum_{h_i\in\chi(T)}(D^{(1)}_{h_i}b)s_{\Sg_P,\,\widehat{x}}(\rho_{h_i}(v)).
\lb{5.1.18}\eea Denote by $\mathcal{S}=\{b\in \mathfrak{M}_T(X, A, J):
\alpha_T(v)=0\}$. Then by (ii-a) of Definition 1.1, $\mathcal{S}$ is smooth
suborbifold of $\mathfrak{M}_T(X, A, J)$. Let $\mathcal{NS}$ denote
the normal bundle of $\mathcal{S}$ in $\mathfrak{M}_T(X, A, J)$.

Suppose $(b, N)\in \mathcal{NS}$ and $v=((b,
N);\;v)\in\widetilde{\mathcal{FT}}_{\delta_K}^\emptyset|_{\widetilde{K}^{(0)}}$.
We consider the second-order expansion of $\langle\langle
\overline{\partial}¡ª_Ju_{((b, N);\;v)},\, R_vV\eta\rangle\rangle_{v,
2}$. As in Theorem 5.1.1, we have
\bea&&\left|\frac{}{}\right.\langle\langle \overline{\partial}¡ª_Ju_{(b,
N);\;v},\, R_vV\eta\rangle\rangle_{v, 2}+ \langle\langle
\sum_{h_i\in\chi(T)}D^{(1)}_{h_i}(b,
N)s^{(1)}_{b,\,\widetilde{x}_{h_i}(v)}((d\phi_{b,\,
{\mathcal{T}(h_i)}}|_{\widetilde{x}_{h_i}(v)})^{-1}\rho_{h_i}(v))
\nn\\&&\quad+\alpha^{(2)}_T((b, 0);\;v), \,V\eta\rangle\rangle_{v,
2}\left|\frac{}{}\right. \le C_K|\rho(v)|^2(|v|+|N|)\|V\eta\|.
\lb{5.1.19}\eea Let $\widetilde{x}_{h_i}((b,
N);\;v)=\widetilde{x}_{h_i}((b,
0);\;v)\equiv\widetilde{x}_{h_i}(v)\in\Sg_P$ for $h_i\in\chi(T)$.
Identify a small neighborhood of $\widehat{x}$ in $\Sg_P$ with a
small neighborhood of $0$ in $T_{\widehat{x}}\Sg_P$ and let
$x_i^\ast(v)=\widetilde{x}_{h_i}(v)-\widetilde{x}_{h_1}(v)
=v_{\widehat{h}}(x_{h_i}-x_{h_1})$ and
$x^\ast(v)=(x^\ast_2(v),\ldots,x^\ast_p(v))$.

Let $\{\psi_j\}$ be an orthonormal basis for
$\mathcal{H}_{\Sg_P}^{0,\,1}$ such that
$\psi_1\in\mathcal{H}_{\Sg_P}^+(\widetilde{x}_{h_1}(v))$,
$\psi_2\in\mathcal{H}_{\Sg_P}^-(\widetilde{x}_{h_1}(v))$ and
$\{V_i\}$ an orthonormal basis for $T_{ev_P(b, N)}X$. Note that
since $\zeta\in \widetilde{\Gamma}^+(v)$, we have \bea
\langle\langle \zeta,\,D_v^\ast R_vV_i\psi_1\rangle\rangle_{v,\,2}=0
\lb{5.2,20}\eea Since
$\psi_2\in\mathcal{H}_{\Sg_P}^-(\widetilde{x}_{h_1}(v))$, we have
$\psi_2(\widetilde{x}_{h_1}(v))=0$. Hence by (\ref{5.10}) we have
\bea |D_v^\ast R_vV_i\psi_2|_{v,\,z}\le
C_K|\psi_2|_{(s,\,t)}|\partial_s| |V_i|\le
C_K|\rho(v)|(|\rho(v)|+|x^\ast(v)|). \lb{5.1.21}\eea Note that we
have \bea &&\left|s^{(1)}_{b,\,\widetilde{x}_{h_i}(v)}((d\phi_{b,\,
{\mathcal{T}(h_i)}}|_{\widetilde{x}_{h_i}(v)})^{-1}\rho_{h_i}(v))
-s^{(1)}_{b,\,\widetilde{x}_{h_1}(v)}(\rho_{h_i}(v))
-s^{(2)}_{b,\,\widehat{x}}(x_i^\ast(v), \rho_{h_i}(v))
\right|\nn\\\le &&C_K(|x_i^\ast(v)||\rho(v)|(|x_i^\ast(v)|+|v|).
\lb{5.1.22}\eea By (\ref{5.1.19})-(\ref{5.1.22}), we have \bea
\pi_{v, -}^{0, 1}(v, \zeta)&\equiv&\pi_{v, -}^{0,
1}(\overline{\partial}_Ju_v+D_v\zeta+N_v\zeta) =D^{(1)}_{h_1}(b,
N)s^{(1)}_{b,\,\widetilde{x}_{h_1}(v)}((d\phi_{b,\,
{\mathcal{T}(h_1)}}|_{\widetilde{x}_{h_1}(v)})^{-1}\rho_{h_1}(v))
\nn\\&&+\sum_{h_i\in\chi(T)\setminus h_1} (D^{(1)}_{h_i}(b,
N)s^{(1)}_{b,\,\widetilde{x}_{h_1}(v)}(\rho_{h_i}(v))
+D^{(1)}_{h_i}(b, 0)s^{(2)}_{b,\,\widehat{x}}(x_i^\ast(v),
\rho_{h_i}(v))) \nn\\&&+\sum_{h_i\in\chi(T)} (D^{(2)}_{h_i}(b,
0)s^{(2, +)}_{b,\,\widehat{x}}(\rho_{h_i}(v)) +D^{(2)}_{h_i}(b,
0)s^{(2, -)}_{b,\,\widehat{x}}(\rho_{h_i}(v))) +\epsilon(v, \zeta),
\lb{5.1.23}\eea where \bea \|\epsilon(v,
\zeta)\|\le&&C_K|\rho(v)|^2(|v|+|N|)
+C_K|\rho(v)|(|\rho(v)|+|x^\ast(v)|)\cdot\|\zeta\|_{v,\,p,\,1}
\nn\\&&+C_K(|x^\ast(v)||\rho(v)|(|x^\ast(v)|+|v|+|N|).\lb{5.1.24}\eea
Note that $\rho_{h_i}(v)=v_{\widehat{h}}v_{h_i}$ and
$x^\ast_i(v)=v_{\widehat{h}}(x_{h_i}-x_{h_1})\neq 0$ for
$h_i\in\chi(T)\setminus h_1$. Hence \be|s^{(2,
-)}_{b,\,\widehat{x}}(\rho_{h_i}(v))| =o(|s^{(2,
-)}_{b,\,\widehat{x}}(x_i^\ast(v), \rho_{h_i}(v))|), \quad\forall
h_i\in\chi(T)\setminus h_1,\lb{5.1.25}\ee where we denote by $o(x)$
the higher order term of $x$ as $x\rightarrow0$.
By changing the order of $\{h_1, \ldots, h_p\}$,
we may assume $\rho_{h_1}(v)=\min_{1\le i\le p}\rho_{h_i}(v)$.
Thus we have
\be|s^{(2, -)}_{b,\,\widehat{x}}(\rho_{h_1}(v))|
=\sum_{h_i\in\chi(T)\setminus h_1}o(|s^{(2,
-)}_{b,\,\widehat{x}}(x_i^\ast(v), \rho_{h_i}(v))|).\lb{5.1.26}\ee
Since $\widehat{x}$ is not one of the six branch points of the canonical
map, we have $s^{(2, -)}_{b,\,\widehat{x}}\neq 0$. Thus in order to
satisfy (\ref{5.21}), we must have \bea
\widetilde{\alpha}_T(v)\equiv\sum_{h_i\in\chi(T)\setminus h_1}
D^{(1)}_{h_i}(b, 0)s^{(2, -)}_{b,\,\widehat{x}}(x_i^\ast(v),
\rho_{h_i}(v))=0. \lb{5.1.27}\eea

Conversely, suppose $\alpha_T(v)=0=\widetilde{\alpha}_T(v)$ We want
to construct  $\phi(v)\in\mathfrak{M}^0_{2, k}(X, A, J)$
converging to $b$. We have the following: In a small neighborhood of
$b$, for any \bea v=((\widetilde{b},
N);\;v)\in\widetilde{\mathcal{FT}}_{\delta_K}^{\emptyset}|_{\widetilde{K}^{(0)}},
\quad \|\zeta\|_{v, p, 1}\le 2C_K|v|^\frac{1}{p} \nn\eea the
equation (\ref{5.1.23}) with $b$ replaced by $\widetilde{b}$ has a
unique small solution $(b^\ast, N^\ast)\in
\widetilde{\mathcal{M}}^{(0)}_{T}$.

In fact, let $\pi^\pm_{\widetilde{x}_{h_1}(v)}$ be the projections
to $\mathcal{H}_{\Sg_P}^\pm(\widetilde{x}_{h_1}(v))$ respectively
and consider the equation \bea
\pi^-_{\widetilde{x}_{h_1}(v)}\circ\pi_{v, -}^{0,
1}(((\widetilde{b}, N);\;v),\; \zeta)
\equiv\widetilde{\alpha}_T(v)+\pi^-_{\widetilde{x}_{h_1}(v)}\circ\epsilon(v,
\zeta)=0. \lb{5.1.28}\eea Note that by (ii-a) of Definition 1.1
and $s^{(2, -)}_{b,\,\widehat{x}}=D_{b,\,\widehat{x}}^{(2)}\psi_2\neq0$, the map
$\widetilde{\alpha}_T$ is transversal to the zero section. This
together with (\ref{5.1.24}) yields a unique solution $(b^\ast, N)$
of (\ref{5.1.28}) for each $N$ provided $|v|$ and $|N|$ are
sufficiently small.

Now consider the equation \bea
&&\pi^+_{\widetilde{x}_{h_1}(v)}\circ\pi_{v, -}^{0, 1}(((b^\ast,
N);\;v),\; \zeta) \equiv D^{(1)}_{h_1}(b^\ast,
N)s^{(1)}_{b^\ast,\,\widetilde{x}_{h_1}(v)}((d\phi_{b^\ast,\,
{\mathcal{T}(h_1)}}|_{\widetilde{x}_{h_1}(v)})^{-1}\rho_{h_1}(v))
\nn\\&&+\sum_{h_i\in\chi(T)\setminus h_1} (D^{(1)}_{h_i}(b^\ast,
N)s^{(1)}_{b^\ast,\,\widetilde{x}_{h_1}(v)}(\rho_{h_i}(v))
+D^{(1)}_{h_i}(b^\ast, 0)s^{(2,
+)}_{b^\ast,\,\widehat{x}}(x_i^\ast(v), \rho_{h_i}(v)))
\nn\\&&+\sum_{h_i\in\chi(T)} D^{(2)}_{h_i}(b^\ast, 0)s^{(2,
+)}_{b^\ast,\,\widehat{x}}(\rho_{h_i}(v))
+\pi^+_{\widetilde{x}_{h}(v)}\circ\epsilon(v, \zeta)=0
\lb{5.1.29}\eea By (ii-a) of Definition 1.1, the map \bea&&\Phi^+:
\widetilde{\mathcal{FT}}_{\delta_K}^{\emptyset} \rightarrow
T_{u(\Sg_P)}X\otimes\mathcal{H}^{+}_{\Sg_P}(\widetilde{x}_{h}(v)),
\nn\\&& N\mapsto D^{(1)}_{h_1}(b^\ast,
N)s^{(1)}_{b^\ast,\,\widetilde{x}_{h_1}(v)}((d\phi_{b^\ast,\,
{\mathcal{T}(h_1)}}|_{\widetilde{x}_{h_1}(v)})^{-1}\rho_{h_1}(v))
+\sum_{h_i\in\chi(T)\setminus h_1}D^{(1)}_{h_i}(b^\ast,
N)s^{(1)}_{b,\,\widetilde{x}_{h_1}(v)}(\rho_{h_i}(v)) \nn\eea is
transversal to the zero section. This together with (\ref{5.1.24})
yields a unique solution $(b^\ast, N^\ast)$ of (\ref{5.1.29})
provided $|v|$  is sufficiently small. Thus $(b^\ast, N^\ast)$ is
the unique solution of (\ref{5.1.23}) as desired. Denote by $\mu(v,
\zeta)=((b^\ast, N^\ast),\, v)$. Then the same argument in Theorem
5.1.1 yields the holomorphic map $\phi(v)=[\exp_{\mu(v,\, \zeta_v)
}\zeta_v]\in\mathfrak{M}^0_{2, k}(X, A, J)$.

For the general case, let $\widehat{h}=\min_{h<h_i,\,
h_i\in\chi(T)}\{h: |\{l\in I: \iota_l=h\}|\ge 2\}$. Then a similar
argument as above obtains the required conditions
$\alpha_T(v)=0=\widetilde{\alpha}_T(v)$. The proof of the theorem is
complete.\hfill\hb

{ \bf Theorem 5.1.4.} {\it Suppose $T$ is a bubble type given by
(iv) of Theorem 2.9 with $\Sg_P$ being smooth together with
$|I_1|=1$ and $|\chi(T)|\ge 2$. Assume the attaching node
$\widehat{x}$ of the bubble tree is one of the six branch points of
the canonical map. Then one of the following must hod:

(i) An element $b\equiv[\mathcal{C},\, u]\in
\mathfrak{M}_T(X, A, J)\cap\overline{\mathfrak{M}}^0_{2, k}(X, A, J)$
if and only if  $\alpha_T(v)=0$ for some $v=(b,
v)\in\widetilde{\mathcal{FT}}_{\delta_K}^\emptyset|_{\widetilde{K}^{(0)}}$
and $\widetilde{\alpha}_T(w, v)=0$, where $\widetilde{\alpha}_T(w,
v)$ is a linear combination of $\{D^{(1)}_{h_i}b,\,D_{h_i} ^{(2)}b,
\,D_{h_i} ^{(3)}b\}_{h_i\in\chi(T)}$ with coefficients depending on
$w\in T_{z_m}\Sg_P$, $v$ and the position of the nodes in bubble
trees. In particular,
$\mathfrak{M}_T(X, A, J)\cap\overline{\mathfrak{M}}^0_{2, k}(X, A, J)$
is a smooth orbifold of dimension at most $\dim
\overline{\mathfrak{M}}^{vir}_{2,k}(X, A)-2$.

(ii) An element $b\equiv[\mathcal{C},\, u]\in
\mathfrak{M}_T(X, A, J)\cap\overline{\mathfrak{M}}^0_{2, k}(X, A, J)$
must satisfy: there exists $h\in\chi(T)$ such that $u_h$ factor through a branched covering $\widetilde{u}: S^2\rightarrow
X$, i.e., there exists a holomorphic branched covering $\phi:
S^2\rightarrow S^2$ such that $u_h=\widetilde{u}\circ\phi$ and
$\deg(\phi)\ge 2$. In particular,
$\mathfrak{M}_T(X, A, J)\cap\overline{\mathfrak{M}}^0_{2, k}(X, A, J)$
is contained in a smooth orbifold of dimension at most $\dim
\overline{\mathfrak{M}}^{vir}_{2,k}(X, A)-2$.
}

{\bf Proof.} We only give the proof of the simplest case that all of
$\{h_1, h_2,\ldots,h_p\}\equiv\chi(T)$ are attached to the bubble
$\widehat{h}\in I_1$. The general cases follow similarly.

As in Theorem 5.1.3, we have (\ref{5.1.18}). Denote by
$\mathcal{S}_m=\{b\in \mathcal{S}: x_{\widehat{h}}=z_m\}$, where
$\mathcal{S}$ is given by Theorem 5.1.3, and $\mathcal{NS}_m$ the
normal bundle of $\mathcal{S}_m$ in $\mathcal{S}$,

Suppose $(b,\,w)\in \mathcal{NS}_m$, $(b, w, N)\in \mathcal{NS}$ and
$v=((b, w,
N);\;v)\in\widetilde{\mathcal{FT}}_{\delta_K}^\emptyset|_{\widetilde{K}^{(0)}}$.
We consider the third-order expansion of $\langle\langle
\overline{\partial}¡ª_Ju_{((b, w, N);\;v)},\, R_vV\eta\rangle\rangle_{v,
2}$. As in Theorem 5.1.2, we have
\bea&&\left|\frac{}{}\right.\langle\langle \overline{\partial}¡ª_Ju_{(b,
w, N);\;v},\, R_vV\eta\rangle\rangle_{v, 2}+ \langle\langle
\sum_{h_i\in\chi(T)}(D^{(1)}_{h_i}(b, w,
N)s^{(1)}_{b,\,\widetilde{x}_{h_i}(v)}((d\phi_{b,\,
{\mathcal{T}(h_i)}}|_{\widetilde{x}_{h_i}(v)})^{-1}\rho_{h_i}(v))
\nn\\&&\quad+D^{(2)}_{h_i}(b, w,
N)s^{(2)}_{b,\,\widetilde{x}_{h_i}(v)}((d\phi_{b,\,
{\mathcal{T}(h_i)}}|_{\widetilde{x}_{h_i}(v)})^{-1}\rho_{h_i}(v))
+\alpha^{(3)}_T((b, 0, 0);\;v), \,V\eta\rangle\rangle_{v,
2}\left|\frac{}{}\right. \nn\\\le&&
C_K|\rho(v)|^3(|v|+|w|+|N|)\|V\eta\|. \lb{5.1.30}\eea Let
$\widetilde{x}_{h_i}(v)=\widetilde{x}_{h_i}((b, w,
N);\;v)=\widetilde{x}_{h_i}((b, w,
0);\;v)\equiv\widetilde{x}_{h_i}(w, v)\in\Sg_P$ for $h_i\in\chi(T)$.
Identify a small neighborhood of $z_m$ in $\Sg_P$ with a small
neighborhood of $0$ in $T_{z_m}\Sg_P$ and let $x^\ast_i(w,
v)=\widetilde{x}_{h_i}(w, v)-\widetilde{x}_{h_1}(w, v)
=v_{\widehat{h}}(x_{h_i}-x_{h_1})$ and $x^\ast(w, v)=(x^\ast_2(w,
v),\ldots,x^\ast_p(w, v))$.

Let $\{\psi_j\}$ be an orthogonal basis for
$\mathcal{H}_{\Sg_P}^{0,\,1}$ such that
$\psi_1\in\mathcal{H}_{\Sg_P}^+(\widetilde{x}_{h_1}(w, v))$,
$\psi_2\in\mathcal{H}_{\Sg_P}^-(\widetilde{x}_{h_1}(w, v))$ and
$\{V_i\}$ an orthogonal basis for $T_{ev_P(b, w, N)}X$. Note that
since $\zeta\in \widetilde{\Gamma}^+(v)$, we have \bea
\langle\langle \zeta,\,D_v^\ast R_vV_i\psi_1\rangle\rangle_{v,\,2}=0
\lb{5.1.31}\eea Since
$\psi_2\in\mathcal{H}_{\Sg_P}^-(\widetilde{x}_{h_1}(w, v))$, we have
$\psi_2(\widetilde{x}_{h_1}(w, v))=0$. Note that  $s^{(2,
-)}(z_m)=0$, thus by (\ref{5.10}) we have (cf. Lemma 4.28 in
\cite{Z1}) \bea |D_v^\ast R_vV_i\psi_2|_{v,\,z}&\le&
C_K|\psi_2|_{(s,\,t)}|\partial_s| |V_i|\le
C_K(|\rho_{h_1}(v)|^2(|\rho_{h_1}(v)|+|\widetilde{x}_{h_1}(w, v)|).
\nn\\&+&\sum_{2\le i\le
p}|\rho_{h_i}(v)|(|v_{\widehat{h}}|^2+|\widetilde{x}_{h_1}(w,
v)||v_{\widehat{h}}|).  \lb{5.1.32}\eea Let
$\pi^\pm_{\widetilde{x}_{h_1}(w, v)}$ be the projections to
$\mathcal{H}_{\Sg_P}^\pm(\widetilde{x}_{h_1}(w, v))$ respectively.
Then by Taylor expansion with respect to $x_i^\ast(w,
v)=\widetilde{x}_{h_i}(w, v)-\widetilde{x}_{h_1}(w, v)$, we have
\bea &&\left|\frac{}{}\right. s^{(1)}_{b,\,\widetilde{x}_{h_i}(w,
v)}((d\phi_{b,\,
{\mathcal{T}(h_i)}}|_{\widetilde{x}_{h_i}(v)})^{-1}\rho_{h_i}(v))
-s^{(1)}_{b,\,\widetilde{x}_{h_1}(w, v)}(\rho_{h_i}(v))
-s^{(2)}_{b,\,\widetilde{x}_{h_1}(w, v)}(x_i^\ast(w, v),
\rho_{h_i}(v)) \nn\\&&\quad-s^{(3)}_{b,\,\widetilde{x}_{h_1}(w,
v)}(x_i^\ast(w, v), x_i^\ast(w, v), \rho_{h_i}(v))
\left|\frac{}{}\right. \le C_K|x_i^\ast(w, v)|^3|\rho_{h_i}(v)|.
\lb{5.1.33}\eea Since $\pi^-_{\widetilde{x}_{h_1}(w,
v)}s_{b,\,\widetilde{x}_{h_1}(w, v)}=0$, $s^{(2, -)}_{b, z_m}=0$ and
$\widetilde{x}_{h_i}(w, v)=w+v_{\widehat{h}}x_{h_i}$, we have \bea
&&\left|\pi^-_{\widetilde{x}_{h_1}(w,
v)}s^{(2)}_{b,\,\widetilde{x}_{h_1}(w, v)}(x_i^\ast(w, v),
\rho_{h_i}(v)) -s^{(3,\,-)}_{b,\,z_m}(\widetilde{x}_{h_1}(w, v),
x_i^\ast(w, v), \rho_{h_i}(v))\right| \nn\\\le
&&C_K|\widetilde{x}_{h_1}(w, v)|^2|\rho_{h_i}(v)||x_i^\ast(w, v)|.
\lb{5.1.34}\\
&&\left|\pi^-_{\widetilde{x}_{h_1}(w,
v)}s^{(3)}_{b,\,\widetilde{x}_{h_1}(w, v)}(x_i^\ast(w, v),
x_i^\ast(w, v), \rho_{h_i}(v)) -s^{(3,\,-)}_{b,\,z_m}(x_i^\ast(w,
v), x_i^\ast(w, v), \rho_{h_i}(v))\right| \nn\\\le
&&C_K|\widetilde{x}_{h_1}(w, v)||\rho_{h_i}(v)||x_i^\ast(w, v)|^2.
\lb{5.1.35} \eea Summing up (\ref{5.1.33})-(\ref{5.1.35}), we obtain
\bea &&\left|\frac{}{}\right.\pi^-_{\widetilde{x}_{h_1}(w,
v)}s^{(1)}_{b,\,\widetilde{x}_{h_i}(w, v)}((d\phi_{b,\,
{\mathcal{T}(h_i)}}|_{\widetilde{x}_{h_i}(v)})^{-1}\rho_{h_i}(v))
-s^{(3, -)}_{b,\,z_m}(\widetilde{x}_i(w, v), x_i^\ast(w, v),
\rho_{h_i}(v))
\left|\frac{}{}\right.\nn\\
\le &&C_K(|x_i^\ast(w, v)|^3|\rho_{h_i}(v)|+|\widetilde{x}_{h_1}(w,
v)||\rho_{h_i}(v)||x_i^\ast(w, v)|(|x_i^\ast(w,
v)|+|\widetilde{x}_{h_1}(w, v)|)). \lb{5.1.36} \eea Similarly, we
have \bea &&\left|\pi^-_{\widetilde{x}_{h_1}(w,
v)}s^{(2)}_{b,\,\widetilde{x}_{h_i}(v)}(d\phi_{b,
\widehat{h}}|_{\widetilde{x}_{h_i}(v)})^{-1}\rho_{h_i}(v))
-s^{(3,\,-)}_{b,\,z_m}(\widetilde{x}_{h_i}(w, v), \rho_{h_i}(v),
\rho_{h_i}(v))\right| \nn\\\le &&C_K(|x_i^\ast(w,
v)|^2|\rho_{h_i}(v)|^2+|\widetilde{x}_{h_1}(w,
v)|^2|\rho_{h_i}(v)|^2+|x_i^\ast(w, v)||\widetilde{x}_{h_1}(w, v)|
|\rho_{h_i}(v)|^2). \lb{5.1.37}\eea By
(\ref{5.1.30})-(\ref{5.1.37}), we have \bea &&\pi_{v, -}^{0, 1}(v,
\zeta)\equiv\pi_{v, -}^{0,
1}(\overline{\partial}_Ju_v+D_v\zeta+N_v\zeta) =D^{(1)}_{h_1}(b, w,
N)s^{(1)}_{b,\,\widetilde{x}_{h_1}(v)}((d\phi_{b,\,
{\mathcal{T}(h_1)}}|_{\widetilde{x}_{h_1}(v)})^{-1}\rho_{h_1}(v))
\nn\\&&+\sum_{h_i\in\chi(T)\setminus h_1} (D^{(1)}_{h_i}(b, w,
N)\pi^+_{\widetilde{x}_{h_1}(w,
v)}s^{(1)}_{b,\,\widetilde{x}_{h_i}(w, v)}((d\phi_{b,\,
{\mathcal{T}(h_i)}}|_{\widetilde{x}_{h_i}(v)})^{-1}\rho_{h_i}(v))
\nn\\&&¡¢\qquad+D^{(1)}_{h_i}(b, 0,
0)s^{(3,\,-)}_{b,\,z_m}(\widetilde{x}_{h_i}(w, v), x_i^\ast(w, v),
\rho_{h_i}(v))) \nn\\&&+\sum_{h_i\in\chi(T)} (D^{(2)}_{h_i}(b, w,
N)\pi^+_{\widetilde{x}_{h_1}(w,
v)}s^{(2)}_{b,\,\widetilde{x}_{h_i}(w, v)}((d\phi_{b,\,
{\mathcal{T}(h_i)}}|_{\widetilde{x}_{h_i}(v)})^{-1}\rho_{h_i}(v))
\nn\\&&\qquad+D^{(2)}_{h_i}(b, 0,
0)s^{(3,\,-)}_{b,\,z_m}(\widetilde{x}_{h_i}(w, v), \rho_{h_i}(v),
\rho_{h_i}(v))) +\alpha^{(3)}_T((b, 0, 0);\;v)+\epsilon(v, \zeta),
\lb{5.1.38}\eea where \bea &&\|\epsilon(v, \zeta)\|\le
C_K(|\rho(v)|^3+|x^\ast(w, v)||\rho(v)||\widetilde{x}(w, v)|+
|\rho(v)|^2|\widetilde{x}(w, v)|)(|v|+|w|+|N|)
\nn\\&&+C_K\left(|\rho_{h_1}(v)|^2(|\rho_{h_1}(v)|+|\widetilde{x}_{h_1}(w,
v)|) +\sum_{2\le i\le
p}|\rho_{h_i}(v)|(|v_{\widehat{h}}|^2+|\widetilde{x}_{h_1}(w,
v)||v_{\widehat{h}}|)\right) \cdot\|\zeta\|_{v,\,p,\,1}
\nn\\&&+C_K(|x_i^\ast(w,
v)|^3|\rho_{h_i}(v)|+|\widetilde{x}_{h_1}(w,
v)||\rho_{h_i}(v)||x_i^\ast(w, v)|(|x_i^\ast(w,
v)|+|\widetilde{x}_{h_1}(w, v)|))\nn\\&& +C_K(|x_i^\ast(w,
v)|^2|\rho_{h_i}(v)|^2+|\widetilde{x}_{h_1}(w,
v)|^2|\rho_{h_i}(v)|^2+|x_i^\ast(w, v)||\widetilde{x}_{h_1}(w, v)|
|\rho_{h_i}(v)|^2). \lb{5.1.39}\eea where $|\widetilde{x}(w,
v)|=\sum_{1\le i\le p}|\widetilde{x}_{h_i}(w, v)|$.

By replacing the order of $\{h_1, \ldots, h_p\}$ if necessary, we
may assume that $|\widetilde{x}_{h_1}(w,
v)|\le|\widetilde{x}_{h_i}(w, v)|$ for $2\le i\le p$. Note that
$\widetilde{x}_{h_i}(w, v)-\widetilde{x}_{h_j}(w, v)
=v_{\widehat{h}}(x_{h_i}-x_{h_j})$ and $x_{h_i}\neq x_{h_j}$ for
$i\neq j$. Thus we have $|v_{\widehat{h}}|\le
C_K|\widetilde{x}_{h_i}(w, v)|$ for $2\le i\le p$. Hence we have
\bea |s^{(3,\,-)}_{b,\,z_m}(\widetilde{x}_{h_i}(w, v),
\rho_{h_i}(v), \rho_{h_i}(v))|
=o(|s^{(3,\,-)}_{b,\,z_m}(\widetilde{x}_{h_i}(w, v),
x_i^\ast(w, v), \rho_{h_i}(v))|) \nn\\
|s^{(3,\,-)}_{b,\,z_m}(\rho_{h_i}(v), \rho_{h_i}(v), \rho_{h_i}(v))|
=o(|s^{(3,\,-)}_{b,\,z_m}(\widetilde{x}_{h_i}(w, v), x_i^\ast(w, v),
\rho_{h_i}(v))|) \nn\eea for $2\le i\le p$.

We have the following two cases:

{\bf Case 1.} {\it We have $D^{(1)}_{h_1}b\neq0$.}

In this case we can represent $v_{h_1}$ as linear combination of
$v_{h_i}$s for $h_i\in\chi(T)\setminus h_1$. In particular, we have
$|v_{h_1}|\le C_K\sum_{h_i\in\chi(T)\setminus h_1}|v_{h_i}|$. Hence
we have \bea |s^{(3,\,-)}_{b,\,z_m}(\widetilde{x}_{h_1}(w, v),
\rho_{h_1}(v), \rho_{h_1}(v))|=\sum_{h_i\in\chi(T)\setminus
h_1}o(|s^{(3,\,-)}_{b,\,z_m}(\widetilde{x}_{h_i}(w, v),
x_i^\ast(w, v), \rho_{h_i}(v))|) \nn\\
|s^{(3,\,-)}_{b,\,z_m}(\rho_{h_1}(v))| =\sum_{h_i\in\chi(T)\setminus
h_1}o(|s^{(3,\,-)}_{b,\,z_m}(\widetilde{x}_{h_i}(w, v), x_i^\ast(w,
v), \rho_{h_i}(v))|) \nn\eea Note that $s^{(3,\,-)}_{b,\,z_m}\neq
0$, thus in order to satisfy (\ref{5.21}), we must have \bea
\widetilde{\alpha}_T(w, v)\equiv\sum_{2\le i\le p} D^{(1)}_{h_i}(b,
0, 0)s^{(3,\,-)}_{b,\,z_m}(\widetilde{x}_{h_i}(w, v), x_i^\ast(w,
v), \rho_{h_i}(v))=0 \lb{5.1.40}\eea by (\ref{5.1.38}) and
(\ref{5.1.39}).

{\bf Case 2.} {\it We have $D^{(1)}_{h_1}b=0$.}

In this case we have \bea &&\widetilde{\alpha}_T(w,
v)\equiv\sum_{2\le i\le p} D^{(1)}_{h_i}(b, 0,
0)s^{(3,\,-)}_{b,\,z_m}(\widetilde{x}_{h_i}(w, v), x_i^\ast(w, v),
\rho_{h_i}(v)) \nn\\&&+D^{(2)}_{h_1}(b, 0,
0)s^{(3,\,-)}_{b,\,z_m}(\widetilde{x}_{h_1}(w, v), \rho_{h_1}(v),
\rho_{h_1}(v)) +D^{(3)}_{h_1}(b, 0,
0)s^{(3,\,-)}_{b,\,z_m}(\rho_{h_1}(v))=0 \lb{5.1.41}\eea provided
$|v|$, $|w|$ and $|N|$ are sufficiently small.

Note that $\alpha_T(v)$ and $\widetilde{\alpha}_T(w, v)$ are linear
independent in Case 1. While in Case 2, we  have
\bea\alpha_T(v)=\sum_{h_i\in\chi(T)\setminus
h_1}(D^{(1)}_{h_i}b)s_{\Sg_P,\,\widehat{x}}(\rho_{h_i}(v))=0.
\nn\eea Thus the dimension of $\alpha_T^{-1}(0)$ is at most $\dim
\overline{\mathfrak{M}}^{vir}_{2,k}(X, A)-2$

The proof of the converse is similar to the previous theorems
and we omit it here. We only remark that in Case 1, (i) holds;
while in Case 2, (i) holds if
 $H^1(\P^1, u_{h_1}^\ast TX\otimes \mathcal{O}_{\P^1}(-4\infty))=0$
 and (ii) holds if  $H^1(\P^1, u_{h_1}^\ast TX\otimes \mathcal{O}_{\P^1}(-4\infty))\neq0$.
These follows by (ii-d) in Definition 1.1. \hfill\hb

\subsection{ There are exactly two bubble trees}

In this subsection, we consider the case $|I_1|=2$.
Note that in this case the rank of $\alpha_T(v)$
may not equal to the dimension of $\mathrm{coker}D_b$.
Thus we need further expansions.

{ \bf Theorem 5.2.1.} {\it Suppose $T$ is a bubble type given by
(iv) of Theorem 2.9 with $\Sg_P$ being smooth and $|I_1|=2$, i.e.,
there are exactly two bubble trees. Assume  $x_1$ and $x_2$ are not
conjugate via the map $s_{\Sg_P}$, where $x_1$ and $x_2$ are
attaching nodes of the two bubble trees. Then an element
$b\equiv[\mathcal{C},\, u]\in \mathfrak{M}_T(X, A, J)\cap\overline{\mathfrak{M}}^0_{2, k}(X, A, J)$
if and only if $\alpha_T(v)=0$ for some $v=(b,
v)\in\widetilde{\mathcal{FT}}_{\delta_K}^\emptyset|_{\widetilde{K}^{(0)}}$.
In particular, $\mathfrak{M}_T(X, A,
J)\cap\overline{\mathfrak{M}}^0_{2, k}(X, A, J)$ is  a smooth
orbifold of dimension at most $\dim
\overline{\mathfrak{M}}^{vir}_{2,k}(X, A)-2$.  }

{\bf Proof.} By Lemma 5.5, an element $b\equiv[\mathcal{C},\, u]\in
\mathfrak{M}_T(X, A, J)\cap\overline{\mathfrak{M}}^0_{2, k}(X, A, J)$
must satisfy $\alpha_T(v)=0$ for some $v=(b,
v)\in\widetilde{\mathcal{FT}}_{\delta_K}^\emptyset|_{\widetilde{K}^{(0)}}$.
Since $x_1$ and $x_2$ are not conjugate via the map $s_{\Sg_P}$, we
have $\mathrm{rank}(\alpha_T)=4n=\dim(\mathrm{coker}D_b)$. Hence
\bea\alpha_T: \widetilde{\mathcal{FT}}_{\delta_K}^{\emptyset}
\rightarrow T_{u(\Sg_P)}X\otimes\mathcal{H}^{0, 1}_{\Sg_P}, \qquad
v\mapsto \alpha_T(v),\nn\eea
 is transversal to the zero section
in a small neighborhood of  $b$ by (ii-a) in Definition 1.1. Hence by (\ref{5.23})
and (\ref{5.24}), in a small neighborhood of  $b$, for any \bea
v=(\widetilde{b},\,
v)\in\widetilde{\mathcal{FT}}_{\delta_K}^{\emptyset}|_{\widetilde{K}^{(0)}},
, \quad \|\zeta\|_{v, p, 1}\le 2C_K|v|^\frac{1}{p}, \nn\eea the
equation \bea \pi_{v, -}^{0, 1}((\widetilde{b},\,v),
\,\zeta)\equiv\alpha_T(v)+\epsilon(v, \zeta)=0 \nn\eea has a unique
small solution $b^\ast\in \widetilde{\mathcal{M}}^{(0)}_{T}$. Denote
by $\mu(v, \zeta)=(b^\ast,\, v)$. Then the theorem follows by the
same argument as in Theorem 5.1.1.
 \hfill\hb

Next we consider the case $x_1$ and $x_2$ are conjugate
via the map $s_{\Sg_P}$.

{ \bf Theorem 5.2.2.} {\it Suppose $T$ is a bubble type given by (iv) of
Theorem 2.9 with $\Sg_P$ being smooth and $|I_1|=|\chi(T)|\equiv|\{h_1, h_2\}|=2$.
Assume the attaching nodes $x_{\mathcal{T}(h_1)}$
and $x_{\mathcal{T}(h_2)}$
differ by the nontrivial holomorphic automorphism of $\Sg_P$
(cf. p.254 of \cite{GH}).
Then one of the following must hold:

(i) An element $b\equiv[\mathcal{C},\, u]\in
\mathfrak{M}_T(X, A, J)\cap\overline{\mathfrak{M}}^0_{2, k}(X, A, J)$
if and only if $\alpha_T(v)=0$ for some $v=(b,
v)\in\widetilde{\mathcal{FT}}_{\delta_K}^\emptyset|_{\widetilde{K}^{(0)}}$
and $\widetilde{\alpha}_T(w, v)=0$, where $\widetilde{\alpha}_T(w, v)$ is
a linear combination of $\{D^{(1)}_{h_i}b, D_{h_i} ^{(2)}b\}_{i=1, 2}$ with
coefficients depending on $w\in T_{x_{\mathcal{T}(h_1)}}\Sg_P$, $v$ and the position of the nodes in bubble trees.
In particular,
$\mathfrak{M}_T(X, A, J)\cap\overline{\mathfrak{M}}^0_{2, k}(X, A, J)$
is a smooth orbifold of dimension at most $\dim
\overline{\mathfrak{M}}^{vir}_{2,k}(X, A)-2$
.

(ii) An element $b\equiv[\mathcal{C},\, u]\in
\mathfrak{M}_T(X, A, J)\cap\overline{\mathfrak{M}}^0_{2, k}(X, A, J)$
must satisfy $u_1(\Sg_{h_1})=u_2(\Sg_{h_2})$.
}

{\bf Proof.} By Lemma 5.5, an element $b\equiv[\mathcal{C},\, u]\in
\mathfrak{M}_T(X, A, J)\cap\overline{\mathfrak{M}}^0_{2, k}(X, A, J)$
must satisfy $\alpha_T(v)=0$ for some $v=(b, v)\in\widetilde{\mathcal{FT}}_{\delta_K}^\emptyset|_{\widetilde{K}^{(0)}}$.
Note that
\be \alpha_T(v)=(D^{(1)}_{h_1}b)s_{\Sg_P,\,x_{\mathcal{T}(h_1)}}(\rho_{h_1}(v))
+(D^{(1)}_{h_2}b)s_{\Sg_P,\,x_{\mathcal{T}(h_2)}}(\rho_{h_2}(v)),
\lb{5.2.1}\ee
and $s_{\Sg_P,\,x_{\mathcal{T}(h_1)}}=s_{\Sg_P,\,x_{\mathcal{T}(h_2)}}$.
Denote by  $-x$ the image of $x\in\Sg_P$ under the
nontrivial holomorphic automorphism $\sg$ of $\Sg_P$
and identify $T_x\Sg_P$ with $T_{-x}\Sg_P$ via $\sg$.
Let $\Sg_P^\ast=\Sg_P\setminus\{z_m\}_{1\le m\le 6}$,
where $\{z_m\}_{1\le m\le 6}$ are the fixed points of
$\sg$. Denote by $\mathcal{S}=\alpha_T^{-1}(0)$.
Then by (ii-a) in Definition 1.1, $\mathcal{S}$ is a complex suborbifold
of $\mathfrak{M}_T(X, A, J)$. Moreover, we have
$\mathcal{S}=\mathcal{S}_0\times \mathcal{S}_1$,
where $\mathcal{S}_0=\{(x_{\mathcal{T}(h_1)},\,-x_{\mathcal{T}(h_1)}: x_{\mathcal{T}(h_1)}\in\Sg_P^\ast\}$
and $\mathcal{S}_1=\{D^{(1)}_{h_2}b=\lambda D^{(1)}_{h_1}b,\,\lambda\in\C\}$.
Let $\mathcal{NS}$ denote the normal bundle of
$\mathcal{S}$ in $\mathfrak{M}_T(X, A, J)$
and identify a small neighborhood
of its zero section with a tubular neighborhood of
$\mathcal{S}$ in $\mathfrak{M}_T(X, A ,J)$.
Then we have
\bea\mathcal{NS}=\mathcal{NS}_0\oplus\mathcal{NS}_1,
\quad{\rm where}\quad
\mathcal{NS}_0=\pi_{\Sg_P,\,x_{\mathcal{T}(h_2)}}^\ast T\Sg_P,\quad
\mathcal{NS}_1= E_1,
\nn\eea
where $\pi_{\Sg_P,\,x_{\mathcal{T}(h_i)}}: \mathcal{S}_0\subset
\Sg_P\times\Sg_P\rightarrow\Sg_P$ is the projection to the $i$-th factor
and $E_1$ is the orthogonal complement of $D^{(1)}_{h_1}b$
in $T_{u(\Sg_P)}X$.

Suppose $(b, w, N)\in \mathcal{NS}$ and
$v=((b, w, N);\;v)\in\widetilde{\mathcal{FT}}_{\delta_K}^\emptyset|_{\widetilde{K}^{(0)}}$.
We consider the second-order expansion of
$\langle\langle \overline{\partial}¡ª_Ju_{((b, w, N);\;v)},\, R_vV\eta\rangle\rangle_{v, 2}$.
Note that we have
\bea ||s^{(2)}_{b,\,\widetilde{x}_h(v)}((d\phi_{b,\, {\mathcal{T}(h)}}|_{\widetilde{x}_h(v)})^{-1}\rho_h(v))
-s^{(2)}_{b,\,x_{\mathcal{T}(h)}}(\rho_h(v))||_2&\le&
C_K|\phi_{b,\, \mathcal{T}(h)}(\widetilde{x}_h(v))|_{b}|\rho_h(v)|^2
\nn\\&\le& C_K|v|\cdot|\rho_h(v)|^2,\nn\\
\sum_{m\ge 3}\left|D_h^{(m)}(b, w, N)\right||\rho_h(v)|^m&\le& C_K|\rho_h(v)|^3,
\nn\\\left|D_h^{(2)}(b, w, N)-D_h^{(2)}(b, 0, 0)\right|&\le& C_K(|N|+|w|),
\nn\eea
for  $h\in\chi(T)$ and $v,\,w,\, N$ sufficiently small
by continuity. Thus by Lemma 5.1 we have
\bea&&\left|\frac{}{}\right.\langle\langle \overline{\partial}¡ª_Ju_{(b, w, N);\;v},\, R_vV\eta\rangle\rangle_{v, 2}+
\langle\langle D^{(1)}_{h_1}(b, w, N)s^{(1)}_{b,\,\widetilde{x}_{h_1}(v)}((d\phi_{b,\, {\mathcal{T}(h_1)}}|_{\widetilde{x}_{h_1}(v)})^{-1}\rho_{h_1}(v))
\nn\\&&\qquad+D^{(1)}_{h_2}(b, w, N)s^{(1)}_{b,\,\widetilde{x}_{h_2}(v)}((d\phi_{b,\, {\mathcal{T}(h_2)}}|_{\widetilde{x}_{h_2}(v)})^{-1}\rho_{h_2}(v))
+\alpha^{(2)}_T((b, 0, 0);\;v), \,V\eta\rangle\rangle_{v, 2}\left|\frac{}{}\right.
\nn\\\le&&C_K|\rho(v)|^2(|v|+|w|+|N|)\|V\eta\|.
\lb{5.2.2}\eea
Let
$\widetilde{x}_{h_i}((b, w, N);\;v)=\widetilde{x}_{h_i}((b, w, 0);\;v)\equiv\widetilde{x}_{h_i}(w, v)\in\Sg_P$
for $i=1,\,2$.
Identify a small neighborhood of $x_{\mathcal{T}(h_i)}$ in $\Sg_P$ with a
small neighborhood of $0$ in $T_{x_{\mathcal{T}(h_i)}}\Sg_P$
and let $x^\ast(w, v)=\widetilde{x}_{h_2}(w, v)-\widetilde{x}_{h_1}(w, v)$,
then we have $|x^\ast(w, v)|\le C_K(|v|+|w|)$.

Let $\{\psi_j\}$ be an orthonormal basis for $\mathcal{H}_{\Sg_P}^{0,\,1}$
such that $\psi_1\in\mathcal{H}_{\Sg_P}^+(\widetilde{x}_{h_1}(w, v))$,
$\psi_2\in\mathcal{H}_{\Sg_P}^-(\widetilde{x}_{h_1}(w, v))$
and $\{V_i\}$ an orthonormal basis for $T_{ev_P(b, w, N)}X$.
Note that since $\zeta\in \widetilde{\Gamma}^+(v)$, we have
\bea \langle\langle \zeta,\,D_v^\ast R_vV_i\psi_1\rangle\rangle_{v,\,2}=0
\lb{5.2.3}\eea
Since $\psi_2\in\mathcal{H}_{\Sg_P}^-(\widetilde{x}_{h_1}(w, v))$,
we have $\psi_2(\widetilde{x}_{h_1}(w, v))=0$.
Note that $\psi_2(-\widetilde{x}_{h_1}(w, v))=0$
since $\psi_2$ is invariant under $\sg$.
Hence by (\ref{5.10}) we have
\bea |D_v^\ast R_vV_i\psi_2|_{v,\,z}\le C_K|\psi_2|_{(s,\,t)}|\partial_s|
|V_i|\le C_K|\rho(v)|(|\rho(v)|+|x^\ast(w, v)|).
\lb{5.2.4}\eea
Note that  we have
\bea &&\left|s^{(1)}_{b,\,\widetilde{x}_{h_2}(w, v)}
((d\phi_{b,\, {\mathcal{T}(h_2)}}|_{\widetilde{x}_{h_2}(v)})^{-1}\rho_{h_2}(v))
-s^{(1)}_{b,\,\widetilde{x}_{h_1}(w, v)}(\rho_{h_2}(v))
-s^{(2)}_{b,\,x_{\mathcal{T}(h_1)}}(x^\ast(w, v), \rho_{h_2}(v))
\right|\nn\\\le &&C_K(|x^\ast(w, v)||\rho(v)|(|x^\ast(w, v)|+|v|).
\lb{5.2.5}\eea
By (\ref{5.2.2})-(\ref{5.2.5}), we have
\bea &&\pi_{v, -}^{0, 1}(v, \zeta)\equiv\pi_{v, -}^{0, 1}(\overline{\partial}_Ju_v+D_v\zeta+N_v\zeta)
\nn\\=&&D^{(1)}_{h_1}(b, w, N)s^{(1)}_{b,\,\widetilde{x}_{h_1}(w, v)}((d\phi_{b,\, {\mathcal{T}(h_1)}}|_{\widetilde{x}_{h_1}(v)})^{-1}\rho_{h_1}(v))
+D^{(1)}_{h_2}(b, w, N)s^{(1)}_{b,\,\widetilde{x}_{h_1}(w, v)}(\rho_{h_2}(v))
\nn\\&&+D^{(1)}_{h_2}(b, 0, 0)s^{(2)}_{b,\,x_{\mathcal{T}(h_1)}}(x^\ast(w, v), \rho_{h_2}(v))
+D^{(2)}_{h_1}(b, 0, 0)s^{(2,\,+)}_{b,\,x_{\mathcal{T}(h_1)}}(\rho_{h_1}(v))
\nn\\&&+D^{(2)}_{h_1}(b, 0, 0)s^{(2, -)}_{b,\,x_{\mathcal{T}(h_1)}}(\rho_{h_1}(v))
+D^{(2)}_{h_2}(b, 0, 0)s^{(2,\,+)}_{b,\,x_{\mathcal{T}(h_2)}}(\rho_{h_2}(v))
\nn\\&&+D^{(2)}_{h_2}(b, 0, 0)s^{(2, -)}_{b,\,x_{\mathcal{T}(h_2)}}(\rho_{h_2}(v))
+\epsilon(v, \zeta)
\lb{5.2.6}\eea
where
\bea &&\|\epsilon(v, \zeta)\|\le
C_K|\rho(v)|^2(|v|+|w|+|N|)
+C_K|\rho(v)|(|\rho(v)|+|x^\ast(w, v)|)\cdot\|\zeta\|_{v,\,p,\,1}
\nn\\&&+C_K(|x^\ast(w, v)||\rho(v)|(|x^\ast(w, v)|+|v|+|w|+|N|).\lb{5.2.7}\eea
Since $x_{\mathcal{T}(h_1)}\in\Sg_P^\ast$, we have $s^{(2, -)}_{b,\,x_{\mathcal{T}(h_1)}}\neq 0$.
Hence in order to satisfies (\ref{5.21}),
we must have
\bea &&\widetilde{\alpha}_T(w, v)\equiv D^{(1)}_{h_2}(b, 0, 0)s^{(2,-)}_{b,\,x_{\mathcal{T}(h_1)}}(x^\ast(w, v), \rho_{h_2}(v))
+D^{(2)}_{h_1}(b, 0, 0)s^{(2, -)}_{b,\,x_{\mathcal{T}(h_1)}}(\rho_{h_1}(v))
\nn\\&&+D^{(2)}_{h_2}(b, 0, 0)s^{(2, -)}_{b,\,x_{\mathcal{T}(h_2)}}(\rho_{h_2}(v))=0
\lb{5.2.8}\eea
provided $|v|$, $|w|$ and $|N|$ are sufficiently small.

Conversely, suppose $\alpha_T(v)=0$
and $\widetilde{\alpha}_T(w, v)=0$.
We want to construct  $\phi(w, v)\in\mathfrak{M}^0_{2, k}(X, A, J)$
converging to $b$. We want to show: In a small neighborhood of $b$, for any
\bea v=((\widetilde{b}, w, N);\;v)\in\widetilde{\mathcal{FT}}_{\delta_K}^{\emptyset}|_{\widetilde{K}^{(0)}},
\quad |w|\le\delta_K,\quad \|\zeta\|_{v, p, 1}\le 2C_K|v|^\frac{1}{p}
\nn\eea
the equation (\ref{5.2.6}) with $b$ replaced by $\widetilde{b}$ has
a unique small solution $(b^\ast, w, N^\ast)\in
\widetilde{\mathcal{M}}^{(0)}_{T}$.
Let $\pi^\pm_{\widetilde{x}_{h_1}(w, v)}$ be the projections
to $\mathcal{H}_{\Sg_P}^\pm(\widetilde{x}_{h_1}(w, v))$ respectively
and consider the following equation:
\bea
&&\pi^-_{\widetilde{x}_{h_1}(w, v)}\circ\pi_{v, -}^{0,
1}(((\widetilde{b}, w, N);\;v),\; \zeta) \equiv
D^{(1)}_{h_2}(\widetilde{b}, 0,
0)s^{(2,-)}_{b,\,x_{\mathcal{T}(h_1)}}(x^\ast(w, v), \rho_{h_2}(v))
\nn\\&&+D^{(2)}_{h_1}(\widetilde{b}, 0, 0)s^{(2,
-)}_{\widetilde{b},\,x_{\mathcal{T}(h_1)}}(\rho_{h_1}(v))
+D^{(2)}_{h_2}(\widetilde{b}, 0, 0)s^{(2,
-)}_{\widetilde{b},\,x_{\mathcal{T}(h_2)}}(\rho_{h_2}(v))
+\pi^-_{\widetilde{x}_{h_1}(w, v)}\circ\epsilon(v, \zeta)=0.\qquad
\lb{5.2.9}\eea
We have the following two cases:

{\bf Case 1.} {\it We have $D^{(1)}_{h_1}b\neq0$.}

In this case we have $|\rho_{h_1}(v)|\le C_K|\rho_{h_2}(v)|$.
Note that by $\alpha_T(v)=0$ and (v-a) in Definition 1.1,
if $u_{h_1}(\Sg_{h_1})\neq u_{h_2}(\Sg_{h_2})$, we may assume
$H^1(\P^1, u_2^\ast TX\otimes \mathcal{O}_{\P^1}(-3z_2))\neq0$,
thus $\widetilde{\alpha}_T(w, v)$ is transversal to the zero section.
By (\ref{5.2.7}), the terms in $\pi^-_{\widetilde{x}_{h_1}(w, v)}\circ\epsilon(v, \zeta)$
contains $\rho_{h_1}(v)$ can be dominated by $|\rho_{h_2}(v)|^2$.
Thus we obtain a unique solution $(b^\ast, w, N)$ of (\ref{5.2.9}) for each
$N$ when $\delta_K$ is sufficiently small.

{\bf Case 2.} {\it We have $D^{(1)}_{h_1}b=0$.}

In this case, we have $D^{(1)}_{h_2}b=0$ by (\ref{5.2.1}).
Hence by (ii-b) in Definition 1.1, $\widetilde{\alpha}_T(w, v)$
is transversal to the zero section.
This together with (\ref{5.2.7})
yields a unique solution $(b^\ast, w, N)$ of (\ref{5.2.9}) for each
$N$ when $\delta_K$ is sufficiently small.

Now consider the equation \bea &&\pi^+_{\widetilde{x}_{h_1}(w,
v)}\circ\pi_{v, -}^{0, 1}(((b^\ast, w, N);\;v),\; \zeta) \equiv
D^{(1)}_{h_1}(b^\ast, w, N)s^{(1)}_{b^\ast,\,\widetilde{x}_{h_1}(w,
v)}((d\phi_{b^\ast,\,{\mathcal{T}(h_1)}}|_{\widetilde{x}_{h_1}(v)})^{-1}\rho_{h_1}(v))
\nn\\&&+D^{(1)}_{h_2}(b^\ast, w,
N)s^{(1)}_{b^\ast,\,\widetilde{x}_{h_1}(w, v)}(\rho_{h_2}(v))
+D^{(1)}_{h_2}(b^\ast, 0, 0)s^{(2,
+)}_{b^\ast,\,x_{\mathcal{T}(h_1)}}(x^\ast(w, v), \rho_{h_2}(v))
\nn\\&&+D^{(2)}_{h_1}(b^\ast, 0, 0)s^{(2,
+)}_{b^\ast,\,x_{\mathcal{T}(h_1)}}(\rho_{h_1}(v)) +D^{(2)}_{h_2}(b^\ast, 0,
0)s^{(2,\,+)}_{b^\ast,\,x_{\mathcal{T}(h_2)}}(\rho_{h_2}(v))\nn\\
&&+\pi^+_{\widetilde{x}_{h_1}(w, v)}\circ\epsilon(v, \zeta)=0.
\lb{5.2.10}\eea By (ii-a) in Definition 1.1, the map \bea &&\Phi^+:
\widetilde{\mathcal{FT}}_{\delta_K}^{\emptyset} \rightarrow
T_{u(\Sg_P)}X\otimes\mathcal{H}^{+}_{\Sg_P}(\widetilde{x}_{h}(v)),
\nn\\&& \quad N\mapsto D^{(1)}_{h_1}(b^\ast, w,
N)s^{(1)}_{b^\ast,\,\widetilde{x}_{h_1}(w, v)}((d\phi_{b,\,
{\mathcal{T}(h_1)}}|_{\widetilde{x}_{h_1}(v)})^{-1}\rho_{h_1}(v))
+D^{(1)}_{h_2}(b^\ast, w,
N)s^{(1)}_{b^\ast,\,\widetilde{x}_{h_1}(w, v)}(\rho_{h_2}(v))
\nn\eea
is transversal to the zero section. Note that
\bea &&\left\|\frac{}{}\right.
D^{(1)}_{h_2}(b^\ast, 0, 0)s^{(2, +)}_{b^\ast,\,x_{\mathcal{T}(h_1)}}(x^\ast(w, v), \rho_{h_2}(v))
+D^{(2)}_{h_1}(b^\ast, 0, 0)s^{(2, +)}_{b^\ast,\,x_{\mathcal{T}(h_1)}}(\rho_{h_1}(v))
\nn\\&&\quad+D^{(2)}_{h_2}(b, 0, 0)s^{(2,\,+)}_{b^\ast,\,x_{\mathcal{T}(h_2)}}(\rho_{h_2}(v))\left\|\frac{}{}\right.
\le C_K|\rho(v)|(|\rho(v)|+|x^\ast(w, v)|). \nn\eea This together with (\ref{5.2.7}) yields a
unique solution $(b^\ast, w, N^\ast)$ of (\ref{5.2.10}) provided
$|v|$ is sufficiently small. Thus $(b^\ast, w, N^\ast)$ is the unique
solution of (\ref{5.2.6}) as desired. Let $\mu(v, \zeta)=((b^\ast,
w, N^\ast),\, v)$. Then the same argument in Theorem 5.1.1 yields
the holomorphic map $\phi(w, v)=[\exp_{\mu(v,\, \zeta_v)
}\zeta_v]\in\mathfrak{M}^0_{2, k}(X, A, J)$. The proof of the
theorem is complete.\hfill\hb

Now we consider the general case.

{ \bf Theorem 5.2.3.} {\it Suppose $T$ is a bubble type given by
(iv) of Theorem 2.9 with $\Sg_P$ being smooth together with
$|I_1|=2$ and $|\chi(T)|\ge 3$. Assume the attaching nodes $x_{\hat
h_1}$ and $x_{\hat h_2}$ of the two bubble trees differ by the
nontrivial holomorphic automorphism of $\Sg_P$. Then an element
$b\equiv[\mathcal{C},\, u]\in \mathfrak{M}_T(X, A, J)\cap\overline{\mathfrak{M}}^0_{2, k}(X, A, J)$
if and only if $\alpha_T(v)=0$ for some $v=(b, v)\in\widetilde{\mathcal{FT}}_{\delta_K}^\emptyset|_{\widetilde{K}^{(0)}}$
and $\widetilde{\alpha}_T(w, v)=0$, where $\widetilde{\alpha}_T(w,
v)$ is a linear combination of $\{D^{(1)}_{h_i}b, D_{h_i}
^{(2)}b\}_{h_i\in\chi(T)}$ with coefficients depending on $w\in
T_{x_{\hat h_1}}\Sg_P$, $v$ and the position of the nodes in bubble
trees. In particular, $\mathfrak{M}_T(X, A, J)\cap\overline{\mathfrak{M}}^0_{2, k}(X, A, J)$
is a smooth orbifold of dimension at most $\dim
\overline{\mathfrak{M}}^{vir}_{2,k}(X, A)-2$. }

{\bf Proof.} Denote by $\chi(T)=\{h_1, h_2,\ldots,h_p\}$. We only
give the proof of the simplest case that all of $\{h_1,
h_2,\ldots,h_{q}\}$ are attached to the bubble $h_\alpha\in
I_1$ and $\{h_{q+1},\ldots, h_p\}$ are attached to $h_\beta\in
I_1$. The
general case follows similarly.

By Lemma 5.5, an element $b\equiv[\mathcal{C},\, u]\in
\mathfrak{M}_T(X, A, J)\cap\overline{\mathfrak{M}}^0_{2, k}(X, A, J)$
must satisfy $\alpha_T(v)=0$ for some $v=(b,
v)\in\widetilde{\mathcal{FT}}_{\delta_K}^\emptyset|_{\widetilde{K}^{(0)}}$,
where
\bea\alpha_T(v)=\sum_{h_i\in\chi(T)}(D^{(1)}_{h_i}b)s_{\Sg_P,\,\widehat{x}}(\rho_{h_i}(v)).
\lb{5.2.11}\eea Here we identify $x_{\hat h_1}$ and $x_{\hat h_2}$
via the nontrivial holomorphic automorphism $\sg$ of $\Sg_P$ and
denote them simply by $\widehat{x}$. Denote by $\mathcal{S}=\{b\in
\mathfrak{M}_T(X, A, J): \alpha_T(v)=0\}$. Then by (ii-a) in Definition 1.1,
$\mathcal{S}$ is complex suborbifold of $\mathfrak{M}_T(X, A, J)$.
Let $\mathcal{NS}$ denote the normal bundle of $\mathcal{S}$ in
$\mathfrak{M}_T(X, A, J)$. As in Theorem 5.2.2, we have
$\mathcal{S}=\mathcal{S}_0\times \mathcal{S}_1$, where
\bea\mathcal{S}_0=\{(x_{\hat h_1},\,-x_{\hat h_1}): x_{\hat
h_1}\in\Sg_P^\ast\},
\qquad\mathcal{S}_1=\left\{\sum_{h_i\in\chi(T)}(D^{(1)}_{h_i}b)s_{\Sg_P,\,x_{\hat
h_1}}(\rho_{h_i}(v))=0\right\}. \nn\eea Then we have
\bea\mathcal{NS}=\mathcal{NS}_0\oplus\mathcal{NS}_1, \quad{\rm
where}\quad \mathcal{NS}_0=\pi_{\Sg_P,\,x_{\mathcal{T}(h_p)}}^\ast
T\Sg_P,\quad \mathcal{NS}_1= E_1, \nn\eea where
$\pi_{\Sg_P,\,x_{\mathcal{T}(h_i)}}: \mathcal{S}_0\subset
\Sg_P\times\Sg_P\rightarrow\Sg_P$ is the projection to the $i$-th
factor and $E_1$ is the orthogonal complement of
$\mathrm{span}_{h_i\in\chi(T)}\{D^{(1)}_{h_i}b\}$ in
$T_{u(\Sg_P)}X$.

Suppose $(b, w, N)\in \mathcal{NS}$ and $v=((b, w,
N);\;v)\in\widetilde{\mathcal{FT}}_{\delta_K}^\emptyset|_{\widetilde{K}^{(0)}}$.
Let $\widetilde{x}_{h_i}((b, w, N);\;v)=\widetilde{x}_{h_i}((b, w,
0);\;v)\equiv\widetilde{x}_{h_i}(w, v)\in\Sg_P$ for $h_i\in\chi(T)$.
Identify a small neighborhood of $x_{\mathcal{T}(h_i)}$ in $\Sg_P$
with a small neighborhood of $0$ in $T_{x_{\mathcal{T}(h_i)}}\Sg_P$
and let $x_i^\ast(w, v)=\widetilde{x}_{h_i}(w,
v)-\widetilde{x}_{h_1}(w, v)$, then we have $|x^\ast(w, v)|\le
C_K(|v|+|w|)$.

We consider the second-order expansion of $\langle\langle
\overline{\partial}¡ª_Ju_{((b, w, N);\;v)},\, R_vV\eta\rangle\rangle_{v,
2}$. Similar to Theorem 5.1.3 and Theorem 5.2.2, we have \bea
\pi_{v, -}^{0, 1}(v, \zeta)&\equiv&\pi_{v, -}^{0,
1}(\overline{\partial}_Ju_v+D_v\zeta+N_v\zeta) =D^{(1)}_{h_1}(b, w,
N)s^{(1)}_{b,\,\widetilde{x}_{h_1}(w, v)}((d\phi_{b,\,
{\mathcal{T}(h_1)}}|_{\widetilde{x}_{h_1}(w, v)})^{-1}\rho_{h_1}(v))
\nn\\&&+\sum_{2\le i\le p} (D^{(1)}_{h_i}(b, w,
N)s^{(1)}_{b,\,\widetilde{x}_{h_1}(w, v)}(\rho_{h_i}(v))
+D^{(1)}_{h_i}(b, 0, 0)s^{(2)}_{b,\,\widehat{x}}(x^\ast_i(w, v),
\rho_{h_i}(v))) \nn\\&&+\sum_{1\le i\le p} (D^{(2)}_{h_i}(b, 0,
0)s^{(2, +)}_{b,\,\widehat{x}}(\rho_{h_i}(v)) +D^{(2)}_{h_i}(b, 0,
0)s^{(2, -)}_{b,\,\widehat{x}}(\rho_{h_i}(v))) +\epsilon(v, \zeta),
\lb{5.2.12}\eea where \bea &&\|\epsilon(v, \zeta)\|\le
C_K|\rho(v)|^2(|v|+|w|+|N|) +C_K|\rho(v)|(|\rho(v)|+|x^\ast(w,
v)|)\cdot\|\zeta\|_{v,\,p,\,1} \nn\\&&+C_K(|x^\ast(w,
v)||\rho(v)|(|x^\ast(w, v)|+|v|+|w|+|N|). \lb{5.2.13}\eea Note that
$\rho_{h_i}(v)=v_\alpha  v_{h_i}$ and $\widetilde{x}_{h_i}(w,
v)=v_\alpha x_{h_i}$ for $1\le i\le q$,
$\rho_{h_i}(v)=v_\beta  v_{h_i}$ and $\widetilde{x}_{h_i}(w,
v)=w+v_\beta x_{h_i}$ for $q+1\le i\le p$. Hence we have
$x^\ast_i(w, v)=v_\alpha(x_{h_i}-x_{h_1})\neq 0$ for $2\le
i\le q$ and $x^\ast_i(w, v)=w+v_\beta x_{h_i}-v_\alpha x_{h_1}$
for $q+1\le i\le p$. Clearly we may assume
$|x^\ast_i(w, v)|\ge\delta_K|v_\beta|$
for $q+1\le i\le p-1$.
Hence we have \bea|s^{(2,
-)}_{b,\,\widehat{x}}(\rho_{h_i}(v))| =o(|s^{(2,
-)}_{b,\,\widehat{x}}(x^\ast_i(w, v), \rho_{h_i}(v))|), \qquad 2\le
i\le p-1\lb{5.2.14}\eea
By changing the order of $\{h_1, \ldots, h_{q}\}$,
we may assume $\rho_{h_1}(v)=\min_{1\le i\le q}\rho_{h_i}(v)$.
Thus we have
\bea|s^{(2, -)}_{b,\,\widehat{x}}(\rho_{h_1}(v))| =\sum_{2\le i\le q}o(|s^{(2,
-)}_{b,\,\widehat{x}}(x^\ast_i(w, v), \rho_{h_i}(v))|).\lb{5.2.15}\eea
We have the following two cases:

{\bf Case 1.} {\it We have $D^{(1)}_{h_p}b\neq0$.}

In this case we have $|\rho_{h_p}(v)|\le C_K\sum_{1\le i\le p-1}|\rho_{h_i}(v)|$
by $\alpha_T(v)=0$. Hence by (\ref{5.2.14}) and (\ref{5.2.15}), we have
\bea|s^{(2, -)}_{b,\,\widehat{x}}(\rho_{h_p}(v))| =\sum_{2\le i\le p-1}o(|s^{(2,
-)}_{b,\,\widehat{x}}(x^\ast_i(w, v), \rho_{h_i}(v))|).\lb{5.2.16}\eea
Since $\widehat{x}\in\Sg_P^\ast$, we have $s^{(2,
-)}_{b,\,\widehat{x}}\neq 0$. Hence in order to satisfies
(\ref{5.21}), we must have \bea \widetilde{\alpha}_T(w,
v)&\equiv&\sum_{2\le i\le p} D^{(1)}_{h_i}(b, 0, 0)s^{(2,
-)}_{b,\,\widehat{x}}(x_i^\ast(w, v), \rho_{h_i}(v))=0 \lb{5.2.17}\eea provided
$|v|$, $|w|$ and $|N|$ are sufficiently small.

{\bf Case 2.} {\it We have $D^{(1)}_{h_p}b=0$.}

In this case,  in order to satisfies
(\ref{5.21}), we must have \bea \widetilde{\alpha}_T(w,
v)&\equiv&\sum_{2\le i\le p} D^{(1)}_{h_i}(b, 0, 0)s^{(2,
-)}_{b,\,\widehat{x}}(x_i^\ast(w, v), \rho_{h_i}(v))
+D^{(2)}_{h_p}(b, 0, 0)s^{(2,
-)}_{b,\,\widehat{x}}(\rho_{h_p}(v))=0 \lb{5.2.18}\eea provided
$|v|$, $|w|$ and $|N|$ are sufficiently small.

The proof of the converse is similar to the previous theorems since
the map $\widetilde{\alpha}_T$ and $\alpha_T$ are transversal to the
zero sections in both cases by (ii-a) and (ii-b) in Definition 1.1. \hfill\hb

\subsection{ There are at least three  bubble trees}

In this subsection, we consider the case $|I_1|\ge 3$. Note that in
this case the rank of $\alpha_T(v)$ is $4n$. Hence we have the following:

{ \bf Theorem 5.3.1.} {\it Suppose $T$ is a bubble type given by
(iv) of Theorem 2.9 with $\Sg_P$ being smooth and $|I_1|\ge 3$,
i.e., there are at least three bubble trees. Then an element
$b\equiv[\mathcal{C},\, u]\in \mathfrak{M}_T(X, A, J)\cap\overline{\mathfrak{M}}^0_{2, k}(X, A, J)$ if and only if
$\alpha_T(v)=0$ for some $v=(b,
v)\in\widetilde{\mathcal{FT}}_{\delta_K}^\emptyset|_{\widetilde{K}^{(0)}}$.
In particular, $\mathfrak{M}_T(X, A, J)
\cap\overline{\mathfrak{M}}^0_{2, k}(X, A, J)$ is  a smooth
orbifold of dimension at most $\dim
\overline{\mathfrak{M}}^{vir}_{2,k}(X, A)-2$.  }

{\bf Proof.} By Lemma 5.5, an element $b\equiv[\mathcal{C},\, u]\in
\mathfrak{M}_T(X, A, J)\cap\overline{\mathfrak{M}}^0_{2, k}(X, A, J)$ must satisfy $\alpha_T(v)=0$ for some $v=(b,
v)\in\widetilde{\mathcal{FT}}_{\delta_K}^\emptyset|_{\widetilde{K}^{(0)}}$.
Note that  in this case we have
$\mathrm{rank}(\alpha_T)=4n=\dim(\mathrm{coker}D_b)$ by (ii-a) in Definition 1.1.
Thus the theorem follows by a similar argument as in Theorem 5.2.1.
 \hfill\hb

\subsection{ The principle component $\Sg_P$ is not smooth and $u_\ast[\Sg_P]=0$}

In this subsection, we consider the case
that the principle component $\Sg_P$ is not smooth
and $u_\ast[\Sg_P]=0$, i.e., $u|_{\Sg_P}=const$.

First we use \S3.2 to smooth out all the nodes in $\Sg_P$, then
we use the methods in \S5.1-5.3 to study the conditions under which
an element $b\equiv[\mathcal{C},\, u]\in
\mathfrak{M}_T(X, A, J)\cap\overline{\mathfrak{M}}^0_{2, k}(X, A, J)$.
Denote by $(\Sg_P, v_P)$ the smooth Riemann surface obtained from $\Sg_P$
by smoothing nodes. By the gluing construction, we have an isomorphism
$R_{v_P}: \mathcal{H}^{0,\, 1}_{\Sg_P} \rightarrow\mathcal{H}^{0,\, 1}_{(\Sg_P,\, v_P)}$
which depends continuously on the parameter $v_P$. In fact,
we can choose  basis of the Hodge bundle
of holomorphic 1-forms over the Deligne-Mumford space
$\overline{\mathcal{M}}_{2, k_P}$ depending continuously
on the parameter $v_P$, where $k_P=|I_1\cup M_P|$.
Then this case follows by a similar  argument as in the previous sections by replacing the term
$R_v$ by $R_{(v_0, v_1)}\circ R_{v_P}$,
where $v_1\equiv \{v^{(l)}\}_{l\in I_1}$.

The following is the main theorem in this section.

{ \bf Theorem 5.4.1.} {\it Suppose $T$ is a bubble type given by (iv) of
Theorem 2.9 with $u_\ast[\Sg_P]=0$. Then one of the following must hold:

(i) An element $b\equiv[\mathcal{C},\, u]\in
\mathfrak{M}_T(X, A, J)\cap\overline{\mathfrak{M}}^0_{2, k}(X, A, J)$
if and only if $\{D^{(1)}_{h}b,\,D_{h} ^{(2)}b, \,D_{h} ^{(3)}b\}_{h\in\chi(T)}$
satisfy a set of linear equations of rank $4n$ whose
coefficients depending on $v=(b, v)\in\widetilde{\mathcal{FT}}_{\delta_K}^\emptyset|_{\widetilde{K}^{(0)}}$
and the position of nodes on $\mathcal{C}$. In particular,
$\mathfrak{M}_T(X, A, J)\cap\overline{\mathfrak{M}}^0_{2, k}(X, A, J)$
is a smooth orbifold of dimension at most $\dim
\overline{\mathfrak{M}}^{vir}_{2,k}(X, A)-2$
.

(ii) An element $b\equiv[\mathcal{C},\, u]\in
\mathfrak{M}_T(X, A, J)\cap\overline{\mathfrak{M}}^0_{2, k}(X, A, J)$
must satisfy: there exists $h\in\chi(T)$ such that $u_h$ factor through a branched covering $\widetilde{u}: S^2\rightarrow
X$, i.e., there exists a holomorphic branched covering $\phi:
S^2\rightarrow S^2$ such that $u_h=\widetilde{u}\circ\phi$ and
$\deg(\phi)\ge 2$.

(iii) An element $b\equiv[\mathcal{C},\, u]\in
\mathfrak{M}_T(X, A, J)\cap\overline{\mathfrak{M}}^0_{2, k}(X, A, J)$
must satisfy: there exist $h_1, h_2\in\chi(T)$ such that
$u_{h_1}(\Sg_{h_1})=u_{h_2}(\Sg_{h_2})$.
}

{\bf Proof.} We only give the proof of two cases:
$\Sg_P$ is a torus with only one node and $\Sg_P$ is
obtained from two tori attached at one node. The other cases follow
by a similar argument, so we omit their proofs here.

{\bf Case 1.} {\it The principle component $\Sg_P$ is a torus with only one node.}

We will construct an isomorphism
$R_{v_P}: \mathcal{H}^{0,\, 1}_{\Sg_P}\rightarrow\mathcal{H}^{0,\, 1}_{(\Sg_P,\, v_P)}$
depending continuously on the parameter $v_P$.
Suppose $(T, x_1, x_2)$ is the normalization of $\Sg_P$.
Then we have
$$ H^1(\Sg_P, \mathcal{O)}\cong H^1(T, \mathcal{O}(-x_1-x_2))\cong H^0(T, \mathcal{O}(x_1+x_2)\otimes K_T)^\ast\cong\C^2
$$
which consists of meromorphic one-forms $\omega$ on $T$ that are holomorphic
on $T\setminus\{x_1, x_2\}$ and have at most simple poles at $x_1$ and $x_2$ together with
$\mathrm{Res}_{x_1}\omega+\mathrm{Res}_{x_2}\omega=0$ (cf. \S22.3 of \cite{MirSym}).
Note that by the Residue Theorem (cf. P.222 of \cite{GH}),
the condition $\mathrm{Res}_{x_1}\omega+\mathrm{Res}_{x_2}\omega=0$ is
automatically satisfied.

Now we study nonzero sections of $H^1(\Sg_P, \mathcal{O)}$.
Note that if $\mathrm{Res}_{x_1}\omega_1=0=\mathrm{Res}_{x_2}\omega_1$,
then the one-form $\omega_1$ is holomorphic on $T$, i.e.,
$\omega_1\in H^0(T, K_T)^\ast\cong H^0(T, \mathcal{O})^\ast\cong\C$.
Thus $\omega_1$ is nowhere vanishing on $T$.
If $\mathrm{Res}_{x_1}\omega_2=-\mathrm{Res}_{x_2}\omega_2\neq 0$,
then the one-form $\omega_2$ has exactly two simple poles on $T$.
We claim that $\omega_2$ has at most two zeros on $T$.
In fact, if $\omega_2$ has $l$ zeros, then we have
$\omega_2\in H^0(T, \mathcal{O}(2-l)\otimes K_T)^\ast=0$
provided $l>2$. Now we fix an $h_\ast\in\chi(T)$ and
Let $\{\psi_{1}, \psi_{2}\}$ be an orthogonal basis of
$\mathcal{H}_{\Sg_P}^{0,\,1}$
such that $\psi_{1}(\widetilde{x}_{h_\ast}(v))\neq0$ and
$\psi_{2}(\widetilde{x}_{h_\ast}(v))=0$.
We extend $\{\psi_{1}, \psi_{2}\}$ to be a basis
$\{\psi_{v_P, 1}, \psi_{v_P, 2}\}$ of $\mathcal{H}_{(\Sg_P, v_P)}^{0,\,1}$
such that $\psi_{v_P, 1}(\widetilde{x}_{h_\ast}(v))\neq0$ and
$\psi_{v_P, 2}(\widetilde{x}_{h_\ast}(v))=0$.
Note that $\{\psi^{(m)}_{v_P, j}\}$ is  continuous with respect
to the parameter $v_P$ outside small neighborhoods $U_i$ of $x_i$
for any $m\ge 0$ and the attaching nodes of bubble trees lies
outside $U_i$. Hence Theorem 5.4.1 holds in this case by a similar
argument as in \S5.1-5.3 by considering expansions of
the approximately $J$-holomorphic curves.

{\bf Case 2.} {\it The principle component $\Sg_P$ is
two tori attached at one node.}

We will construct an isomorphism
$R_{v_P}: \mathcal{H}^{0,\, 1}_{\Sg_P} \rightarrow\mathcal{H}^{0,\, 1}_{(\Sg_P,\, v_P)}$
depending continuously on the parameter $v_P$.
Suppose $(T_1, x_1)$, $(T_2, x_2)$ is the normalization of $\Sg_P$.
Then we have
\bea &&H^1(\Sg_P, \mathcal{O)}\cong H^1(T_1, \mathcal{O}(-x_1))\oplus H^1(T_2, \mathcal{O}(-x_2))
\nn\\ \cong &&H^0(T_1, \mathcal{O}(x_1)\otimes K_{T_1})^\ast
\oplus H^0(T_2, \mathcal{O}(x_2)\otimes K_{T_2})^\ast\cong\C^2
\nn\eea
which consists of meromorphic one-forms $\omega$ on $T_1\cup T_2$ that are holomorphic
on $T_1\setminus\{x_1\}$ and $T_2\setminus\{x_2\}$ and have at most
simple poles at $x_1$ and $x_2$ together with
$\mathrm{Res}_{x_1}\omega+\mathrm{Res}_{x_2}\omega=0$.
Note that by the Residue Theorem applied to each $T_i$,
we have  $\mathrm{Res}_{x_1}\omega=0=\mathrm{Res}_{x_2}\omega$.
Hence $\omega$ is holomorphic on each $T_i$.
Since $H^0(T_i, K_{T_i})^\ast\cong H^0(T_i, \mathcal{O})^\ast\cong\C$,
we have $\omega|_{T_i}$ is either nowhere vanishing or identically zero.
We choose a basis $\{\psi_1, \psi_2\}$
of $\mathcal{H}^{0,\, 1}_{\Sg_P}$ such that
$\psi_1|_{T_1}\neq 0$,  $\psi_1|_{T_2}= 0$
and $\psi_2|_{T_1}= 0$,  $\psi_2|_{T_2}\neq 0$.
Then we extend $\{\psi_1, \psi_2\}$ to be
a basis $\{\psi_{v_P, 1}, \psi_{v_P, 2}\}$
of $\mathcal{H}_{(\Sg_P, v_P)}^{0,\,1}$.
Since $\psi_2|_{T_1}= 0$, we need to study the behavior
of $\psi_{v_P, 2}|_{T_1}$ carefully in order to derive the
conditions for $b\equiv[\mathcal{C},\, u]\in
\mathfrak{M}_T(X, A, J)\cap\overline{\mathfrak{M}}^0_{2, k}(X, A, J)$.
Write the local coordinate at the node as $\{(z, w)\in\C^2, zw=v_P\}$.
Assume in a small neighborhood of the node $T_1$ and $T_2$ are given by
the coordinate planes $\{z=0\}$ and $\{w=0\}$ respectively.
Suppose $\psi_2|_{T_2}= dz$ and $\psi_2|_{T_1}=0$.
We construct an approximately holomorphic one-form
$\widehat{\psi}_{v_P, 2}$ on $(\Sg_P, v_P)$ as follows:
First note that the one-form  $\psi_2=dz$
has the form $\psi_2=-\frac{v_P}{w^2} dw$ in the $w$-coordinate.
Then we choose a meromorphic one-form $\omega$ on $T_1$ with only one
pole at $x_1$ of order two such that its principle part is
$-\frac{v_P}{w^2} dw$. By the Residue Theorem applied to
$(T_1, \omega)$, we have the expansion of $\omega$ near $x_1$ as
$\omega=v_P(-\frac{1}{w^2}+a_0+a_1w+\cdots) dw$, where $a_i\in\C$
is independent of $v_P$ for each $i\ge 0$.
Now we define
\be \widehat{\psi}_{v_P, 2}=(1-\beta_{|v_P|}(|z|))\psi_2+\beta_{|v_P|}(|z|)\omega,\lb{5.4.1}
\ee
where $\beta$ is given by (\ref{3.5}). Then $\widehat{\psi}_{v_P, 2}$
is a globally defined one form on $(\Sg_P, v_P)$ and
it is holomorphic outside the annulus
$A\equiv\{|v_P|^{\frac{1}{2}}\le |z|\le2|v_P|^{\frac{1}{2}}\}
=\{\frac{1}{2}|v_P|^{\frac{1}{2}}\le |w|\le|v_P|^{\frac{1}{2}}\}$.
Inside the annulus $A$, we have
\be \widehat{\psi}_{v_P, 2}=v_P\beta_{|v_P|}\left(\frac{|v_P|}{|w|}\right)(a_0+a_1w+\cdots) dw.\lb{5.4.2}\ee
Since
$\omega\in H^0(T_1, \mathcal{O}(2)\otimes K_{T_1})^\ast$,
it  has at most two zeros on $T_1$.
In fact if $\omega$ has $l$ zeros, then we have
$\omega\in H^0(T_1, \mathcal{O}(2-l)\otimes K_{T_1})^\ast=0$
provided $l>2$.
Note that we have $\|\overline{\partial}\widehat{\psi}_{v_P, 2}\|_{(\Sg_P, v_P), \,L^p}
\le C|v_P|^{1+\frac{1}{2p}}$ for some constant $C$ independent of
$v_P$ and $p>2$, where near the node $x$ we use the cylinder-like metric.
Since $\overline{\partial}$ is a first-order elliptic operator,
it follows from the standard elliptic estimate that we can find a
one-form $\epsilon$ on $\Sg_P$ such that $\|\epsilon\|_{(\Sg_P, v_P), \,W^{1, p}}\le C_K|v_P|^{1+\frac{1}{2p}}$
and $\psi_{v_P, 2}\equiv\widehat{\psi}_{v_P, 2}+\epsilon$ is holomorphic.
Thus as $v_P$ converges to $0$, we have
$\frac{1}{v_P}\psi_{v_P, 2}$ will $C^\infty$-converge to
$\frac{1}{v_P}\widehat{\psi}_{v_P, 2}$ outside small neighborhood $U$
of the node $x_1$ in $T_1$. In fact, this holds since
both $\frac{1}{v_P}\psi_{v_P, 2}$ and $\frac{1}{v_P}\widehat{\psi}_{v_P, 2}$
are holomorphic outside $U$ and two holomorphic functions
are $C^0$-close implies they are $C^\infty$-close.
We construct $\psi_{v_P, 1}$ similarly.

Now we fix an $h_\ast\in\chi(T)$ and let
$\{\psi_{v_P, 1}, \psi_{v_P, 2}\}$ be a basis of $\mathcal{H}_{\Sg_P, v_P}^{0,\,1}$
constructed as above such that $\psi_{v_P, 1}(\widetilde{x}_{h_\ast}(v))\neq0$ and
$\psi_{v_P, 2}(\widetilde{x}_{h_\ast}(v))=0$.
Note that  in the proof of the theorems
in \S5.1-5.3, we only use $\frac{\psi_{v_P, 2}^{(m+l)}(y)}{\psi_{v_P, 2}^{(m)}(y)}\le C_K$
where $m$ is the order of $y$ as a zero of $\psi_{v_P, 2}$ and $l\in\N$.
By the above construction,
$\frac{\psi_{v_P, 2}^{(m+l)}(y)}{\psi_{v_P, 2}^{(m)}(y)}$
will $C^0$-converge to
$\frac{\widehat{\psi}_{v_P, 2}^{(m+l)}(y)}{\widehat{\psi}_{v_P, 2}^{(m)}(y)}$
outside small neighborhood $U$ of the node $x_1$ in $T_1$. While
$\frac{\widehat{\psi}_{v_P, 2}^{(m+l)}(y)}{\widehat{\psi}_{v_P, 2}^{(m)}(y)}$
is independent of $v_P$.
Hence Theorem 5.4.1 holds in this case by a similar
argument as in \S5.1-5.3 by considering expansions of
the approximately $J$-holomorphic curves.

At last we sketch the proof of the general case.
Firstly we glue certain rational components in $\Sg_P$
to obtain $\Sg_P^\prime$ which belongs to one of the above cases.
Then we use the above to show the theorem holds.
\hfill\hb

\subsection{ The principle component $\Sg_P$ is not smooth and $u_\ast[\Sg_P]\neq0$}

In this subsection, we consider the case
that the principle component $\Sg_P$ is not smooth
and $u_\ast[\Sg_P]\neq0$.
The simplest example in this case is given by the second figure
in Figure 2.1.

We denote by $\Sg$ the union of components
of $\Sg_P$ which are mapped to constants such that each connected component
of $\Sg$ has genus greater than zero. Since $\deg(u|_{\Sg_P})\neq0$,
$\Sg$ consists of exactly two connected components $\Sg_1$ and $\Sg_2$,
each one is mapped to a constant, e.g. the two tori in Figure 2.1.
Note that we have
$$ H^1_{\overline{\partial}}(\mathcal{C}, u^\ast TX)\cong (T_{ev(\Sg_1)}X\otimes \mathcal{H}^{0, 1}_{\Sg_1})\oplus( T_{ev(\Sg_2)}X\otimes\mathcal{H}^{0, 1}_{\Sg_2})\cong\C^{2n}.$$
Let $\widetilde{\Sg}_1\equiv\overline{\Sg_P\setminus\Sg_2}\supset\Sg_1$ and
$\widetilde{\Sg}_2\equiv\overline{\Sg_P\setminus\Sg_1}\supset\Sg_2$,
e.g. the left torus with the central sphere or the right torus with the central sphere
respectively in Figure 2.1.
Denote by
$$\Sg^{(1)}\equiv\widetilde{\Sg}_1\cup\{T_B^{(h)}\}_{x_h\in\widetilde{\Sg}_1},
\qquad
\Sg^{(2)}\equiv\widetilde{\Sg}_1\cup\{T_B^{(h)}\}_{x_h\in\widetilde{\Sg}_2},
$$
where $x_h$ is the attaching node of the bubble tree $T_B^{(h)}$,
e.g. the genus-one curve obtained form $\mathcal{C}$
by removing the right torus together with all the bubbles descendent from it
or removing the left torus together with all the bubbles descendent from it
respectively in Figure 2.1..
Then both $\Sg^{(1)}$ and $\Sg^{(2)}$ are nodal Riemann
surfaces of genus one. Denote by $\Sg^{(1)}_P$ and $\Sg^{(2)}_P$
the principle components  of $\Sg^{(1)}$ and $\Sg^{(2)}$ respectively,
e.g. the two tori in Figure 2.1.

We only consider the simplest case that both  $\Sg^{(1)}_P$
and $\Sg^{(2)}_P$ are smooth tori, i.e.,
as illustrated in  the second figure in Figure 2.1, the general case follows similarly
as explained in \S5.4.
For $v=(b, v)\in\widetilde{\mathcal{FT}}_{\delta_T}^\emptyset$
sufficiently small and
$(V_1\otimes\psi_1, V_2\otimes\psi_2)\in
(T_{ev(\Sg_1)}X\otimes \mathcal{H}^{0, 1}_{\Sg^{(1)}_P})\oplus( T_{ev(\Sg_2)}X\otimes\mathcal{H}^{0, 1}_{\Sg^{(2)}_P})$,
we define $R_v(V_1\otimes\psi_1, V_2\otimes\psi_2)\in\Gamma^{0, 1}(u_v)$ as the following:
if $z\in\Sg_v$ is such that $q_v(z)\in\Sg_{b,\, h}$
for some $h\in\chi(\Sg^{(1)})$ as defined in (\ref{3.33}) and $|q_S^{-1}(q_v(z))|\le 2\delta_T(b)$,
we define $\overline{u}_v(z)\in T_{u(\Sg_1)}X$ by
$\exp_{u(\Sg_1)}\overline{u}_v(z)=u_v(z)$.
Given $z\in\Sg_v$,
let $h_z$ be such that $q_v(z)\in \Sg_{b,\, h_z}$. If $w\in T_z\Sg_v$, put
\be R_vV_1\psi_1|_zw=
\left\{\matrix{0 ,&&\quad {\rm if}\quad  \chi_{\Sg^{(1)}}h_z=2;\cr
\beta(\delta_T(b)|q_vz|)(\psi_1|_zw)\Pi_{\overline{u}_v(z)}V_1,&&\quad {\rm if}\quad \chi_{\Sg^{(1)}}h_z=1;\cr
(\psi_1|_zw)V_1,&&\quad {\rm if}\quad \chi_{\Sg^{(1)}}h_z=0,\cr}\right.\lb{5.5.1}\ee
where $\chi_{\Sg^{(1)}}$ is defined similar to (\ref{5.7})
with respect to the bubble type $\Sg^{(1)}$
and $\Pi_{\overline{u}_v(z)}$ is the parallel transport along the
geodesic $t\mapsto\exp_{u(\Sg_1)}t\overline{u}_v(z)$ with respect
to the Levi-Civita connection of the metric $g_{X,\,\{u(\Sg_1), u(\Sg_2)\}}$,
which are flat near $(u(\Sg_1)$ and $u(\Sg_2)$.
We define $R_vV_2\psi_2$ similarly.
Note that since  $u_\ast[\Sg_P]\neq0$,
there must be a rational component $h$ between $\Sg_1$ and $\Sg_2$
such that $u_\ast[\Sg_h]\neq0$,
thus the map $R_v(V_1\otimes\psi_1, V_2\otimes\psi_2)$  is well-defined.
Let $\Gamma_\pm(v)$ be given by the formula in \S4.

Comparing with the previous sections, we have the following:

{ \bf Lemma 5.5.1.} {\it Suppose $T$ is a bubble type given by
(iv) of Theorem 2.9 and $u_\ast[\Sg_P]\neq0$. Then for every precompact open subset $K$ of $\mathfrak{M}_T(X, A, J)$,
there exist $\delta_K, C_K\in\R^+$ and an open neighborhood
$U_K$ of $K$ in $\mathfrak{X}_{2, k}(X, A, J)$ with the following properties:

(i) For every $[\widetilde{b}]\in\mathfrak{X}^0_{2, k}(X, A, J)\cap U_K$,
there exist $v=(b, v)\in\widetilde{\mathcal{FT}}_{\delta_K}^\emptyset|_{\widetilde{K}^{(0)}}$,
and $\zeta\in\Gamma_+(v)$ such that $\|\zeta\|_{v, p, 1}<\delta_K$
and $[\exp_{u_v}\zeta]=[\widetilde{b}]]$.

(ii) For every $v=(b, v)\in\widetilde{\mathcal{FT}}_{\delta_K}^\emptyset|_{\widetilde{K}^{(0)}}$,
we have
\bea C_K^{-1}\|\xi\|_{v, p, 1}\le\|D_v\xi\|_{v, p}\le C_K\|\xi\|_{v, p, 1}, \quad\forall \xi\in\Gamma_+(v),
\lb{5.5.2}\eea

(iii) For every $v=(b, v)\in\widetilde{\mathcal{FT}}_{\delta_K}^\emptyset|_{\widetilde{K}^{(0)}}$
and $\eta\equiv(V_1\otimes\psi_1, V_2\otimes\psi_2)\in H^1_{\overline{\partial}}(\mathcal{C}, u^\ast TX)$, we have
\bea ||D_v^\ast R_v\eta||_{v,\,C^0}
\le C(b)\left(\sum_{h\in\chi(\Sg^{(1)})}|\rho_h^{(1)}(v)|
+\sum_{h\in\chi(\Sg^{(2)})}|\rho_h^{(2)}(v)|\right)\|\eta||_2.
\lb{5.5.3}\eea

(iv) For every $v=(b, v)\in\widetilde{\mathcal{FT}}_{\delta_K}^\emptyset|_{\widetilde{K}^{(0)}}$
and $\eta\equiv(V_1\otimes\psi_1, V_2\otimes\psi_2)\in H^1_{\overline{\partial}}(\mathcal{C}, u^\ast TX)$, we have
\bea&&\left |\langle\langle \overline{\partial}_Ju_v,\, R_v\eta\rangle\rangle_{v, 2}+
\sum_{h\in\chi(\Sg^{(1)})}\langle D^{(1)}_hb, \,V_1\rangle
\overline{\psi_1(\rho_h^{(1)}(v))}
+\sum_{h\in\chi(\Sg^{(2)})}\langle D^{(2)}_hb, \,V_2\rangle
\overline{\psi_2(\rho_h^{(2)}(v))}
\right|\nn\\
\le &&C_K|v|\cdot|\rho(v)|\cdot\|\eta\|,
\lb{5.5.4}\eea
where $D^{(1)}_hb$ is given by (\ref{5.3})
and $\rho_h^{(i)}(v)$ is given by (\ref{3.32})
with respect to the bubble type $\Sg^{(i)}$.}

{\bf Proof.} (i) and (ii) hold by a similar argument in \cite{Z2}.

We prove (iii). Note that by \S3.2, the gluing construction
(\ref{3.30}) in the principle component $\Sg_P$ coincide
with the gluing construction in $\Sg^{(1)}$ or $\Sg^{(2)}$
as  bubble types of genus-one by replacing
the term $q_S(p_{h, (x_h, v_h)}(z))$ in (\ref{3.17})
by $\phi_{x_h, 1}^{-1}(\widetilde{p}_{h, (x_h, v_h))}(z))$,
i.e., we have $\phi_{x_h, 1}^{-1}(\xi)=q_S(\overline{\xi})$.
Thus the proof of Lemma 5.2 remains valid, which yields (iii).

We prove (iv).  By
the construction of $q_v$ and $R_v$, we have
$\langle\overline{\partial}_Ju_v,\, R_v\eta\rangle=0$ outside the
annuli $A^-_{v, h}(|v_h|)$ for $h\in\chi(\Sg^{(1)})\bigsqcup\chi(\Sg^{(2)})$.
Hence we have \be\langle\langle
\overline{\partial}_Ju_v,\, R_v\eta\rangle\rangle_{v,\,2}=
\sum_{h\in\chi(\Sg^{(1)})\bigsqcup\chi(\Sg^{(2)})}\int_{A^-_{v, h}(|v_h|)}
\langle\overline{\partial}_Ju_v,\, R_v\eta\rangle.\lb{5.5.5} \ee
Now we consider  $\Sg^{(1)}$ and $\Sg^{(2)}$  separately.
Then the proof of  Lemma 4.3 in \cite{Z1} can be used to obtain an
expansion
\bea\langle\langle \overline{\partial}_Ju_v,\, R_v\eta\rangle\rangle_{v,\,2}
= &&-\sum_{m\ge1,\,h\in\chi(\Sg^{(1)})}\langle D^{(m)}_hb, \,V_1\rangle
\overline{\left(\{D_{b,\,\widetilde{x}_h^{(1)}(v)}^{(m)}\psi_1\}
((d\phi_{b,\, {\mathcal{T}^{(1)}(h)}}|_{\widetilde{x}_h^{(1)}(v)})^{-1}\rho_h^{(1)}(v))\right)}
\nn\\&&-\sum_{m\ge1,\,h\in\chi(\Sg^{(2)})}\langle D^{(m)}_hb, \,V_2\rangle
\overline{\left(\{D_{b,\,\widetilde{x}_h^{(2)}(v)}^{(m)}\psi_2\}
((d\phi_{b,\, {\mathcal{T}^{(2)}(h)}}|_{\widetilde{x}_h^{(2)}(v)})^{-1}\rho_h^{(2)}(v))\right)}
\nn\eea
where $\mathcal{T}^{(i)}(h)$, $\widetilde{x}_h^{(i)}(v)$,
$\rho_h^{(i)}(v)$  are the maps as in Lemma 5.1 with respect to
the bubble type $\Sg^{(i)}$ for $i=1$ or $2$.
Hence (iv) holds by the same argument as in Lemma 5.4.
\hfill\hb

For any $v=(b, v)\in\widetilde{\mathcal{FT}}^{\emptyset}$,
let
\be \alpha_{T,\, 1}(v)=\sum_{h\in\chi(\Sg^{(1)})}
(D^{(1)}_hb)\overline{\psi_1(\rho_h^{(1)}(v))},
\quad \alpha_{T,\, 2}(v)=\sum_{h\in\chi(\Sg^{(2)})}
(D^{(1)}_hb)\overline{\psi_2(\rho_h^{(2)}(v))}
\lb{5.5.6}\ee

The following is the main result in this subsection.

{ \bf Theorem 5.5.2.} {\it Suppose $T$ is a bubble type given by
(iv) of Theorem 2.9 and $u_\ast[\Sg_P]\neq0$.
Then one of the following must hold:

(i) An element $b\equiv[\mathcal{C},\, u]\in
\mathfrak{M}_T(X, A, J)\cap\overline{\mathfrak{M}}^0_{2, k}(X, A, J)$
if and only if $\alpha_{T, 1}(v)=0$ and $\alpha_{T, 2}(v)=0$ for some $v=(b,
v)\in\widetilde{\mathcal{FT}}_{\delta_K}^\emptyset|_{\widetilde{K}^{(0)}}$. In particular,
$\mathfrak{M}_T(X, A, J)\cap\overline{\mathfrak{M}}^0_{2, k}(X, A, J)$
is a smooth orbifold of dimension at most $\dim
\overline{\mathfrak{M}}^{vir}_{2,k}(X, A)-2$.

(ii) An element $b\equiv[\mathcal{C},\, u]\in
\mathfrak{M}_T(X, A, J)\cap\overline{\mathfrak{M}}^0_{2, k}(X, A, J)$
must satisfy: there exists $\Sg_h\subset\Sg_P$ such that $u_h$ factor through a branched covering $\widetilde{u}: S^2\rightarrow
X$, i.e., there exists a holomorphic branched covering $\phi:
S^2\rightarrow S^2$ such that $u_h=\widetilde{u}\circ\phi$ and
$\deg(\phi)\ge 2$.
 }

{\bf Proof.} Let $\delta_K$ be given by Lemma 5.5.1.
For each $v=(b, v)\in\widetilde{\mathcal{FT}}_{\delta_K}^\emptyset|_{\widetilde{K}^{(0)}}$,,
we define the homomorphism
\be \pi_{v, -}^{0, 1}:
\Gamma^{0, 1}(v)\rightarrow\Gamma^{0, 1}_-(b_P),
\qquad \pi_{v, -}^{0, 1}\xi=-\sum_{1\le i\le n,\,1\le j\le 2}\langle \xi, R_ve^j_i\psi_j\rangle e_i^j\psi_j\in\Gamma^{0, 1}_-(b_P),
\lb{5.5.7}\ee
where  $\psi_1$ and $\psi_2$ are orthonormal basis
for $\mathcal{H}^{0, 1}_{\Sg_1}$ and $\mathcal{H}^{0, 1}_{\Sg_2}$
and
$\{e_i^j\}_{1\le i\le n,\, 1\le j\le2}$ are orthonormal basis for
$ T_{u(\Sg_j)}X$. Denote the kernel of  $\pi_{v, -}^{0, 1}$ by $\Gamma_+^{0, 1}(v)$.
Then by Lemmas 5.5.1 and the same argument as in Lemma 5.5, we have
\bea \pi_{v, -}^{0, 1}(v, \zeta)\equiv\pi_{v, -}^{0, 1}(\overline{\partial}_Ju_v+D_v\zeta+N_v\zeta)
=\alpha_{T, 1}(v)+\alpha_{T, 2}(v)+\epsilon(v, \zeta),
\nn\eea
and
\bea \|\epsilon(v, \zeta)\|\le C_K(|v|+\|\zeta\|_{v, p, 1})|\rho(v)|
\le C_K|v|^\frac{1}{p}|\rho(v)|,
\nn\eea
where we use notations in Lemma 5.5.
Hence in order to satisfies (\ref{5.21}),
we must have $\alpha_{T, 1}(v)=0$ and $\alpha_{T, 2}(v)=0$  provided
$|v|$ is sufficiently small.

The converse is similar to the previous by using (ii-a), (ii-b) and (ii-d)
in Definition 1.1. \hfill\hb

\setcounter{equation}{0}
\section{ Study for $\mathfrak{M}_T(X, A, J)$ in (iii) of Theorem 2.9}

In this case we have $\mathrm{coker}D_b=H^1_{\overline{\partial}}(\mathcal{C}, u^\ast TX)=\C^n$
by Theorem 2.10. The simplest examples in this case are illustrated in
Figure 2.3.

We denote by $\Sg_1$ the union of components
of $\Sg_P$ which are mapped to constants such that
$\Sg_1$ is connected and has genus one, e.g. the  tori
on the left hand side of the two figures in Figure 2.3. Then  we have
$$ H^1_{\overline{\partial}}(\mathcal{C}, u^\ast TX)\cong T_{ev(\Sg_1)}X\otimes \mathcal{H}^{0, 1}_{\Sg_1}.$$
Denote by
$\Sg^{(1)}\equiv\Sg_1\cup\{T_b^{(h)}\}_{x_h\in\Sg_1}$.
Then  $\Sg^{(1)}$ is a nodal Riemann surface of genus one.
Let $\widehat{\chi}(T)=\{h: |\Sg_h\subset\overline{\Sg_P\setminus\Sg_1}, \,nod(h)\cap(\Sg_1\cap\overline{\Sg_P\setminus\Sg_1})\neq\emptyset\}$,
where $nod(h)$ denotes the nodes on $h$.
For $h\in \widehat{\chi}(T)$,  $x\in \Lambda(h)\equiv nod(h)\cap(\Sg_1\cap\overline{\Sg_P\setminus\Sg_1})$
and $m\in\mathds{N}$, define
\be D^{(m)}_{x}b=\frac{2}{(m-1)!}\frac{D^{m-1}}{ds^{m-1}}
\frac{d}{ds}(u_h\circ \phi_{x, 1}^{-1}(\overline{z}))\left|\frac{}{}\right._{z=(s, \,t)=0},\lb{6.1}\ee
where $\phi_{x, 1}$ is given by (\ref{3.30}) and the covariant derivatives are taken with respect to
the standard metric $s+it\in\C$ and a metric $g_{X,\,u(\Sg_1)}$ on $X$
which is flat near $u(\Sg_1)$.
Here $\Sg_{x, 1}=h$ and  $\Sg_{x, 0}$ is the component in $\Sg_1$
containing $x$, we identify $ T_{x, 0}\Sg_{x, 0}$
and $T_{x, 1}\Sg_{x, 1}$ with $\C$ via $\phi_{x, 0}$ and
$\phi_{x, 1}$ respectively. For example, $h$ is the right
sphere in Figure 2.3 with $|\Lambda(h)|=2$ or $1$ respectively
in Figure 2.3.

We study the case that the principle component $\Sg^{(1)}_P$ of
$\Sg^{(1)}$ is a smooth torus, the general cases follow similarly as
explained in \S5.4. For $v=(b,
v)\in\widetilde{\mathcal{FT}}_{\delta_T}^\emptyset$ sufficiently
small and $V\psi\in T_{u(\Sg_P^{(1)})}X\otimes\mathcal{H}^{0,
1}_{\Sg^{(1)}_P}$, define $R_vV\psi\in\Gamma^{0, 1}(u_v)$ as
follows. If $z\in\Sg_v$ is such that $q_v(z)\in\Sg_{b,\, h}$ for
some $h\in\chi(\Sg^{(1)})$ and $|q_S^{-1}(q_v(z))|\le 2\delta_T(b)$
or $h\in\widehat{\chi}(T)$, $x\in \Lambda(h)$ and $|\phi_{x,
1}(q_v(z))|\le 2\delta_T(b)$, we define $\overline{u}_v(z)\in
T_{u(\Sg_1)}X$ by $\exp_{u(\Sg_1)}\overline{u}_v(z)=u_v(z)$.
Given $z\in\Sg_v$, let $h_z$ be such that $q_v(z)\in \Sg_{b,\,
h_z}$. If $w\in T_z\Sg_v$, put \be R_vV\psi|_zw= \left\{\matrix{0
,&&\quad {\rm if}\quad  \chi_{\Sg^{(1)}}h_z=2;\cr
\beta(\delta_T(b)|q_vz|)(\psi|_zw)\Pi_{\overline{u}_v(z)}V,&&\quad
{\rm if}\quad h_z\in\chi_{\Sg^{(1)}}\cup\widehat{\chi}(T);\cr
(\psi|_zw)V,&&\quad {\rm if}\quad
\chi_{\Sg^{(1)}}h_z=0,\cr}\right.\lb{6.2}\ee where
$\chi_{\Sg^{(1)}}$ is defined similar to (\ref{5.7}) with respect to
the bubble type $\Sg^{(1)}$ and $\Pi_{\overline{u}_v(z)}$ is the
parallel transport along the geodesic
$t\mapsto\exp_{u(\Sg_1)}t\overline{u}_v(z)$ with respect to the
Levi-Civita connection of the metric $g_{X,\,u(\Sg_1)}$ given by
Lemma 3.4. Let $\Gamma_\pm(v)$ be given by the formula in \S4.

Comparing with \S5, we have the following lemma.

{ \bf Lemma 6.1.} {\it Suppose $T$ is a bubble type given by (iii)
of Theorem 2.9. Then for every precompact open subset $K$ of
$\mathfrak{M}_T(X, A, J)$, there exist $\delta_K, C_K\in\R^+$ and an
open neighborhood $U_K$ of $K$ in $\mathfrak{X}_{2, k}(X, A, J)$
with thr following property:

(i) For every $[\widetilde{b}]\in\mathfrak{X}^0_{2, k}(X, A, J)
\cap U_K$, there exist $v=(b,
v)\in\widetilde{\mathcal{FT}}_{\delta_K}^\emptyset|_{\widetilde{K}^{(0)}}$
and $\zeta\in\Gamma_+(v)$ such that $\|\zeta\|_{v, p, 1}<\delta_K$
and $[\exp_{u_v}\zeta]=[\widetilde{b}]]$.

(ii) For every $v=(b,
v)\in\widetilde{\mathcal{FT}}_{\delta_K}^\emptyset|_{\widetilde{K}^{(0)}}$,
we have \bea C_K^{-1}\|\xi\|_{v, p, 1}\le\|D_v\xi\|_{v, p}\le
C_K\|\xi\|_{v, p, 1}, \quad\forall \xi\in\Gamma_+(v), \lb{6.3}\eea

(iii) For every $v=(b,
v)\in\widetilde{\mathcal{FT}}_{\delta_K}^\emptyset|_{\widetilde{K}^{(0)}}$
and $V\otimes\psi\in T_{ev(\Sg_1)}X\otimes \mathcal{H}^{0,
1}_{\Sg_P^{(1)}}$, we have
\bea &&||D_v^\ast R_v(V\otimes\psi)||_{v,\,C^0}\nn\\
\le &&C(b) \left(\sum_{h\in\chi(\Sg^{(1)})}|\rho_h^{(1)}(v)|
+\sum_{h\in\widehat{\chi}(T),\,
x\in\Lambda(h)}|\rho_h^{(x)}(v)|\right)||V\otimes\psi||_2.
\lb{6.4}\eea

(iv) For every $v=(b,
v)\in\widetilde{\mathcal{FT}}_{\delta_K}^\emptyset|_{\widetilde{K}^{(0)}}$
and $V\otimes\psi\in T_{ev(\Sg_1)}X\otimes \mathcal{H}^{0,
1}_{\Sg_P^{(1)}}$, we have \bea&&\left |\langle\langle
\overline{\partial}_Ju_v,\, R_v(V\otimes\psi)\rangle_{v, 2}+
\sum_{h\in\chi(\Sg^{(1)})}\langle D^{(1)}_hb, \,X\rangle
\overline{\psi(\rho_h^{(1)}(v))} +\sum_{h\in\widehat{\chi}(T),\,
x\in\Lambda(h)}\langle D^{(1)}_xb, \,X\rangle
\overline{\psi(\rho_h^{(x)}(v))}
\right|\nn\\
\le &&C_K|v|\cdot|\rho(v)|\cdot\|V\otimes\psi\|, \lb{6.5}\eea where
$\rho_h^{(1)}(v)$ is given by (\ref{3.32}) with respect to the
bubble type $\Sg^{(1)}$ for $h\in\chi(\Sg^{(1)})$ and
$\rho_h^{(x)}(v)=\rho_{\iota_h}(v)v_x$  for
$h\in\widehat{\chi}(T),\, x\in\Lambda(h)$. Here $\iota_h$ denotes
the component in $\Sg_1$ which contains the node $x$ and
$\rho_{\iota_h}(v)$ is given by (\ref{3.32}) with respect to the
bubble type $\Sg^{(1)}$.}

{\bf Proof.} (i) and (ii) hold by a similar argument in \cite{Z2}.

As in Lemma 5.5.1, (iii) follows from  the proof of Lemma 5.2.

We prove (iv).  By the construction of $q_v$ and $R_v$, we have
$\langle\overline{\partial}_Ju_v,\, R_v(V\otimes\psi)\rangle=0$
outside the annuli $A^-_{v, h}(|v_h|)$ for
$h\in\chi(\Sg^{(1)})\cup\hat\chi(T)$£¬ where $A^-_{v,
h}(|v_h|)=\cup_{x\in\Lambda(h)}q_v^{-1}(\{|\phi_{x, 1}(q_v(z))|\le
2|v_x|^{\frac{1}{2}}\})$ for $h\in\hat\chi(T)$. Hence we have
\be\langle\langle \overline{\partial}_Ju_v,\,
R_v(V\otimes\psi)\rangle\rangle_{v,\,2}=
\sum_{h\in\chi(\Sg^{(1)})\cup\hat\chi(T)}\int_{A^-_{v, h}(|v_h|)}
\langle\overline{\partial}_Ju_v,\, R_v(V\otimes\psi)\rangle.\lb{6.6}
\ee Then the proof of  Lemma 4.3 in \cite{Z1} can be uesd to obtain
an expansion \bea\langle\langle \overline{\partial}_Ju_v,\,
R_v(V\otimes\psi)\rangle\rangle_{v,\,2}
&=&-\sum_{m\ge1,\,h\in\chi(\Sg^{(1)})}\langle D^{(m)}_hb, \,X\rangle
\overline{\left(\{D_{b,\,\widetilde{x}_h^{(1)}(v)}^{(m)}\psi\}
((d\phi_{b,\,
{\mathcal{T}^{(1)}(h)}}|_{\widetilde{x}_h^{(1)}(v)})^{-1}\rho_h^{(1)}(v))\right)}
\nn\\&&-\sum_{m\ge1,\,\atop h\in\widehat{\chi}(T),\,
x\in\Lambda(h)}\langle D^{(m)}_xb, \,X\rangle
\overline{\left(\{D_{b,\,\widetilde{x}_h^{(x)}(v)}^{(m)}\psi\}
((d\phi_{b,\,
{\mathcal{T}^{(x)}(h)}}|_{\widetilde{x}_h^{(x)}(v)})^{-1}\rho_h^{(x)}(v))\right)}
\nn\eea where $\mathcal{T}^{(1)}(h)$, $\widetilde{x}_h^{(1)}(v)$,
$\rho_h^{(1)}(v)$  are the maps as in Lemma 5.1 with respect to the
bubble type $\Sg^{(1)}$ and $\mathcal{T}^{(x)}(h)$,
$\widetilde{x}_h^{(x)}(v)$, $\rho_h^{(x)}(v)$ are defined similarly
with respect to the node $x$. Hence (iv) holds by the same argument
as in Lemma 5.4. \hfill\hb

For any $v=(b, v)\in\widetilde{\mathcal{FT}}^{\emptyset}$, let \be
\alpha_T(v)=\sum_{h\in\chi(\Sg^{(1)})}
(D^{(1)}_hb)\overline{\psi(\rho_h^{(1)}(v))}
+\sum_{h\in\widehat{\chi}(T),\, x\in\Lambda(h)}(D^{(1)}_xb)
\overline{\psi(\rho_h^{(x)}(v))}. \lb{6.7}\ee

The following is the main result in this section.

{ \bf Theorem 6.2.} {\it Suppose $T$ is a bubble type given by (iii)
of Theorem 2.9. Then one of the following must hold:

(i) An element $b\equiv[\mathcal{C},\, u]\in
\mathfrak{M}_T(X, A, J)\cap\overline{\mathfrak{M}}^0_{2, k}(X, A, J)$
if and only if $\alpha_{T}(v)=0$ for some $v=(b,
v)\in\widetilde{\mathcal{FT}}_{\delta_K}^\emptyset|_{\widetilde{K}^{(0)}}$. In particular,
$\mathfrak{M}_T(X, A, J)\cap\overline{\mathfrak{M}}^0_{2, k}(X, A, J)$
is a smooth orbifold of dimension at most $\dim
\overline{\mathfrak{M}}^{vir}_{2,k}(X, A)-2$
.

(ii) An element $b\equiv[\mathcal{C},\, u]\in
\mathfrak{M}_T(X, A, J)\cap\overline{\mathfrak{M}}^0_{2, k}(X, A, J)$
must satisfy: there exists $\Sg_h\subset\Sg_P$ such that
$u_\ast[\Sg_P]=u_\ast[\Sg_h]$ and $u_h$ factor through a branched covering $\widetilde{u}: \Sg^\prime\rightarrow
X$, i.e., there exists a holomorphic branched covering $\phi:
\Sg_h\rightarrow \Sg^\prime$ such that $u_h=\widetilde{u}\circ\phi$ and
$\deg(\phi)\ge 2$.
 }

{\bf Proof.} Let $\delta_K$ be given by Lemma 6.1. For each $v=(b,
v)\in\widetilde{\mathcal{FT}}_{\delta_K}^\emptyset|_{\widetilde{K}^{(0)}}$,,
we define the homomorphism \be \pi_{v, -}^{0, 1}: \Gamma^{0,
1}(v)\rightarrow\Gamma^{0, 1}_-(b_P), \qquad \pi_{v, -}^{0,
1}\xi=-\sum_{1\le i\le n}\langle \xi, R_ve_i\psi\rangle
e_i\psi\in\Gamma^{0, 1}_-(b_P), \lb{6.8}\ee where  $\psi$ is an
orthonormal basis for $\mathcal{H}^{0, 1}_{\Sg^{(1)}_P}$ and
$\{e_i\}_{1\le i\le n}$ is an orthonormal basis for $
T_{u(\Sg_1)}X$. Denote the kernel of  $\pi_{v, -}^{0, 1}$ by
$\Gamma_+^{0, 1}(v)$. Then by Lemmas 6.1 and the same argument as in
Lemma 5.5, we have \bea \pi_{v, -}^{0, 1}(v, \zeta)\equiv\pi_{v,
-}^{0, 1}(\overline{\partial}_Ju_v+D_v\zeta+N_v\zeta)
=\alpha_{T}(v)+\epsilon(v, \zeta), \nn\eea and \bea \|\epsilon(v,
\zeta)\|\le C_K(|v|+\|\zeta\|_{v, p, 1})|\rho(v)| \le
C_K|v|^\frac{1}{p}|\rho(v)|, \nn\eea where we use notations in Lemma
5.5. Hence in order to satisfies (\ref{5.21}), we must have
$\alpha_T(v)=0$   provided $|v|$ is sufficiently small.

The converse is similar to the previous and we only sketch it here.
We have the following two cases:

{\bf Case 1.} {\it  If $|\{x\in\Lambda(h)\}_{h\in\widehat{\chi}(T)}|\equiv|\{x_1,\,x_2\}|=2$,
e.g. the first figure in Figure 2.3.}

If $|\widehat{\chi}(T)|\equiv|\{h\}|=1$,  then by
(ii-e) in Definition 1.1, we have (i) or (ii).
In other cases, we have (i) by (ii-a), (ii-b)
in Definition 1.1.

{\bf Case 2.} {\it  If $|\{x\in\Lambda(h)\}_{h\in\widehat{\chi}(T)}|\equiv|\{x\}|=1$,
e.g. the second figure in Figure 2.3.}

In this case, by (ii-e) and (iii-b) in Definition 1.1,
we have (i) or (ii).\hfill\hb

\setcounter{equation}{0}
\section{ Proof of the main theorems}

In this section we give the proofs of the main theorems.

{ \bf Proof of Theorem 1.2.}
The theorem follows by \S2, 4, 5 and 6.\hfill\hb

{ \bf Proof of Theorem 1.3.} If $b\equiv[\mathcal{C}, u]$
belongs to (i) of Theorem 1.2, the boundary component
$\mathfrak{M}_T(X, A, J)\cap\overline{\mathfrak{M}}^0_{2, k}(X, A, J)$
is a smooth orbifold of dimension at most $\dim
\overline{\mathfrak{M}}^{vir}_{2,k}(X, A)-2$
.

If $b\equiv[\mathcal{C}, u]$ belongs to (ii) of Theorem 1.2,
then we have
$$ev(\mathfrak{M}_T(X, A, J)\cap\overline{\mathfrak{M}}^0_{2, k}(X, A, J))\subset ev(\mathfrak{M}_{T^\prime}(X, A^\prime, J))$$
where $T^\prime$ is the bubble type of genus $g$
corresponding to $\mathcal{C^\prime}$.
By (i) and (iii-a)  in Definition 1.1, we have
\be\dim \overline{\mathfrak{M}}_{T^\prime}(X, A^\prime. J)
=2\langle c_1(TX), A^\prime\rangle+2(n-3)(1-g)+2k-2n_{nod(\mathcal{C^\prime})}.
\lb{7.1}\ee
On the other hand, by (\ref{1.1}) we have
\be \dim \overline{\mathfrak{M}}^{vir}_{2, k}(X, A)
=2\langle c_1(TX), A\rangle-2(n-3)+2k.
\lb{7.2}\ee
By (ii) of Theorem 1.2, we have
$0\neq A-A^\prime=mA_2$,
where $[\mathcal{C}_2, u_2]\in{\mathfrak{M}}^0_{g,
0}(X, A_2, J)$. Thus we have
\bea \dim \mathfrak{M}_{T^\prime}(X, A^\prime, J)
\le \dim \overline{\mathfrak{M}}^{vir}_{2,k}(X, A)-2\lb{7.3}\eea
provided $\dim X\equiv 2n<\min\{N+6, 2N+6\}$.
Hence $\mathfrak{M}_T(X, A, J)\cap\overline{\mathfrak{M}}^0_{2, k}(X, A, J)$
serves as a boundary component of
$\mathfrak{M}^0_{2, k}(X, A, J)$
in the sense of pseudocycle, cf. Definition 6.5.1 of \cite{MS}.
\hfill\hb

{ \bf Proof of Theorem 1.4.}
Suppose $u:\Sg\rightarrow X$ factors
through an $m$-fold cover $\Sg^\prime\rightarrow X$, where
$m\ge 2$. Denote the space of equivalence classes of such maps by
$\mathfrak{M}_{\Sg^\prime, m}$. Then we have
\bea ev(\mathfrak{M}_{\Sg^\prime, m})\subset ev\left(\mathfrak{M}^0_{g, k}\left(X , \frac{A}{m}, J\right)\right)
\lb{7.4}\eea
where $g$ is the genus of $\Sg^\prime$.
Note that we have
\bea \dim \mathfrak{M}^0_{g, k}\left(X , \frac{A}{m}, J\right)
\le \dim \overline{\mathfrak{M}}^{vir}_{2,k}(X, A)-2
\lb{7.5}\eea
provided $\dim X\equiv 2n<\min\{N+6, 2N+6\}$ since $g\le 1$. Hence $ev(\mathfrak{M}_{\Sg^\prime, m})$ will not
intersect $\mu_1\times\cdots\times\mu_k$ in general position.
\hfill\hb

\medskip

\bibliographystyle{abbrv}

\medskip

\end{document}